\documentclass{amsproc}

\usepackage{amsmath,amssymb}
\newtheorem{theorem}{Theorem}[section]
\newtheorem{lemm}[theorem]{Lemma}
\newtheorem{prop}[theorem]{Proposition}
\newtheorem{coro}[theorem]{Corollary}

\theoremstyle{definition}
\newtheorem{defi}[theorem]{Definition}

\theoremstyle{remark}

\numberwithin{equation}{section}

\def\de{\delta}

\def\vn{\varepsilon}

\def\ep{\epsilon}
\def\ot{\otimes}
\def\b{\overline}

\def\om{\omega}

\def\dim{\hbox{dim}}

\def\a{\alpha}

\def\b{\beta}

\def\mod{\hbox{mod}}

\newfont{\df}{eufm10}

\def\ep{\epsilon}

\def\ot{\otimes}

\def\de{\delta}
\def\dim{\hbox{\rm dim}\,}

\def\ot{\otimes}

\begin{document}

\title[Restricted two-parameter quantum groups]
{Convex PBW-Type Lyndon Bases and Restricted\\ Two-parameter Quantum
Groups of Type $B$}

\author[Hu]{Naihong Hu$^\star$}
\address{Department of Mathematics, East China Normal University,
Min Hang Campus, Dong Chuan Road 500, Shanghai 200241, PR China}
\email{nhhu@math.ecnu.edu.cn}
\thanks{$^\star$N.H., Corresponding Author,
supported in part by the NNSF (Grant 10728102), the PCSIRT and the
from the MOE, the National \& Shanghai Leading Academic Discipline
Projects (Project Number: B407).}
\thanks{$^*$X.L., supported by the Nankai Research-encouraging Fund for the PhD-Teachers, and a fund from LPMC}

\author[Wang]{Xiuling Wang$^*$}
\address{School of Mathematical Sciences and  LPMC, Nankai University,
Tianjin  300071, PR China}\email{xiulingwang@nankai.edu.cn}

\subjclass{Primary 17B37, 81R50; Secondary 17B35}
\date{June 18, 2005}


\keywords{Convex PBW-type Lyndon basis, restricted 2-parameter
quantum groups, integrals, ribbon Hopf algebra.}
\begin{abstract}
We construct convex PBW-type Lyndon bases for two-parameter quantum
groups $U_{r,s}(\mathfrak{so}_{2n+1})$ with detailed commutation
relations. 
It turns out that under a certain condition, the restricted
two-parameter quantum group $\mathfrak u_{r,s}(\mathfrak
{so}_{2n+1})$ ($r, s$ are roots of unity) is of Drinfeld double. All
of Hopf isomorphisms of $\mathfrak u_{r,s}(\mathfrak {so}_{2n+1})$,
as well as $\mathfrak u_{r,s}(\mathfrak {sl}_n)$ are determined.
Finally, necessary and sufficient conditions for $\mathfrak
u_{r,s}(\mathfrak{so}_{2n+1})$ to be a ribbon Hopf algebra are
singled out by describing the left and right integrals.
\end{abstract}

\maketitle
\section{Introduction}

\noindent {\it 1.1.} In 2001,  from the down-up algebras approach,
Benkart-Witherspoon in \cite{BW1} reobtained Takeuchi's definition
of two-parameter quantum groups of type $A$. Since then, a
systematic study of the two-parameter quantum groups has been going
on, for instance, \cite{BW2, BW3} for type $A$; \cite{BGH1, BGH2}
for type $B, C, D$; \cite{HS, BH} for types $G_2, E_6, E_7, E_8$.
For a
unified definition, see \cite{HP}. 
Hu-Rosso-Zhang \cite{HRZ} defined the affine type $A$ case, obtained
its Drinfeld realization in two-parameter quantum version and
constructed the quantum affine Lyndon basis. Furthermore, the vertex
operator representations of level $1$ for the two-parameter quantum
affine algebras $U_{r,s}(\widehat {\mathfrak g})$ can be established
(cf. \cite{HZ}, etc.). While, in the root of unity case, Hu-Wang
\cite{HW} generalized the work \cite{BW3} to the restricted type
$G_2$ case. The goal of this paper is to solve the restricted type
$B$ case with overcoming some
peculiar difficulties itself in this case. 

\medskip\noindent
{\it 1.2.}  In the study of Lusztig's symmetry property for the
two-parameter quantum groups in question, 
Theorem 3.1 \cite{BGH1} says that the Lusztig's symmetries only
exist as $\mathbb Q$-isomorphisms between $U_{r,s}(\mathfrak g)$ and
its associated object $U_{s^{-1},r^{-1}}(\mathfrak g)$ when
$\textrm{rank}(\mathfrak g)=2$, and in the case when
$\textrm{rank}(\mathfrak g)>2$, the ``{\it iff} " condition for the
existence of Lusztig's symmetries forces $U_{r,s}(\mathfrak g)$ to
take the ``one-parameter" form $U_{q,q^{-1}}(\mathfrak g)$, where
$r=q=s^{-1}$. This means that one cannot conveniently write out the
convex PBW-type basis for $U_{r,s}(\mathfrak g)$ like in the
one-parameter cases using Lusztig's braid group actions in a typical
fashion (see \cite{Lu}). This is one of difficulties encountered in
the two-parameter setting. Note that the combinatorial construction
of the PBW-type bases in the quantum cases is available but
nontrivial, as it depends not only on the choice of a convex
ordering (\cite{B, Ro}) on a positive root system, but also the
inserting manner of the $\bold q$-bracketings with suitable
structure constants (see \cite{Ro, HRZ}, etc.) into the good Lyndon
words (cf. \cite{LR, Ro}). The constructions of convex PBW-type
bases in the two-parameter cases have been given for type $A$ in
\cite{BW3}, types $E_6, E_7, E_8$ in \cite{BH}, for types $B_2$,
$C_2$, $D_4$ in \cite{H}, and for type $G_2$ in \cite{HW}. Indeed,
in order to study the properties of restriced quantum groups as
finite-dimensional Hopf algebras, it is necessary to find the nice
bases for the corresponding nonrestricted quantum groups at generic
cases. 
The first thing of the article is to present a direct construction
for a convex PBW-type Lyndon basis of type $B$ (for arbitrary rank)
and provide explicit information on commutation relations among
basis elements that is decisive to single out the left (or right)
integrals of $\mathfrak{u}_{r,s}(\mathfrak{so}_{2n+1})$. Note that
the complexity of Lyndon bases corresponding to the nonsimply-laced
Dynkin diagrams causes more inconvenience when treated with the type
$B$ case here.

\medskip
\noindent{\it 1.3.} The paper is organized as follows. Section 2 is
to give the construction of Lyndon bases
for $U_{r,s}(\mathfrak{so}_{2n+1})$. 
In Section 3, we contribute more efforts to make the possible
commutation relations clearly. These give rise to central elements
of degree $\ell$ in the case when parameters $r$, $s$ are roots of
unity, which generate a Hopf ideal. The quotient is the restricted
two-parameter quantum groups $\mathfrak
u_{r,s}(\mathfrak{so}_{2n+1})$ in Section 4. In Section 5, we show
that these Hopf algebras are pointed, and determine all Hopf
isomorphisms of $\mathfrak u_{r,s}(\mathfrak {so}_{2n+1})$, as well
as $\mathfrak u_{r,s}(\mathfrak {sl}_n)$ in terms of the set of
(left) right skew-primitive elements. Section 6 is to prove that
$\mathfrak u_{r,s}(\mathfrak{so}_{2n+1})$ is of Drinfel'd double
under a certain condition. In Section 7, the left and right
integrals of $\mathfrak{b}$ are singled out based on the discussion
in Section 3, which are applied to determine the necessary and
sufficient conditions for $\mathfrak u_{r,s}(\mathfrak{so}_{2n+1})$
to be ribbon in Section 8.

\section{$U_{r,s}(\mathfrak{so}_{2n+1})$ and its convex PBW-type Lyndon basis}

\medskip
\noindent{\it 2.1.} Let ${\Bbb K}={\Bbb Q}(r,s)$ 
be a subfield of $\mathbb C$, where $r, s\in\mathbb C^*$ with
assumptions $r^3\neq s^3$ and $r^4\neq s^4$. Let $\Phi$ be a root
system of $\mathfrak{so}_{2n+1}$ with $\Pi$ a base of simple roots,
which is a finite subset of a Euclidean space $E = {\Bbb R}^3$ with
an inner product $(\,,\,)$. Let $\epsilon _{1},\cdots,\epsilon_{n}$
denote an orthonormal basis of $E$, then $\Pi = \{\alpha_{i} =
\epsilon_{i}-\epsilon_{i+1} \mid 1\leq i< n\}\cup \{\alpha_{n} =
\epsilon_{n}\}$, $\Phi =\{ \pm \epsilon_{i}\pm \epsilon_{j} \mid
1\leq i\neq j\leq n\}\cup \{ \pm \epsilon_{i} \mid 1\leq i \leq
n\}$. In this case, set $\displaystyle r_i =
r^{(\alpha_i,\,\alpha_i)}$, $s_i = s^{(\alpha_i,\,\alpha_i)}$, so
that $r_1= \cdots =r_{n-1}=r^2,r_n=r$ and $s_1= \cdots
=s_{n-1}=s^2,s_n=s$.

Given two sets of symbols $W=\{\omega_1,\cdots,\omega_n\}$,
$W'=\{\omega_1',\cdots,\omega_n'\}$. Define the structural constants
matrix $(\langle \omega_i',\,\omega_j \rangle)_{n\times n}$ of type
$B$ by
\[ \begin{pmatrix} r^2s^{-2} & r^{-2}&1  &\cdots &1&1\\
s^2& r^2s^{-2}& r^{-2} & \cdots &1&1 \\
1&s^2& r^2s^{-2}&  \cdots&1&1\\
\cdots&\cdots&\cdots&\cdots& \cdots&\cdots\\
1&1&1&\cdots& r^2s^{-2}& r^{-2}\\
1&1&1&\cdots&s^2& rs^{-1}
\end{pmatrix}.\]

\begin{defi} (\cite{BGH1})
Let $U=U_{r,\,s}(\mathfrak{so}_{2n+1})$ be the associative algebra
over $\mathbb{K}$ generated by  $e_i,\;f_i,\;\omega_i^{\pm 1},\;
\omega_i'^{\pm 1} \;(1\leq i \leq n)$ subject to relations
$(B1)-(B5)$:
\begin{gather*}
\text{The $\omega_i^{\pm 1}, \omega_j'^{\pm 1}$ all commute with one
another and  $\omega_i\omega_i^{-1} =1= \omega_j'\omega_j'^{-1}$}. \tag{\text{$B1$}}\\
\omega_{j}\,e_{i}\,\omega_{j}^{-1} =\langle \om_i',\om_j\rangle e_i,
\qquad \omega_{j}\,f_{i}\,\omega_{j}^{-1} =\langle
\om_i',\om_j\rangle^{-1} f_i. \tag{\text{$B2$}}\\
\omega_{j}'\,e_{i}\,\omega_{j}'^{-1} =\langle
\om_j',\om_i\rangle^{-1} e_i, \qquad
\omega_{j}'\,f_{i}\,\omega_{j}'^{-1} =\langle \om_j',\om_i\rangle
f_i.\tag{\text{$B3$}}\\
[\,e_i,
f_j\,]=\delta_{ij}\frac{\om_i-\om_i'}{r_i-s_i}.\tag{\text{$B4$}}
\end{gather*}
\noindent $(B5)$ \quad $(r,s)$-Serre relations
$$
(\text{ad}_\ell e_i)^{1-a_{ij}}(e_j)=0,\qquad (\text{ad}_r
f_i)^{1-a_{ij}}(f_j)=0
$$
are given by its Hopf algebra structure on
 $U_{r,s}(\mathfrak{so}_{2n+1})$ with
the comultiplication, the counit and the antipode as follows
\begin{gather*}
\Delta(\om_i^{\pm1})=\om_i^{\pm1}\ot\om_i^{\pm1}, \qquad
\Delta({\om_i'}^{\pm1})={\om_i'}^{\pm1}\ot{\om_i'}^{\pm1},\\
\Delta(e_i)=e_i\ot 1+\om_i\ot e_i, \qquad \Delta(f_i)=1\ot
f_i+f_i\ot \om_i',\\
\vn(\om_i^{\pm})=\vn({\om_i'}^{\pm1})=1, \qquad
\vn(e_i)=\vn(f_i)=0,\\
S(\om_i^{\pm1})=\om_i^{\mp1}, \qquad
S({\om_i'}^{\pm1})={\om_i'}^{\mp1},\\
S(e_i)=-\om_i^{-1}e_i,\qquad S(f_i)=-f_i\,{\om_i'}^{-1}.
\end{gather*}
\end{defi}

When $r=q=s^{-1}$, $U_{r,\,s}(\mathfrak{so}_{2n+1})/I\cong
U_q(\mathfrak{so}_{2n+1})$, the standard Drinfel'd-Jimbo quantum
group, where $I=(\omega_i' - \omega_i^{-1}, \forall\;i)$ is a Hopf
ideal.

\smallskip
\noindent {\it 2.2.} $U$ has a triangular decomposition
$U\cong U^-\otimes U^0\otimes U^+$, where $U^0$ is the subalgebra
generated by $\omega_i^{\pm 1}, {\omega'_i}^{\pm 1}$, and $U^+$
(resp. $U^-$) is the subalgebra generated by $e_i$
 (resp. $f_i$). Let $\mathcal{B}$ (resp. $\mathcal{B}'$)
denote the Hopf subalgebra of $U$ generated by $e_j, \omega_j^{\pm
1}$ (resp. $f_j, {\omega'_j}^{\pm 1}$) with $1\leq j\leq n$.

\begin{prop} $($\cite{BGH1}$)$
There exists a unique skew-dual pairing $\langle\,,\rangle:$ $
\mathcal{B}'\times \mathcal{B}\rightarrow \mathbb{K}$ of the Hopf
subalgebras $\mathcal{B}$ and $\mathcal{B}'$ such that
\begin{gather*}
\langle f_i,\,e_j \rangle = \delta_{ij} \frac{1}{s_i - r_i},
\qquad 1 \leq i,\,j \leq n,\\
\langle {\omega'}^{\pm 1}_i, \omega^{-1}_j\rangle=\langle
{\omega'}^{\pm 1}_i, \omega_j\rangle^{-1}=
 \langle {\omega'}_i, \omega_j\rangle^{\mp 1}, \quad 1 \leq i,\,j \leq n
\end{gather*}
and all other pairs of generators are $0$. Moreover, we have
$\langle S(a),\,S(b) \rangle = \langle a,\,b \rangle$ for $a\in
\mathcal{B}', b\in \mathcal{B}$.
\end{prop}

\noindent {\it 2.3.} Write the positive root system
$\Phi^+=\{\alpha_{ij},\, \epsilon_{\ell},\, \beta_{ij}\mid 1\le i<
j\le n, 1\le\ell\le n \}$, where
$\alpha_{ij}:=\alpha_i+\alpha_{i+1}+\cdots+\alpha_{j-1}
=\epsilon_i-\epsilon_j, \,
~\epsilon_{\ell}=\alpha_\ell+\alpha_{\ell+1}+\cdots+\alpha_{n}, \,
~\beta_{ij}:=\alpha_i+\alpha_{i+1}+\cdots+2\alpha_{n}+\alpha_{n-1}+\cdots+
\alpha_j=\epsilon_i+\epsilon_j$. To a reduced expression of the
maximal length element of Weyl group $W$ for type $B_n$ (see
\cite{CX}) taken as
$$w_0=(s_1s_2\cdots s_{n-1}s_ns_{n-1}\cdots s_2s_1)
\cdots (s_{n-1}s_ns_{n-1}) (s_{n}),$$ there corresponds a convex
ordering on $\Phi^+$ below:
\begin{equation*}
\begin{split}
\alpha_{12}, \quad  \alpha_{13}, \quad \cdots,\quad \alpha_{1n},
\quad &\epsilon_1,\quad \beta_{1n} ,\quad \cdots,\quad
 \beta_{13}, \quad \beta_{12}, \\
\alpha_{23}, \quad \cdots,\quad \alpha_{2n}, \quad &\epsilon_2,\quad
\beta_{2n}, \quad \cdots, \quad
 \beta_{23},  \\
&\cdots \\
\alpha_{n{-}1,n},\quad  &\epsilon_{n{-}1},\quad \beta_{n{-}1,n}, \\
&\epsilon_n.
\end{split}
\end{equation*}
{\it Remark:} The above ordering also corresponds to the standard
Lyndon tree of type $B_n$ (\cite{LR}):

\unitlength 0.78mm 
\linethickness{0.1pt}
\ifx\plotpoint\undefined\newsavebox{\plotpoint}\fi 
\begin{picture}(50,10)(0,10)
\put(9.25,14.5){$B_n:$} \put(29.5,14.5){\circle{2}}
\put(39.5,14.5){\circle{2}} \put(59.5,14.5){\circle{2}}
\put(69.5,14.5){\circle{2}} \put(79.5,14.5){\circle{2}}
\put(89.5,14.5){\circle{2}} \put(109.5,14.5){\circle{2}}
\put(119.5,14.5){\circle{2}}
\put(30.5,14.5){\line(1,0){8}}
\put(40.5,14.5){\line(1,0){4}}
\put(54.5,14.5){\line(1,0){4}}
\put(60.5,14.5){\line(1,0){8}}
\put(70.5,14.5){\line(1,0){8}}
\put(80.5,14.5){\line(1,0){8}}
\put(90.5,14.5){\line(1,0){4}}
\put(104.5,14.5){\line(1,0){4}}
\put(110.5,14.5){\line(1,0){8}}
\multiput(44.5,14.5)(.975,0){11}{{\rule{.4pt}{.4pt}}}
\multiput(94.5,14.5)(.975,0){11}{{\rule{.4pt}{.4pt}}}
\put(29.75,8.5){i} \put(39.25,8.5){i{+}1} \put(59.25,8.5){n{-}1}
\put(69.25,8.5){n} \put(79.25,8.5){n} \put(89.25,8.5){n{-}1}
\put(109.25,8.5){i{+}2} \put(119.25,8.5){i{+}1}
\end{picture}\\

Denote by $\mathcal E_\alpha$ the quantum root vector in
$U^+$ 
for each $\alpha\in\Phi^+$.
Briefly, we set $\mathcal
E_{i,j}:=\mathcal E_{\alpha_{i,j+1}}$ for $1\le i\le j< n$;
$\mathcal E_{i,n}:=\mathcal E_{\epsilon_i}$ for $1\le i\le n$;
$\mathcal E_{i,j'}:=\mathcal E_{\beta_{i,j}}$ for $1\le i<j\le n$;
in particular, $\mathcal E_{i,i}:=\mathcal E_{\alpha_i}=e_i$ for
$1\le i\le n$, where in our notation $\alpha_i=\alpha_{i,i+1}$
($i<n$) and $\alpha_n=\epsilon_n$.

\medskip
\noindent{\it 2.4.} With the above ordering and notation, following
the grammatical rule of good Lyndon words with inserting
$(r,s)$-bracketings, we can make an inductive definition starting
from the bottom to up of the above root-ordering graph as follows.

\begin{defi} For $\alpha\in \Phi^+$, define $\mathcal E_{\alpha}$ inductively
\begin{eqnarray}
\mathcal{E}_{i,i}=e_i, \qquad\qquad\qquad 1\leq i\leq n,\label{2.1} \\
\mathcal{E}_{i,j}=e_i\mathcal{E}_{i+1,j}-r^2\mathcal{E}_{i+1,j}e_i,
\qquad 1\leq i< j \leq n, \label{2.2}\\
\mathcal{E}_{i,n'}=\mathcal{E}_{i,n}e_n-rse_n\mathcal{E}_{i,n},
\qquad1\leq i \leq n-1, \label{2.3}\\
\mathcal{E}_{i,j'}=\mathcal{E}_{i,j+1'}e_j-s^{-2}e_j\mathcal{E}_{i,j+1'},
\qquad 1\leq i <j \leq n-1.\label{2.4}
\end{eqnarray}
\end{defi}
Then the $(r,s)$-Serre relations for $U^{+}$ in  $(B5)$ can be
reformulated as
\begin{eqnarray}
e_i\mathcal{E}_{i,i+1}=s^2\mathcal{E}_{i,i+1}e_i, \label{2.5}\\
\mathcal{E}_{j,j+1}e_{j+1}=s^2e_{j+1}\mathcal{E}_{j,j+1},\label{2.6}\\
\mathcal{E}_{n-1,n'}e_n=s^2e_n\mathcal{E}_{n-1,n'},\label{2.7}
\end{eqnarray}
for $1\leq i \leq n-1$, $1\le j<n-1$.

\medskip
\noindent{\it Remark:} (i) From \cite{Ro}, the set of quantum Lie
brackets (or say, $\bold q$-bracketings) of all good Lyndon words
consists of all quantum root vectors of $U^+$ (for definitions to
see {\it Remark} 3.7 (3) in \cite{HRZ}, pp. 467).

\smallskip
(ii) As for how to insert $(r,s)$-bracketings into each good Lyndon
word in type $B$ case, in \cite{H} an observation was obtained via a
representation theoretic analysis in rank $2$ case, that is, we have
to obey the defining rule as $\mathcal E_\gamma:=[\mathcal
E_\alpha,\mathcal E_\beta]_{\langle
\omega_{\alpha}',\omega_{\beta}\rangle^{-1}}=\mathcal
E_\alpha\mathcal E_\beta-\langle
\omega_{\alpha}',\omega_{\beta}\rangle^{-1}\mathcal E_\beta\mathcal
E_\alpha$ for $\alpha, \gamma, \beta\in \Phi^+$ such that
$\alpha<\gamma<\beta$ in the convex ordering, and
$\gamma=\alpha+\beta$.

\medskip
\noindent{\it 2.5.} 
According to \cite{H, K, Ro, LR}, we have

\begin{theorem}
$\Bigl\{\,\mathcal
{E}_{n,n}^{c_{n^2}}\cdot\,\mathcal{E}_{n{-}1,n'}^{c_{n^2{-}1}}
\mathcal{E}_{n{-}1,n}^{c_{n^2{-}2}}\mathcal{E}_{n{-}1,n{-}1}^{c_{n^2{-}3}}\cdots
\mathcal{E}_{2,3'}^{c_{4n{-}4}}\cdots
\mathcal{E}_{2,n'}^{c_{3n{-}1}}\mathcal{E}_{2,n}^{c_{3n{-}2}} \cdots
\mathcal{E}_{2,2}^{c_{2n}}$

\noindent $\left.\mathcal{E}_{1,2'}^{c_{2n{-}1}}\cdots
\mathcal{E}_{1,n'}^{c_{n{+}1}}\mathcal{E}_{1,n}^{c_n}
\cdots\mathcal{E}_{1,2}^{c_2}\mathcal{E}_{1,1}^{c_1}\; \right|\;
(c_1,\cdots,c_{n^2})\in \Bbb Z_+^{n^2}\Bigr\}$ forms a convex
PBW-type Lyndon basis of the algebra $U^+$.\hfill\qed
\end{theorem}

\noindent{\it Remark:} (1) One can take a lexicographical ordering
on the set of positive quantum root vectors provided that one puts
$j:=2n-j'+1$ for $1\le j'\le n$. That is, the set
$$\left.\left\{\;\mathcal{E}_{i_p,j_p}\cdots\mathcal{E}_{i_2,j_2}\mathcal{E}_{i_1,j_1}\;\right|\;
(i_1,j_1)\leq (i_2,j_2)\leq \cdots \leq (i_p,j_p) \quad
\textit{lexicographically}\;\right \},$$ where $1\le i_l\le n$ and
$1\le j_l\le 2n$, is a basis of the algebra $U^+$.

\smallskip
(2) The set of elements
$\mathcal{E}_{1,1}^{c_1}\mathcal{E}_{1,2}^{c_2}\mathcal{E}_{1,3}^{c_3}\cdots\mathcal{E}_{1,n}^{c_n}
\mathcal{E}_{1,n'}^{c_{n+1}}\cdots
\mathcal{E}_{1,2'}^{c_{2n-1}}\mathcal{E}_{2,2}^{c_{2n}}\mathcal{E}_{2,3}^{c_{2n+1}}$
$\cdots \mathcal{E}_{2,3'}^{c_{4n-4}}$
$\mathcal{E}_{3,3}^{c_{4n-3}}\cdots\mathcal{E}_{n-1,n-1}^{c_{n^2-3}}
 \mathcal{E}_{n-1,n}^{c_{n^2-2}}\mathcal{E}_{n-1,n'}^{c_{n^2-1}}\mathcal{E}_{n,n}^{c_{n^2}}
 \ \Bigl( (c_1,\cdots,c_{n^2})\in N^{n^2}\Bigr)$ forms a
basis of $U^+$.

\medskip
\noindent{\it 2.6.} The following $\tau$ plays a key role in describing Lyndon basis (see \cite{HRZ}). 

\begin{defi}
Let $\tau$ be the $\mathbb{Q}$-algebra anti-automorphism of
$U_{r,s}(\mathfrak {so}_{2n+1})$ such that $\tau(r)=s$, $\tau(s)=r$,
$\tau(\langle \om_i',\om_j\rangle^{\pm1})=\langle
\om_j',\om_i\rangle^{\mp1}$, and
\begin{gather*}
\tau(e_i)=f_i, \quad \tau(f_i)=e_i, \quad \tau(\om_i)=\om_i',\quad
\tau(\om_i')=\om_i.
\end{gather*}
Then $\mathcal B'=\tau(\mathcal B)$ with those induced defining
 relations from $\mathcal B$, and those cross relations in
$(\textrm{B2})$---$(\textrm{B4})$ are antisymmetric with respect to
$\tau$.\hfill\qed
\end{defi}

Using $\tau$ to $U^+$, we can get those negative quantum root
vectors in $U^-$. Define $\mathcal{F}_{i,i}=\tau(\mathcal
E_{i,i})=f_i$ for $1\leq i\leq n$, and
\begin{eqnarray}
\mathcal{F}_{i,j}=\tau(\mathcal
E_{i,j})=\mathcal{F}_{i+1,j}f_i-s^2f_i\mathcal{F}_{i+1,j},
\qquad 1\leq i< j \leq n, \\
\mathcal{F}_{i,n'}=\tau(\mathcal
E_{i,n'})=f_n\mathcal{F}_{i,n}-rs\mathcal{F}_{i,n}f_n,
\qquad 1\leq i \leq n-1, \\
\mathcal{F}_{i,j'}=\tau(\mathcal
E_{i,j'})=f_j\mathcal{F}_{i,j+1'}-r^{-2}\mathcal{F}_{i,j+1'}f_j,
\quad 1\leq i <j \leq n-1.
\end{eqnarray}

\begin{coro}
$\Bigl\{\,\mathcal{F}_{1,1}^{c_1}\mathcal{F}_{1,2}^{c_2}
\mathcal{F}_{1,3}^{c_3}\cdots\mathcal{F}_{1,n}^{c_n}
\mathcal{F}_{1,n'}^{c_{n+1}}\cdots
\mathcal{F}_{1,2'}^{c_{2n-1}}\cdot\,\mathcal{F}_{2,2}^{c_{2n}}
\mathcal{F}_{2,3}^{c_{2n+1}}\cdots\mathcal{F}_{2,n}^{c_{3n{-}2}} $

\noindent
$\left.\mathcal{F}_{2,n'}^{c_{3n{-}1}}\cdots\mathcal{F}_{2,3'}^{c_{4n-4}}\cdot\cdots\mathcal{F}_{n-1,n-1}^{c_{n^2-3}}
 \mathcal{F}_{n-1,n}^{c_{n^2-2}}\mathcal{F}_{n-1,n'}^{c_{n^2-1}}\cdot\mathcal{F}_{n,n}^{c_{n^2}}
\;\right| \; (c_1,\cdots,c_{n^2}) \in \Bbb Z_+^{n^2}\Bigr\}$  forms
a convex PBW-type basis of the algebra $U^-$.\hfill\qed
\end{coro}


\section{Commutation relations and homogeneous central elements}

\medskip
\noindent {\it 3.1. Commutation relations in $U^{+}$.} We will make
more efforts to explicitly calculate all possible commutation
relations, which are useful for determining central elements in
subsection {\it 3.2} and integrals in Section 7. We point out the
fact that these commutation relations or say, $(r,s)$-identities
hold in $U^+$ arises essentially from the $(r,s)$-Serre relations
$(B5)$. First of all, we pay attention to the following

\medskip {\bf Question:} Is the defining result of $(r,s)$-bracketings
in {\it Remark} (ii) unique for those non-unique decomposition of
$\gamma=\alpha'+\beta'$ such that $\alpha', \beta'\in\Phi^+$ and
$\alpha'<\gamma<\beta'$? i.e., does it have $[\mathcal
E_{\alpha'},\mathcal E_{\beta'}]_{\langle
\omega_{\alpha'}',\omega_{\beta'}\rangle^{-1}}=[\mathcal
E_{\alpha},\mathcal E_{\beta}]_{\langle
\omega_{\alpha}',\omega_{\beta}\rangle^{-1}}$ for such a $\gamma$?

\smallskip
Lemma 3.1 (3) \& (4) below have a positive answer to this question
provided that we work under the structure constants matrix listed in
subsection {\it 2.1}. We will see that Lemma 3.1 (3) \& (4) are not
only the most basic ones in all commutation relations involved but
also well-adopted to our calculating technique.

\medskip For simplicity, we write $(r,s)$-bracketings by
$[\mathcal E_{\alpha},\mathcal E_{\beta}]_{\bullet}=[\mathcal
E_{\alpha},\mathcal E_{\beta}]_{\langle
\omega_{\alpha}',\omega_{\beta}\rangle^{-1}}$ briefly.

\begin{lemm} The following relations hold in $U^{+}:$

\smallskip
$(1)\quad
\mathcal{E}_{i,j}\mathcal{E}_{k,l}=\mathcal{E}_{k,l}\mathcal{E}_{i,j},
\hskip2.73cm\quad\ \, i\leq j,\ j+1<k\leq l\le n;$

$(2)\quad
\mathcal{E}_{i,j}\mathcal{E}_{k,l'}=\mathcal{E}_{k,l'}\mathcal{E}_{i,j},
\hskip2.73cm\quad i\leq j,\ j+1<k< l\le n;$

$(3)\quad
\mathcal{E}_{i,j}=\mathcal{E}_{i,l-1}\mathcal{E}_{l,j}-r^2\mathcal{E}_{l,j}
\mathcal{E}_{i,l-1},\ \;\; \quad\hskip0.73cm i< l\le j\leq n;$

$(4)\quad
\mathcal{E}_{i,j'}=\mathcal{E}_{i,l-1}\mathcal{E}_{l,j'}-r^2\mathcal{E}_{l,j'}
\mathcal{E}_{i,l-1},\quad\,\hskip0.72cm i<l<j\leq n.$
\end{lemm}
\begin{proof}
(1) \& (2)  follow directly from the $(r,s)$-Serre relation $(B5)$.

(3): We fix $l$ and $j$ with $l\le j$ and use induction on $i$. If
$i=l-1$, this is just the defining formula (2.2) of
$\mathcal{E}_{i,j}$. Assume that $i< l-1$, then we have
\begin{equation*}
\begin{split}
\mathcal{E}_{i,l-1}\mathcal{E}_{l,j}  & =
(e_i\mathcal{E}_{i+1,l-1}-r^2\mathcal{E}_{i+1,l-1}e_i)\mathcal{E}_{l,j}\\
& = e_i\mathcal{E}_{i+1,l-1}\mathcal{E}_{l,j}
-r^2\mathcal{E}_{i+1,l-1}\mathcal{E}_{l,j}e_i\\
& =
e_i(r^2\mathcal{E}_{l,j}\mathcal{E}_{i+1,l-1}+\mathcal{E}_{i+1,j})
-r^2(r^2\mathcal{E}_{l,j}\mathcal{E}_{i+1,l-1}+\mathcal{E}_{i+1,j})e_i\\
&=r^2\mathcal{E}_{l,j}e_i\mathcal{E}_{i+1,l-1}+e_i\mathcal{E}_{i+1,j}-
r^4\mathcal{E}_{l,j}\mathcal{E}_{i+1,l-1}e_i-r^2\mathcal{E}_{i+1,j}e_i
\\ & =r^2\mathcal{E}_{l,j}\mathcal{E}_{i,l-1}+\mathcal{E}_{i,j}
\end{split}
\end{equation*}
by (1), the induction hypothesis and the defining formula (2.2).

(4): The following expression can be easily verified:
$$
\mathcal{E}_{i,n'}=e_i\mathcal{E}_{i+1,n'}-r^2\mathcal{E}_{i+1,n'}e_i,
\quad i\leq n-2.  \eqno (*)  $$ Indeed, for $i\leq n-2$, we have
$e_ie_n=e_ne_i$. By the defining  relation (\ref{2.3}), we get
\begin{equation*}
\begin{split}
e_{i}\mathcal{E}_{i+1,n'}  & =
e_i(\mathcal{E}_{i+1,n}e_n-rse_n\mathcal{E}_{i+1,n})\\
& =
(\mathcal{E}_{i,n}+r^2\mathcal{E}_{i+1,n}e_i)e_n-rse_n(\mathcal{E}_{i,n}+r^2\mathcal{E}_{i+1,n}e_i)\\
& =\mathcal{E}_{i,n'}+r^2\mathcal{E}_{i+1,n'}e_i.
\end{split}
\end{equation*}

We suppose initially  that $j=n$ and $i=l-1$, this is just the
relation in $(*)$.

Now assume (4) is true for $\mathcal E_{i+1,n'}$. This means
$i{+}1<l$. Then using (2), $e_i\mathcal E_{l,n'}=\mathcal
E_{l,n'}e_i$. By the induction hypothesis and ($*$), we obtain
\begin{equation*}
\begin{split}
\mathcal{E}_{i,l-1}\mathcal{E}_{l,n'}  & =
(e_i\mathcal{E}_{i+1,l-1}-r^2\mathcal{E}_{i+1,l-1}e_i)\mathcal{E}_{l,n'}\\
& =
e_i(\mathcal{E}_{i+1,n'}+r^2\mathcal{E}_{l,n'}\mathcal{E}_{i+1,l-1})
-r^2\mathcal{E}_{i+1,l-1}\mathcal{E}_{l,n'}e_i\\
& =
e_i\mathcal{E}_{i+1,n'}+r^2\mathcal{E}_{l,n'}e_i\mathcal{E}_{i+1,l-1}-r^2
\mathcal{E}_{i+1,n'}e_i-r^4\mathcal{E}_{l,n'}\mathcal{E}_{i+1,l-1}e_i\\
&=\mathcal{E}_{i,n'}+r^2\mathcal{E}_{l,n'}\mathcal{E}_{i,l-1}.
\end{split}
\end{equation*}

Finally, assume (4) is true for $\mathcal E_{i,j+1'}$ with $i<l$.
Since $i\le l{-}1<l<j$,
$e_j\mathcal{E}_{i,l-1}=\mathcal{E}_{i,l-1}e_j$ by (1). Using the
defining formula (2.4) and the induction hypothesis, we get
\begin{equation*}
\begin{split}
\mathcal{E}_{i,j'} & =
\mathcal{E}_{i,j+1'}e_j{-}s^{-2}e_j\mathcal{E}_{i,j+1'}\\
& = (\mathcal{E}_{i,l-1}\mathcal{E}_{l,j+1'}
{-}r^2\mathcal{E}_{l,j+1'}\mathcal{E}_{i,l-1})e_j{-}s^{-2}
e_j(\mathcal{E}_{i,l-1}\mathcal{E}_{l,j+1'}
{-}r^2\mathcal{E}_{l,j+1'}\mathcal{E}_{i,l-1})\\
&
=\mathcal{E}_{i,l{-}1}\mathcal{E}_{l,j{+}1'}e_j{-}r^2\mathcal{E}_{l,j{+}1'}e_j
\mathcal{E}_{i,l{-}1}{-}s^{{-}2}
\mathcal{E}_{i,l{-}1}e_j\mathcal{E}_{l,j{+}1'}
{+}s^{{-}2}r^2e_j\mathcal{E}_{l,j{+}1'}\mathcal{E}_{i,l{-}1}\\
&=\mathcal{E}_{i,l{-}1}\mathcal{E}_{l,j'}-r^2\mathcal{E}_{l,j'}\mathcal{E}_{i,l{-}1}.
\end{split}
\end{equation*}

Thus, the proof is complete.
\end{proof}

\begin{lemm} The following relations hold in $U^{+}:$

\smallskip
$(1)\quad e_i\mathcal E_{i,j}=s^2\mathcal E_{i,j}e_i,
\;\,\qquad\qquad 1\le i<j\le n;$

$(2)\quad e_i\mathcal{E}_{i,j'}=s^2\mathcal{E}_{i,j'}e_i,\qquad
\qquad i+1<j\leq n;$

$(3)\quad\mathcal E_{i,j}e_j=s^2e_j\mathcal E_{i,j}, \qquad\;\qquad
1\le i<j< n;$

$(4)\quad\mathcal{E}_{i,n'}e_n=s^2e_n\mathcal{E}_{i,n'},
\qquad\quad\; 1\le i< n;$

$(5)\quad \mathcal E_{i,j'}e_j=r^{-2}e_j\mathcal
E_{i,j'},\qquad\quad\hskip0.05cm 1\le i<j<n.$
\end{lemm}
\begin{proof}
(1) follows directly from the $(r,s)$-Serre relations $(B5)$, (2.5)
and Lemma 3.1 (1) \& (3) as one can decompose $\mathcal
E_{i,j}=[\mathcal E_{i,i+1},\mathcal E_{i+2,j}]_{r^2}$.

(2), (3) \& (4) can be proved similarly using the $(r,s)$-Serre
relations $(B5)$, (2.6) or (2.7), together with Lemma 3.1.

(5) will be proved in Lemma 3.6 (2).
\end{proof}

\begin{lemm} The following relations hold in $U^{+}:$

\smallskip
$(1)\quad e_l\mathcal E_{i,j}=\mathcal E_{i,j}e_l,
\hskip3.31cm\qquad\qquad\ i<l<j\le n;$

$(2)\quad
\mathcal{E}_{i,l}\mathcal{E}_{l,j}-(rs)^2\mathcal{E}_{l,j}\mathcal{E}_{i,l}
=(s^2-r^{2})e_l\mathcal{E}_{i,j},\qquad\hskip0.41cm i<l<j\le n;$

$(3)\quad
\mathcal{E}_{i,j}\mathcal{E}_{k,l}=\mathcal{E}_{k,l}\mathcal{E}_{i,j},\quad\hskip3.25cm\qquad
i<k\le l<j\le n;$

$(4)\quad
\mathcal{E}_{i,j}\mathcal{E}_{l,j}=s^2\mathcal{E}_{l,j}\mathcal{E}_{i,j},\qquad
\hskip2.98cm\quad i<l\le j<n;$

$(5)\quad
\mathcal{E}_{i,l}\mathcal{E}_{i,j}=s^2\mathcal{E}_{i,j}\mathcal{E}_{i,l},\qquad
\hskip2.98cm\quad\, i\le l<j\le n;$

$(6)\quad e_l\mathcal E_{i,j'}=\mathcal E_{i,j'}e_l,
\hskip3.3cm\qquad\quad\,\ \;i<l, l+1<j\le n;$

$(7)\quad
\mathcal{E}_{i,l}\mathcal{E}_{l,j'}-(rs)^2\mathcal{E}_{l,j'}\mathcal{E}_{i,l}
=(s^2-r^{2})e_l\mathcal{E}_{i,j'},\hskip0.13cm\qquad i<l, l+1<j\le
n;$

$(8)\quad
\mathcal{E}_{i,j'}\mathcal{E}_{k,l}=\mathcal{E}_{k,l}\mathcal{E}_{i,j'},\hskip3.42cm\qquad
i<k\le l, l+1<j\le n;$

$(9)\quad
\mathcal{E}_{i,l}\mathcal{E}_{i,j'}=s^2\mathcal{E}_{i,j'}\mathcal{E}_{i,l},\hskip3.2cm\qquad
i\le l, l+1<j\le n.$
\end{lemm}
\begin{proof} \
(1) \& (2): Using Lemma 3.1 (3) to $\mathcal E_{l,j}$ and Lemma 3.2
(3), we get
\begin{equation*}
\begin{split}
\mathcal E_{i,l}&\mathcal E_{l,j}-(rs)^2\mathcal
E_{l,j}\mathcal E_{i,l}\\
&=\mathcal E_{i,l}(e_l\mathcal E_{l+1,j}-r^2\mathcal
E_{l+1,j}e_l)-(rs)^2(e_l\mathcal E_{l+1,j}-r^2\mathcal
E_{l+1,j}e_l)\mathcal
E_{i,l}\\
&=s^2e_l\mathcal E_{i,l}\mathcal E_{l+1,j}-r^2\mathcal
E_{i,l}\mathcal E_{l+1,j}e_l-(rs)^2e_l\mathcal E_{l+1,j}\mathcal
E_{i,l}+r^4\mathcal E_{l+1,j}\mathcal
E_{i,l}e_l\\
&=s^2e_l\mathcal E_{i,j}-r^2\mathcal E_{i,j}e_l.
\end{split}
\end{equation*}

Again, using Lemma 3.1 (3) to $\mathcal E_{i,l}$ and Lemma 3.2 (1),
we have
\begin{equation*}
\begin{split}
\mathcal E_{i,l}&\mathcal E_{l,j}-(rs)^2\mathcal
E_{l,j}\mathcal E_{i,l}\\
&=(\mathcal E_{i,l-1}e_l-r^2e_l\mathcal E_{i,l-1})\mathcal
E_{l,j}-(rs)^2\mathcal E_{l,j}(\mathcal
E_{i,l-1}e_l-r^2e_l\mathcal E_{i,l-1})\\
&=s^2\mathcal E_{i,l-1}\mathcal E_{l,j}e_l-r^2e_l\mathcal
E_{i,l-1}\mathcal E_{l,j} -(rs)^2\mathcal E_{l,j}\mathcal
E_{i,l-1}e_l+r^4e_l\mathcal E_{l,j}\mathcal
E_{i,l-1}\\
&=s^2\mathcal E_{i,j}e_l-r^2e_l\mathcal E_{i,j}.
\end{split}
\end{equation*}

Combining the above results, we have $(r^2+s^2)e_l\mathcal
E_{i,j}=(r^2+s^2)\mathcal E_{i,j}e_l$. By assumption $r^2+s^2\ne 0$,
we obtain (1), hence (2).

(3) is a direct consequence of (1).

(4) follows from (3) and Lemma 3.2 (3) since $\mathcal
E_{l,j}=[\mathcal E_{l,j-1},e_j]_{r^2}$ by Lemma 3.1 (3).

(5) can be proved similarly to (4), by (3), Lemma 3.1 (3) and Lemma
3.2 (1).

(6), (7): can be proved similarly to (1) \& (2), respectively.

(8) follows from (6).

(9) follows from (8), together with Lemma 3.2 (2).
\end{proof}

\begin{lemm} The following relations hold in $U^{+}:$

\smallskip
$(1)\quad
\mathcal{E}_{i,n}\mathcal{E}_{j,n}-(rs)\mathcal{E}_{j,n}\mathcal{E}_{i,n}
=\mathcal{E}_{j,n-1}\mathcal{E}_{i,n'}-r^2\mathcal{E}_{i,n'}\mathcal{E}_{j,n-1},
\qquad\hskip0.5cm i{<}j{<}n;$

$(2)\quad
\mathcal{E}_{j,k}\mathcal{E}_{i,k{+}1'}{-}r^2\mathcal{E}_{i,k{+}1'}\mathcal{E}_{j,k}
=\mathcal{E}_{i,k{+}2'}
\mathcal{E}_{j,k{+}1}{-}s^{-2}\mathcal{E}_{j,k{+}1}\mathcal{E}_{i,k{+}2'},
\ i{<}j{<}k{<}n{-}1;$

$(3)\quad
\mathcal{E}_{n-1,n}\mathcal{E}_{n-1,n'}=s^2\mathcal{E}_{n-1,n'}\mathcal{E}_{n-1,n};$

$(4)\quad \mathcal{E}_{i,j'}e_n=(rs)^2e_n\mathcal{E}_{i,j'},
\quad\quad\,\qquad\hskip4.25cm i{<}j{<}n;$

$(5)\quad
\mathcal{E}_{i,n'}\mathcal{E}_{n-1,n}=\mathcal{E}_{n-1,n}\mathcal{E}_{i,n'},
\qquad\quad\hskip4.06cm i{<}n{-}1;$

$(6)\quad
\mathcal{E}_{i,n'}\mathcal{E}_{n-1,n'}=s^2\mathcal{E}_{n-1,n'}\mathcal{E}_{i,n'},
\qquad\quad\hskip3.548cm i{<}n{-}1;$

$(7)\quad
\mathcal{E}_{i,n-1'}\mathcal{E}_{n-1,n}=s^2\mathcal{E}_{n-1,n}\mathcal{E}_{i,n-1'},
\qquad\quad\hskip3.02cm i{<}n{-}1;$

$(8)\quad
\mathcal{E}_{i,n-1'}\mathcal{E}_{n-1,n'}=(rs^2)^2\mathcal{E}_{n-1,n'}\mathcal{E}_{i,n-1'},
\qquad\quad\hskip2.24cm i{<}n{-}1.$
\end{lemm}
\begin{proof}
(1): By Lemma 3.3 (3), we have $\mathcal E_{i,n}\mathcal
E_{j,n-1}=\mathcal E_{j,n-1}\mathcal E_{i,n}$. So by Lemma 3.1 (3),
we obtain
\begin{equation*}
\begin{split}
\mathcal{E}_{i,n}\mathcal{E}_{j,n} &=
\mathcal{E}_{i,n}(\mathcal{E}_{j,n-1}e_{n}
-r^2e_{n}\mathcal{E}_{j,n-1})
 \\
&=\mathcal{E}_{j,n-1}\mathcal{E}_{i,n}e_{n}
-r^2\mathcal{E}_{i,n}e_{n}\mathcal{E}_{j,n-1}\\ &=
 \mathcal{E}_{j,n-1}\mathcal{E}_{i,n'}+rs\mathcal{E}_{j,n-1}e_{n}\mathcal{E}_{i,n}
-r^2\mathcal{E}_{i,n'}\mathcal{E}_{j,n-1}
-r^3se_{n}\mathcal{E}_{i,n}\mathcal{E}_{j,n-1} \\
&=\mathcal{E}_{j,n-1}\mathcal{E}_{i,n'}+rs\mathcal{E}_{j,n}\mathcal{E}_{i,n}
-r^2\mathcal{E}_{i,n'}\mathcal{E}_{j,n-1}.
\end{split}
\end{equation*}

(2): By Lemma 3.3 (8), we have $\mathcal E_{j,k}\mathcal
E_{i,k+2'}=\mathcal E_{i,k+2'}\mathcal E_{j,k}$. Using $\mathcal
E_{j,k+1}=\mathcal E_{j,k}e_{k+1}-r^2e_{k+1}\mathcal E_{j,k}$, we
obtain
\begin{equation*}
\begin{split}
\mathcal{E}_{j,k}\mathcal{E}_{i,k+1'} &=
\mathcal{E}_{j,k}\Bigl(\mathcal{E}_{i,k+2'}e_{k+1}
-s^{-2}e_{k+1}\mathcal{E}_{i,k+2'}\Bigr)\\
&=\mathcal{E}_{i,k+2'}\mathcal{E}_{j,k}e_{k+1}
-s^{-2}\mathcal{E}_{j,k}e_{k+1}\mathcal{E}_{i,k+2'})\\
&=\mathcal{E}_{i,k+2'}\mathcal{E}_{j,k+1}+r^2\mathcal{E}_{i,k+2'}e_{k+1}\mathcal{E}_{j,k}\\
&\quad -s^{-2}\mathcal{E}_{j,k+1}\mathcal{E}_{i,k+2'}
-r^2s^{-2}e_{k+1}\mathcal{E}_{i,k+2'}\mathcal{E}_{j,k} \\
&=r^2\mathcal{E}_{i,k+1'}\mathcal{E}_{j,k}+\mathcal{E}_{i,k+2'}
\mathcal{E}_{j,k+1}-s^{-2}\mathcal{E}_{j,k+1}\mathcal{E}_{i,k+2'}.
\end{split}
\end{equation*}

(3) can be proved by the same argument used in checking the
commutation relation (6) of Lemma 3.6 in \cite{BGH1} (where the
identity $\mathcal E_{1'2}\mathcal E_{12}=r^2\mathcal E_{12}\mathcal
E_{1'2}$ was actually verified for type $C_2$ case by only using
$(r,s)$-Serre relations).

(4): It suffices to check the case when $j=n-1$, i.e., $\mathcal
E_{i,n-1'}e_n=(rs)^2e_n\mathcal E_{i,n-1'}$, since for any $j\le
n-2$, $e_ne_j=e_je_n$ and $\mathcal E_{i,j'}=[\cdots[\,\mathcal
E_{i,n-1'}, e_{n-2}]_{s^{-2}},\cdots,e_j\,]_{s^{-2}}$, then
$\mathcal E_{i,j'}e_n=(rs)^2e_n\mathcal E_{i,j'}$.

Observe that (4) in the case when $j=n{-}1$ is equivalent to (5):
\begin{equation*}
\begin{split}
&\mathcal E_{i,n-1'}e_n{-}(rs)^2e_n\mathcal E_{i,n-1'}\\
&\quad=\bigl(\mathcal E_{i,n'}e_{n-1}{-}s^{-2}e_{n-1}\mathcal
E_{i,n'}\bigr)e_n{-}(rs)^2e_n\bigl(\mathcal
E_{i,n'}e_{n-1}{-}s^{-2}e_{n-1}\mathcal
E_{i,n'}\bigr)\\
&\quad=\mathcal E_{i,n'}e_{n-1}e_n-e_{n-1}e_n\mathcal E_{i,n'}
{-}r^2\mathcal E_{i,n'}e_ne_{n-1}{+}r^2e_ne_{n-1}\mathcal E_{i,n'}\\
&\quad=\mathcal E_{i,n'}\mathcal E_{n-1,n}-\mathcal
E_{n-1,n}\mathcal E_{i,n'}.
\end{split}
\end{equation*}
Since for $i<n-2$, we have the decompositions $\mathcal
E_{i,n-1'}=[\,\mathcal E_{i,n-3},\mathcal E_{n-2,n-1'}]_{r^2}$ and
$\mathcal E_{i,n'}=[\,\mathcal E_{i,n-3},\mathcal
E_{n-2,n'}]_{r^2}$, by Lemma 3.1 (4). Thus, this reduces to verify
$\mathcal E_{n-2,n-1'}e_n=(rs)^2e_n\mathcal E_{n-2,n-1'}$ or
equivalently, $\mathcal E_{n-2,n'}\mathcal E_{n-1,n}$ $=\mathcal
E_{n-1,n}\mathcal E_{n-2,n'}$. While the latter is equivalent to (3)
under assumption $r^3\ne s^3$.

Indeed, using (2.3) and (1), we at first have
\begin{equation*}
\begin{split}
&\mathcal{E}_{n-2,n}\mathcal{E}_{n-1,n'}{-}(rs)^2\mathcal{E}_{n-1,n'}\mathcal{E}_{n-2,n}\\
&\quad=
\bigl(\mathcal{E}_{n-2,n}\mathcal{E}_{n-1,n}\bigr)e_n{-}rs\bigl(\mathcal{E}_{n-2,n}e_{n}\bigr)\mathcal{E}_{n-1,n}
{-}(rs)^2\mathcal{E}_{n-1,n'}\mathcal{E}_{n-2,n}\quad\text{(by (1),
(2.3))}\\
 &\quad=\Bigl(rs\mathcal{E}_{n-1,n}\mathcal{E}_{n-2,n}
{+}e_{n-1}\mathcal{E}_{n-2,n'}-r^2\mathcal{E}_{n-2,n'}e_{n-1}\Bigr)e_n\\
&\qquad{-}rs\mathcal{E}_{n-2,n'}\mathcal{E}_{n-1,n}
{-}(rs)^2e_n\bigl(\mathcal{E}_{n-2,n}\mathcal{E}_{n-1,n}\bigr)
{-}(rs)^2\mathcal{E}_{n-1,n'}\mathcal{E}_{n-2,n}
\end{split}
\end{equation*}
\begin{equation*}
\begin{split}
&\quad=rs\mathcal{E}_{n-1,n}\bigl(\mathcal{E}_{n-2,n'}{+}\underline{rse_n\mathcal
E_{n-2,n}}\bigr){+}s^2e_{n-1}e_n\mathcal{E}_{n-2,n'}
{-}r^2\mathcal{E}_{n-2,n'}e_{n-1}e_n\\
&\qquad{-}rs\mathcal{E}_{n-2,n'}\mathcal{E}_{n-1,n}
{-}(rs)^2e_n\Bigl(\underline{rs\mathcal{E}_{n-1,n}\mathcal{E}_{n-2,n}}
{+}e_{n-1}\mathcal{E}_{n-2,n'}{-}r^2\mathcal{E}_{n-2,n'}e_{n-1}\Bigr)\\
&\qquad{-}\underline{(rs)^2\mathcal{E}_{n-1,n'}\mathcal{E}_{n-2,n}} \qquad\qquad\text{(applying
(2.3) to the terms underlined) }\\
&\quad=rs\mathcal{E}_{n-1,n}\mathcal{E}_{n-2,n'}{+}s^2e_{n-1}e_n\mathcal{E}_{n-2,n'}{-}
r^2\mathcal{E}_{n-2,n'}e_{n-1}e_n\\
&\qquad {-}rs\mathcal{E}_{n-2,n'}\mathcal{E}_{n-1,n}
{-}(rs)^2e_ne_{n-1}\mathcal{E}_{n-2,n'}{+}r^4\mathcal{E}_{n-2,n'}e_ne_{n-1}\\
 &\quad=(rs{+}s^2)\mathcal{E}_{n-1,n}\mathcal{E}_{n-2,n'}{-}(r^2{+}rs)\mathcal{E}_{n-2,n'}\mathcal{E}_{n-1,n}.
\end{split}
\end{equation*}

Next, using Lemma 3.1 (3), together with (3): $\mathcal
E_{n-1,n}\mathcal E_{n-1,n'}=s^2\mathcal E_{n-1,n'}\mathcal
E_{n-1,n}$ and Lemma 3.1 (4), we have another expansion
\begin{equation*}
\begin{split}
&\mathcal{E}_{n-2,n}\mathcal{E}_{n-1,n'}{-}(rs)^2\mathcal{E}_{n-1,n'}\mathcal{E}_{n-2,n}\\
&\quad=e_{n-2}\mathcal E_{n-1,n}\mathcal E_{n-1,n'}-r^2\mathcal
E_{n-1,n}e_{n-2}\mathcal
E_{n-1,n'}\\
&\qquad-(rs)^2\mathcal{E}_{n-1,n'}e_{n-2}\mathcal
E_{n-1,n}+(r^2s)^2\mathcal{E}_{n-1,n'}\mathcal E_{n-1,n}e_{n-2}\\
&\quad=s^2e_{n-2}\mathcal E_{n-1,n'}\mathcal E_{n-1,n}-r^2\mathcal
E_{n-1,n}e_{n-2}\mathcal E_{n-1,n'}\\
&\qquad-(rs)^2\mathcal E_{n-1,n'}e_{n-2}\mathcal
E_{n-1,n}+r^4\mathcal E_{n-1,n}\mathcal E_{n-1,n'}e_{n-2}\\
&\quad=s^2\mathcal E_{n-2,n'}\mathcal E_{n-1,n}-r^2\mathcal
E_{n-1,n}\mathcal E_{n-2,n'}.
\end{split}
\end{equation*}

Combining the above two expansions, we obtain
$$
\mathcal{E}_{n-2,n'}\mathcal{E}_{n-1,n}=\mathcal{E}_{n-1,n}\mathcal{E}_{n-2,n'},
$$ under assumption $r^2+rs+s^2\neq 0$.

(6) \& (7) follow from (5), owing to (2.5) \& (2.7).

(8) follows from (7) \& (4).
\end{proof}

\begin{lemm} The following relations hold in $U^{+}:$

\smallskip
$(1)\quad
\mathcal{E}_{i,n'}\mathcal{E}_{j,n}=\mathcal{E}_{j,n}\mathcal{E}_{i,n'},
\qquad\hskip3.05cm\qquad 1\le i<j<n;$

$(2)\quad
\mathcal{E}_{i,n}\mathcal{E}_{j,n'}-(rs)^2\mathcal{E}_{j,n'}\mathcal{E}_{i,n}
=(s^2-r^2)\mathcal{E}_{i,n'}\mathcal{E}_{j,n},\qquad 1\le i<j<n;$

$(3)\quad \mathcal E_{i,n}\mathcal E_{i,j'}=s^2\mathcal
E_{i,j'}\mathcal E_{i,n}, \qquad\hskip3.56cm 1\le i<j\le n;$

$(4)\quad
\mathcal{E}_{i,n'}\mathcal{E}_{j,n'}=s^2\mathcal{E}_{j,n'}\mathcal{E}_{i,n'},
\qquad\hskip3.26cm 1\le i<j<n;$

$(5)\quad \mathcal E_{i,l'}\mathcal E_{j,n}=\mathcal E_{j,n}\mathcal
E_{i,l'},\hskip3.82cm\;\, \qquad i<j<l\le n;$

$(6)\quad \mathcal E_{i,l'}\mathcal E_{j,n'}=(rs)^2\mathcal
E_{j,n'}\mathcal E_{i,l'},\hskip3.38cm \quad i<j<l< n.$
\end{lemm}
\begin{proof}
(1): Note that $\mathcal E_{i,n'}\mathcal E_{n-1,n}=\mathcal
E_{n-1,n}\mathcal E_{i,n'}$ by Lemma 3.4 (5). So for $j\le n{-}2$,
$\mathcal E_{i,n'}\mathcal E_{j,n-2}=\mathcal E_{j,n-2}\mathcal
E_{i,n'}$ since $\mathcal E_{i,n'}=[\,\mathcal E_{i,n},e_n]_{rs}$
and both $\mathcal E_{i,n}$ and $e_n$ commute with $\mathcal
E_{j,n-2}$ (by Lemma 3.3 (3)). Hence, Lemma 3.1 (3) gives us a
decomposition $\mathcal E_{j,n}=[\,\mathcal E_{j,n-2},\mathcal
E_{n-1,n}]_{r^2}$, from which we obtain the required result.

(2): By (2.3), and using Lemma 3.2 (4), Lemma 3.4 (1), we have
\begin{equation*}
\begin{split}
\mathcal{E}_{i,n}&\mathcal{E}_{j,n'}{-}(rs)^2\mathcal{E}_{j,n'}\mathcal{E}_{i,n}
\\
&=(\mathcal E_{i,n}\mathcal E_{j,n})e_n{-}rs(\mathcal
E_{i,n}e_n)\mathcal E_{j,n}{-}(rs)^2\mathcal E_{j,n}(e_n\mathcal
E_{i,n})+(rs)^3e_n\mathcal E_{j,n}\mathcal E_{i,n}\\
&=rs\mathcal E_{j,n}\mathcal E_{i,n}e_n{+}\mathcal E_{j,n-1}\mathcal
E_{i,n'}e_n{-}r^2\mathcal E_{i,n'}\mathcal
E_{j,n-1}e_n\\
&\qquad{-}rs(\mathcal E_{i,n'}{+}rse_n\mathcal E_{i,n})\mathcal
E_{j,n}{-}(rs)^2\mathcal E_{j,n}(e_n\mathcal
E_{i,n}){+}(rs)^3e_n\mathcal E_{j,n}\mathcal E_{i,n}
\end{split}
\end{equation*}
\begin{equation*}
\begin{split}&=rs\mathcal E_{j,n}\mathcal E_{i,n'}-rs\mathcal E_{i,n'}\mathcal
E_{j,n}{+}s^2\mathcal E_{j,n-1}e_n\mathcal E_{i,n'}{-}r^2\mathcal
E_{i,n'}\mathcal E_{j,n-1}e_n\\
&\qquad{-}(rs)^2e_n(\mathcal E_{i,n}\mathcal E_{j,n}{-}rs\mathcal
E_{j,n}\mathcal E_{i,n})\\
&=rs\mathcal E_{j,n}\mathcal E_{i,n'}-rs\mathcal E_{i,n'}\mathcal
E_{j,n}{+}s^2\mathcal E_{j,n-1}e_n\mathcal E_{i,n'}{-}r^2\mathcal
E_{i,n'}\mathcal E_{j,n-1}e_n\\
&\qquad{-}(rs)^2e_n(\mathcal E_{j,n-1}\mathcal E_{i,n'}-r^2\mathcal
E_{i,n'}\mathcal E_{j,n-1})\\
&=rs[\mathcal E_{j,n},\mathcal E_{i,n'}]{+}s^2\mathcal
E_{j,n}\mathcal E_{i,n'}-r^2\mathcal E_{i,n'}\mathcal E_{j,n}
\end{split}
\end{equation*}
which gives the required result using (1).

(3): First of all, we claim that $\mathcal E_{i,n}\mathcal
E_{i,n'}=s^2\mathcal E_{i,n'}\mathcal E_{i,n}$. When $i=n-1$, this
is Lemma 3.4 (3). For $i<n-1$, it follows from (1) since $\mathcal
E_{i+1,n}\mathcal E_{i,n'}=\mathcal E_{i,n'}\mathcal E_{i+1,n}$ and
$e_i\mathcal E_{i,n'}=s^2\mathcal E_{i,n'}e_i$, as well as $\mathcal
E_{i,n}=[\,e_i,\mathcal E_{i+1,n}]_{r^2}$.

Next, noting  $\mathcal E_{i,j'}=[\cdots[\,\mathcal E_{i,n'},
e_{n-1}]_{s^{-2}},\cdots,e_j]_{s^{-2}}$ and $\mathcal
E_{i,n}e_k=e_k\mathcal E_{i,n}$ for $i<j\le k<n$, together with the
above assertion, we obtain the required result.

(4) follows from (1), Lemma 3.2 (4) as well as $\mathcal
E_{j,n'}=[\,\mathcal E_{j,n}, e_n]_{rs}$.

(5) follows from (1) and $\mathcal E_{i,l'}=[\cdots[\,\mathcal
E_{i,n'},e_{n-1}]_{s^{-2}},\cdots,e_{l}]_{s^{-2}}$, as well as
$e_k\mathcal E_{j,n}$ $=\mathcal E_{j,n}e_k$ for any $k$ with
$j<l\le k<n$.

(6) follows from (5) \& Lemma 3.4 (4).
\end{proof}

\begin{lemm} The following relations hold in $U^{+}:$

\smallskip
$(1)\quad  [\,\mathcal{E}_{i,j'},e_{i+1}]_{\bullet}=0,\qquad\ \,
\;\hskip3.66cm 1\le i< j\le n, \ i\ne j-2;$

$(2)\quad \mathcal E_{i,j'}e_j=r^{-2}e_j\mathcal
E_{i,j'},\qquad\hskip3.55cm 1\le i<j<n;$

$(3)\quad
[\,\mathcal{E}_{i,i+1'},e_{i+2}]_{\bullet}=0,\qquad\hskip3.63cm 1\le
i< n-1;$

$(4)\quad \mathcal{E}_{i,j'}e_l=e_l\mathcal{E}_{i,j'}, \ \qquad
\hskip4.08cm   i<j<l< n;$

$(5)\quad \mathcal{E}_{i,l'}\mathcal{E}_{i,j'}=r^{-2}\mathcal
E_{i,j'}\mathcal E_{i,l'},\hskip3.07cm \qquad i<j<l\le n;$

$(6)\quad \mathcal E_{i,l'}\mathcal E_{j,k'}=(rs)^2\mathcal
E_{j,k'}\mathcal E_{i,l'},\qquad\hskip2.76cm i<j<l<k\le n;$

$(7)\quad \mathcal E_{i,l'}\mathcal E_{j,l'}=s^2\mathcal
E_{j,l'}\mathcal E_{i,l'},\hskip3.335cm\qquad i<j<l< n;$

$(8)\quad
\mathcal{E}_{i,n-1}\mathcal{E}_{i,n'}-(rs)^2\mathcal{E}_{i,n'}\mathcal{E}_{i,n-1}=s(s{-}r)
\mathcal{E}_{i,n}^2,
  \qquad 1\leq i<n.$
\end{lemm}
\begin{proof}
See the appendix.
\end{proof}

\begin{lemm} The following relations hold in $U^{+}:$

\smallskip

$(1)\quad \mathcal
E_{i,k}\mathcal{E}_{i,k+1'}{-}(rs)^2\mathcal{E}_{i,k+1'}\mathcal
E_{i,k} =s^2\mathcal{E}_{i,k+2'}\mathcal{E}_{i,k+1}{-}
s^{-2}\mathcal{E}_{i,k+1}\mathcal{E}_{i,k+2'},\ \, k<n{-}1;$

$(2)\quad
\mathcal{E}_{j,j+1}\mathcal{E}_{j,j+1'}=(rs^2)^2\mathcal{E}_{j,j+1'}\mathcal{E}_{j,j+1},
\hskip3.66cm\qquad j<n-1;$

$(3)\quad
\mathcal{E}_{j-1,j'}\mathcal{E}_{j,j+1'}=(rs^2)^2\mathcal{E}_{j,j+1'}\mathcal{E}_{j-1,j'},
\hskip3.43cm\qquad j<n-1;$

$(4)\quad
\mathcal{E}_{i,j'}\mathcal{E}_{j,j+1'}=(rs^2)^2\mathcal{E}_{j,j+1'}\mathcal{E}_{i,j'},
\hskip4cm\quad i<j<n-1;$

$(5)\quad \mathcal E_{j,k}\mathcal E_{j,j+1'}=(rs^2)^2\mathcal
E_{j,j+1'}\mathcal E_{j,k}, \qquad\hskip4.25cm j<k<n;$

$(6)\quad \mathcal E_{j-1,j}\mathcal E_{i,j'}=(rs)^2\mathcal
E_{i,j'}\mathcal E_{j-1,j}, \hskip3.9cm\qquad i+1<j<n;$

$(7)\quad \mathcal E_{l,k}\mathcal E_{i,j'}=(rs)^2\mathcal
E_{i,j'}\mathcal E_{l,k}, \hskip4cm\qquad i<l<j\le k<n;$

$(8)\quad \mathcal E_{i,k}\mathcal E_{i,j'}=(rs^2)^2\mathcal
E_{i,j'}\mathcal E_{i,k}, \hskip4.4cm\qquad i<j\le k<n.$
\end{lemm}
\begin{proof}
(1): The proof is similar to that of Lemma 3.4 (2) and noting that
$\mathcal E_{i,k}\mathcal E_{i,k+2'}=s^2\mathcal E_{i,k+2'}\mathcal
E_{i,k}$.

(2): Using (1) and Lemma 3.2 (5), we obtain
\begin{equation*}
\begin{split}
\mathcal{E}_{j,j+1}\mathcal{E}_{j,j+1'}&=(e_je_{j+1}-r^2e_{j+1}e_j)\mathcal{E}_{j,j+1'}\\
&=r^2e_j\mathcal{E}_{j,j+1'}e_{j+1}-r^2e_{j+1}e_j\mathcal{E}_{j,j+1'}\\
&=r^2\Bigl((rs)^2\mathcal{E}_{j,j+1'}e_j+s^2\mathcal{E}_{j,j+2'}\mathcal{E}_{j,j+1}-
s^{-2}\mathcal{E}_{j,j+1}\mathcal{E}_{j,j+2'}\Bigr)e_{j+1}\\
&\quad
-r^2e_{j+1}\Bigl((rs)^2\mathcal{E}_{j,j+1'}e_j+s^2\mathcal{E}_{j,j+2'}\mathcal{E}_{j,j+1}-
s^{-2}\mathcal{E}_{j,j+1}\mathcal{E}_{j,j+2'}\Bigr)\\
&=(r^2s)^2\mathcal{E}_{j,j+1'}\mathcal{E}_{j,j+1}+
(rs^2)^2\mathcal{E}_{j,j+1'}\mathcal{E}_{j,j+1}-(rs^{-1})^2\mathcal{E}_{j,j+1}\mathcal{E}_{j,j+1'}.
\end{split}
\end{equation*}
So we have
$$(1+r^2s^{-2})\mathcal{E}_{j,j+1}\mathcal{E}_{j,j+1'}=(rs^2)^2(1+r^2s^{-2})\mathcal{E}_{j,j+1'}\mathcal{E}_{j,j+1},$$
since $1+r^2s^{-2}\neq 0$, we get
$$\mathcal{E}_{j,j+1}\mathcal{E}_{j,j+1'}=(rs^2)^2\mathcal{E}_{j,j+1'}\mathcal{E}_{j,j+1}.$$

(3): By Lemma 3.6 (2), (4) \& Lemma 3.4 (4), we have
\begin{equation*}
\begin{split}
\mathcal{E}_{j-1,j'}\mathcal{E}_{j,j+2'}&=\mathcal{E}_{j-1,j'}\bigl(e_{j}\mathcal{E}_{j+1,j+2'}
-r^{2}\mathcal{E}_{j+1,j+2'}e_{j}\bigr)\\
&=(rs^2)^2\bigl(e_{j}\mathcal{E}_{j+1,j+2'}
-r^{2}\mathcal{E}_{j+1,j+2'}e_{j}\bigr)\mathcal{E}_{j-1,j'}\\
&=(rs^2)^2\mathcal{E}_{j,j+2'}\mathcal{E}_{j-1,j'}.
 \end{split}
\end{equation*}
Again, by Lemma 3.6 (4), we arrive at
\begin{equation*}
\begin{split}
\mathcal{E}_{j-1,j'}\mathcal{E}_{j,j+1'}&=\mathcal{E}_{j-1,j'}\bigl(\mathcal{E}_{j,j+2'}e_{j+1}
-s^{-2}e_{j+1}\mathcal{E}_{j,j+2'}\bigr)\\
&=(rs^2)^2\bigl(\mathcal{E}_{j,j+2'}e_{j+1}
-s^{-2}e_{j+1}\mathcal{E}_{j,j+2'}\bigr)\mathcal{E}_{j-1,j'}\\
&=(rs^2)^2\mathcal{E}_{j,j+1'}\mathcal{E}_{j-1,j'}.
 \end{split}
\end{equation*}

(4) follows from (3), together with $e_l\mathcal E_{j,j+1'}=\mathcal
E_{j,j+1'}e_l$ for $i\le l<j-1$.

(5) follows from (2), together with Lemma 3.6 (4).

(6) An observation when $i=j-2$ comes from the proof of Lemma 3.6
(3), for $j<n$. For $i<j-2$, it follows from Lemma 3.1 (1).

(7) follows from (6), together with Lemmas 3.3 (6) \& 3.6 (4).

(8) follows from (7), together with Lemma 3.2 (2) for $i+1<j$, but
the case for $i+1=j$ follows from (2) as well as Lemma 3.6 (4).
\end{proof}

\noindent {\it 3.2. Central elements.} Restrict two parameters $r$
and $s$ to be roots
 of $1$: $r$ is a primitive $d$th root of unity, $s$ is a primitive
 $d'$th root of unity and $\ell$ is the least common multiple of $d$ and
 $d'$. From now on, we assume that $\Bbb{K}$ contains a primitive $\ell $th
 root of unity.

The following Lemma is useful to derive some commutator relations.

\begin{lemm}\label{3.8} Let $x, y, z$ be elements of a  $\mathbb{K}$-algebra such that  $yx=\a
xy+z$  for some $\a \in \mathbb{K}$ and $m$ a natural number. Then
the following assertions hold

$(1)$ If  $zx=\b xz$ for some $\b(\neq \a)  \in \mathbb{K}$, then
$yx^m=\a^mx^my+\frac{\a^m-\b^m}{\a-\b}x^{m-1}z$;

$(2)$ If  $yz=\b zy$ for some $\b(\neq \a)  \in \mathbb{K}$, then
$y^mx=\a^mxy^m+\frac{\a^m-\b^m}{\a-\b}zy^{m-1}$.
\end{lemm}

Based on Lemma 3.8 and by induction, one can easily check the
following Lemmas, where the $(r,s)$-integers, factorials and
binomial coefficients are defined for positive integers $c$ and $d$
by
$$[c]_i:=\frac{r_i^c-s_i^c}{r_i-s_i},\quad [c]_i!:=[c]_i[c-1]_i\cdots
[2]_i[1]_i,\quad
 \left[c\atop d\right]_i:=\frac{[c]_i!}{[d]_i![c-d]_i!}. $$
Write in brief $[m]=[m]_i$ for $i<n$. By convention $[0]_i=0$ and
$[0]_i!=1$. Denote
$\a_m=\frac{(r^m-s^m)(r^{m-1}-s^{m-1})}{(r-s)(r^2-s^2)}$,
$\beta_m=\frac{r^m-s^m}{r-s}$, then
$\a_{m+1}=s^2\a_m+r^{m-1}\beta_m$ and $\beta_{m+1}=s\beta_m+r^m$.

\medskip For $\mathcal E_{i,j}^m$ ($i\le j\le n$) and $\mathcal
E_{i,n'}^m$ ($i<n$) with $m\in\mathbb Z_+$, we have

\begin{lemm}\label{3.9} For $m\in \mathbb Z_+$ and $i\leq j<n$, the
following equalities hold

\smallskip
$(1)\quad
e_{i-1}\mathcal{E}_{i,j}^m=r^{2m}\mathcal{E}_{i,j}^me_{i-1} +[m]\
\mathcal{E}_{i,j}^{m-1}\mathcal{E}_{i-1,j};$

\smallskip
$(2)\quad
\mathcal{E}_{i,j}^me_{j+1}=r^{2m}e_{j+1}\mathcal{E}_{i,j}^m +[m]\
\mathcal{E}_{i,j+1}\mathcal{E}_{i,j}^{m-1};$

\smallskip
$(3)\quad e_{n-1}e_n^m= \a_m
e^{m-2}_n\mathcal{E}_{n-1,n'}+r^{m-1}\b_m
e^{m-1}_n\mathcal{E}_{n-1,n}+r^{2m}e_n^me_{n-1};$

\smallskip
$(4)\quad \mathcal{E}_{i,n}^me_{n}=(rs)^{m}e_n\mathcal{E}_{i,n}^m
+s^{m-1}\b_m\ \mathcal{E}_{i,n'}\mathcal{E}_{i,n}^{m-1};$

\smallskip
$(5)$\quad $e_{i-1}\mathcal{E}_{i,n}^m=r^{m-1}\b_m\
\mathcal{E}_{i,n}^{m-1}\mathcal{E}_{i-1,n}  + \a_m\
\mathcal{E}_{i,n}^{m-2}\mathcal{E}_{i,n-1}\mathcal{E}_{i-1,n'}$

\smallskip
$\qquad \qquad \qquad \quad -r^2\a_m\
\mathcal{E}_{i,n}^{m-2}\mathcal{E}_{i-1,n'} \mathcal{E}_{i,n-1}+
r^{2m}\mathcal{E}_{i,n}^me_{i-1};$

\smallskip
$(6)\quad e_{n-1}\mathcal{E}_{n-1,n'}^m= s^{2m-1}(s-r)[m]\
\mathcal{E}_{n-1,n'}^{m-1}\mathcal{E}_{n-1,n}^2+(rs)^{2m}\mathcal{E}_{n-1,n'}^me_{n-1};$

\smallskip
$(7)\quad \mathcal{E}_{i,n'}^me_{n-1}= s^{2(1-m)}[m]\
\mathcal{E}_{i,n'}^{m-1}\mathcal{E}_{i,n-1'}+s^{-2m}e_{n-1}\mathcal{E}_{i,n'}^m,
\qquad i<n-1;$

\smallskip
$(8)\quad  e_{i-1}\mathcal{E}_{i,n'}^m= [m]\
\mathcal{E}_{i,n'}^{m-1}\mathcal{E}_{i-1,n'}+r^{2m}\mathcal{E}_{i,n'}^me_{i-1},\qquad
i<n.$
\end{lemm}

For $\mathcal E_{i,j'}^m$ with $m\in\mathbb Z_+$, we have

\begin{lemm}\label{3.9}
\ For $m\in\mathbb Z_+$, the following equalities hold

\smallskip
$(1)\quad
e_{i-1}\mathcal{E}_{i,j'}^m=r^{2m}\mathcal{E}_{i,j'}^me_{i-1} +[m]\
\mathcal{E}_{i,j'}^{m-1}\mathcal{E}_{i-1,j'},\hskip2.8cm\quad i<j\le
n;$

\smallskip
$(2)$\quad
$e_{j-1}\mathcal{E}_{j-1,j'}^m=(rs)^{2m}\mathcal{E}_{j-1,j'}^me_{j-1}
+s^{2(m-1)}[m]\ \mathcal{E}_{j-1,j'}^{m-1}\times$

\smallskip
$\qquad\qquad\qquad\qquad\quad
\times\Bigl(s^2\mathcal{E}_{j-1,j+1'}\mathcal{E}_{j-1,j}-s^{-2}\
\mathcal{E}_{j-1,j}\mathcal{E}_{j-1,j+1'}\Bigr),\quad \qquad j<n;$

\smallskip
$(3)\quad
\mathcal{E}_{i,j'}^me_{j-1}=s^{-2m}e_{j-1}\mathcal{E}_{i,j'}^m
+(rs)^{-2(m-1)}[m]\
\mathcal{E}_{i,j-1'}\mathcal{E}_{i,j'}^{m-1},\qquad i{+}1<j\le n;$

\smallskip
$(4)\quad \mathcal E_{i,j{+}1'}^m\mathcal E_{i,j}\equiv
(rs)^{{-}2m}\Bigl(\mathcal E_{i,j}\mathcal
E_{i,j{+}1'}^m+[m]\,\mathcal E_{i,j{+}1}\mathcal
E_{i,j{+}2'}\mathcal E_{i,j{+}1'}^{m{-}1}\Bigr), \ \, i\le j<n{-}1;$

\smallskip
$(5)\quad \mathcal E_{i,n'}^m\mathcal
E_{i,n-1}=(rs)^{-2m}\Bigl(\mathcal E_{i,n-1}\mathcal
E_{i,n'}^m+s^3(r{-}s)[m]\,\mathcal E_{i,n}^2\mathcal
E_{i,n'}^{m-1}\Bigr);$

\smallskip
$(6)\quad \mathcal E_{i,j+2'}e_{j+1}^m=s^{-2m}e_{j+1}^m\mathcal
E_{i,j+2'}+(rs)^{2(1-m)}[m]\,e_{j+1}^{m-1}\mathcal E_{i,j+1'}, \ \,
i\le j<n{-}1.$
\end{lemm}

Next, we consider the commutation relations of $\mathcal E_{i,j}^m$
with $f_i$ ($1\le i\le n$).

\begin{lemm}
\ For $m\in\mathbb Z_+$, the following equalities hold:

\smallskip
\smallskip
$(1)\quad \mathcal{E}_{i,j}^mf_{i}=f_{i}\mathcal{E}_{i,j}^m -[m]\
\mathcal{E}_{i+1,j}\mathcal{E}_{i,j}^{m-1}\omega_{i},\qquad\hskip3.6cm
i<j<n;$

\smallskip
$(2)\quad \mathcal{E}_{i,j}^mf_{j}=f_{j}\mathcal{E}_{i,j}^m
+(rs)^{-2(m-1)}[m]\
\mathcal{E}_{i,j-1}\mathcal{E}_{i,j}^{m-1}\omega_{j}',\quad\hskip2.1cm
i<j<n;$

\smallskip
$(3)
\quad e_i^mf_i=f_ie_i^m+(r_i-s_i)^{-1}[m]_i\
e_i^{m-1}(s_i^{-m+1}\omega_i -r_i^{-m+1}\omega_i'),\qquad 1\leq i
\leq n;$

\smallskip
$(4)\quad  \mathcal{E}_{i,n}^mf_i=f_i\mathcal{E}_{i,n}^m
-r^{m-1}\b_m \mathcal{E}_{i+1,n}\mathcal{E}_{i,n}^{m-1}\omega_i
-\a_m
\mathcal{E}_{i+1,n-1}\mathcal{E}_{i,n'}\mathcal{E}_{i,n}^{m-2}\omega_i$

$\qquad\qquad\qquad\qquad +\,r^2\a_m \mathcal{E}_{i,n'}
\mathcal{E}_{i+1,n-1}\mathcal{E}_{i,n}^{m-2}\omega_i,\hskip2.5cm
\qquad i<n-1;$

\smallskip
$(5)\quad \mathcal{E}_{n-1,n}^mf_{n-1}=f_{n-1}\mathcal{E}_{n-1,n}^m
-r^{m-1}\b_me_n\mathcal{E}_{n-1,n}^{m-1}\omega_{n-1} -\a_m
\mathcal{E}_{n-1,n'}\mathcal{E}_{n-1,n}^{m-2}\omega_{n-1};$

\smallskip
$(6)\quad \mathcal{E}_{i,n}^mf_n=f_n\mathcal{E}_{i,n}^m
+(r+s)s^{-m-1}\,\beta_m\,
\mathcal{E}_{i,n-1}\omega_n'\mathcal{E}_{i,n}^{m-1},\qquad\hskip2.03cm
i<n.$
\end{lemm}

\begin{lemm}\label{3.12}
\ For $m\in\mathbb Z_+$, the following equalities hold

\smallskip
$(1)\quad \mathcal{E}_{i,n'}^mf_i=f_i\mathcal{E}_{i,n'}^m -[m]\
\mathcal{E}_{i+1,n'}\mathcal{E}_{i,n'}^{m-1}\omega_i,
 \qquad\hskip3cm  i<n-1;$

\smallskip
$(2)\quad
\mathcal{E}_{n-1,n'}^mf_{n-1}=f_{n-1}\mathcal{E}_{n-1,n'}^m
+(r{-}s)\,s^{2m-1}[m]\,
e_n^2\mathcal{E}_{n-1,n'}^{m-1}\omega_{n-1};$

\smallskip
$(3)\quad \mathcal{E}_{i,n'}^mf_{n}=f_{n}\mathcal{E}_{i,n'}^m
+(rs)^{-2m+1}(r{+}s)\,\b_{2m}\,
\mathcal{E}_{i,n}\mathcal{E}_{i,n'}^{m-1}\omega_{n}',\hskip1.4cm
i\leq n-1;$

\smallskip
$(4)\quad \mathcal{E}_{i,j'}^mf_i=f_i\mathcal{E}_{i,j'}^m -[m]\
\mathcal{E}_{i+1,j'}\mathcal{E}_{i,j'}^{m-1}\omega_i,
 \qquad\hskip3.2cm i<j-1;$

\smallskip
$(5)$\quad
$\mathcal{E}_{j-1,j'}^mf_{j-1}=f_{j-1}\mathcal{E}_{j-1,j'}^m
+s^{2(m-1)}[m]\times$

$\qquad\qquad\qquad\qquad\times\bigl(s^{-2}e_j\mathcal{E}_{j,j+1'}-\,s^{2}\mathcal{E}_{j,j+1'}e_j\bigr)\,
\mathcal{E}_{j-1,j'}^{m-1}\omega_{j-1},
 \qquad\ \;  j\le n-1;$

\smallskip
$(6)\quad \mathcal{E}_{i,j'}^mf_j=f_j\mathcal{E}_{i,j'}^m
+s^{-2}[m]\
\mathcal{E}_{i,j+1'}\mathcal{E}_{i,j'}^{m-1}\omega_{j}',\qquad
\hskip2.5cm i<j<n.$
\end{lemm}

\begin{prop}\label{3.14}
\ $\mathcal{E}_{k,j}^{\ell}(k\leq j\leq n),\
\mathcal{E}_{k,m'}^{\ell} (k<m\leq n) $ commute  with $f_i$ for
$1\leq i\leq n$.
\end{prop}
\begin{proof} The proof is carried out in three steps.

Step 1. We want to prove $\mathcal{E}_{k,j}^{\ell} \ \ (k\leq j\leq
n)$ commutes with $f_i$ for $1\leq i\leq n$. If $i<k$ or $i>j$, this
statement is immediate from $(B4)$; while if $i=k$ or $i=j$, it
follows directly from Lemma 3.12. We may suppose that $k<i<j$, using
$\omega_{i}'\mathcal{E}_{i+1,j}=r^2\mathcal{E}_{i+1,j}\omega_{i}'$
and Lemmas 3.1 (3) \& 3.12 (2), we obtain
\begin{equation*}
\begin{split}
\mathcal{E}_{k,j}f_i &=(\mathcal{E}_{k,i}\mathcal{E}_{i+1,j}-r^{2}\mathcal{E}_{i+1,j}\mathcal{E}_{k,i})f_i\\
&=(f_i\mathcal{E}_{k,i}+s^{-2}\mathcal{E}_{k,i-1}\omega_i')\mathcal{E}_{i+1,j}
 -r^{2}\mathcal{E}_{i+1,j}(f_i\mathcal{E}_{k,i}+s^{-2}\mathcal{E}_{k,i-1}\omega_i')\\
 &=f_i\mathcal{E}_{k,j}.
\end{split}
\end{equation*}

Step 2. We want to prove $\mathcal{E}_{k,n'}^{\ell} \ (k<n)$
commutes with $f_i$ for $1\leq i\leq n$. If $i<k$, this statement is
immediate from $(B4)$; while if $i=k$ or $i=n$, it follows directly
from Lemma \ref{3.12}; if $k<i<n$, it is a similar argument as Step
1.

\smallskip
Step 3. We want to prove $\mathcal{E}_{k,m'}^{\ell} \ (k<m<n)$
commutes with $f_i$ for $1\leq i\leq n$. If $i<k$, this statement is
immediate from $(B4)$; while if $i=k$ or $i=m$, it follows directly
from Lemma 3.13. We may suppose that $k\leq i\leq n$. If $k<i<m$,
using
$\omega_{i}'\mathcal{E}_{i+1,m'}=r^2\mathcal{E}_{i+1,m'}\omega_{i}'$
and Lemma 3.1 (2), we obtain
\begin{equation*}
\begin{split}
\mathcal{E}_{k,m'}f_i &=(\mathcal{E}_{k,i}\mathcal{E}_{i+1,m'}-r^{2}\mathcal{E}_{i+1,m'}\mathcal{E}_{k,i})f_i\\
&=(f_i\mathcal{E}_{k,i}+s^{-2}\mathcal{E}_{k,i-1}\omega_i')\mathcal{E}_{i+1,m'}
 -r^{2}\mathcal{E}_{i+1,m'}(f_i\mathcal{E}_{k,i}+s^{-2}\mathcal{E}_{k,i-1}\omega_i')\\
 &=f_i\mathcal{E}_{k,m'}.
\end{split}
\end{equation*}The last case is $m<i\leq n$. By Lemma 3.13 (6),
we have $[\,\mathcal{E}_{k,i'},
f_i]=s^{-2}\mathcal{E}_{k,i+1'}\omega_i'$. Noting that $[\,
\mathcal{E}_{k,i+1'}\omega_i', e_{i-1}]_{s^{-2}}=0$, we have
\begin{equation*}
\begin{split}
[\, \mathcal{E}_{k,j'}, f_i\,]&=[\,[\cdots [\, \mathcal{E}_{k,i'},
e_{i-1}]_{s^{-2}}, \cdots, e_j\,]_{s^{-2}}, f_i\,]\\
&=[\cdots,[\,[\,\mathcal{E}_{k,i'}, f_i\,], e_{i-1}\,]_{s^{-2}},
\cdots,
e_j\,]_{s^{-2}}\\
&=0.
\end{split}
\end{equation*}

This completes the proof.
\end{proof}
\begin{theorem}\label{3.15}
\ The elements $\mathcal{E}_{k,j}^{\ell}\ (1\leq k\leq j\leq n),\
\mathcal{E}_{k,j'}^{\ell}\ (1\leq k<j\leq n),\
\mathcal{F}_{k,j}^{\ell}\ (1\leq k\leq j\leq n),\
\mathcal{F}_{k,j'}^{\ell}\ (1\leq k<j\leq n)$ and
$\omega_k^{\ell}-1,\ (\omega_k')^{\ell}-1 \ (1\leq k\leq n)$ are
central in $U_{r,s}(\mathfrak{so}_{2n+1})$.
\end{theorem}
\begin{proof}
It follows from Propositions 3.11 \& 3.14 that
$\mathcal{E}_{k,j}^{\ell}$ $(1\leq k\leq j\leq n),\
\mathcal{E}_{k,n'}^{\ell}\ (1\leq k<n),\ \mathcal{E}_{k,j'}^{\ell}$
$(1\leq k<j\leq n-1)$ are central in
$U_{r,s}(\mathfrak{so}_{2n+1})$. Applying $\tau$ to $\mathcal
E_{k,j}^\ell$ \& $\mathcal E_{k,j'}^\ell$, we see that
$\mathcal{F}_{k,j}^{\ell}\ (1\leq k\leq j\leq n),\
\mathcal{F}_{k,j'}^{\ell}\ (1\leq k<j\leq n)$ are central. Easy to
see that $\omega_k^{\ell}-1,\ (\omega_k')^{\ell}-1 \ (1\leq k\leq
n)$ are central, too.
\end{proof}

\noindent {\it Remark:}\label{3.16} If $\ell=2\ell'$, then
$\mathcal{E}_{k,j}^{\ell'}\ (1\leq k\le j<n),\,
\mathcal{E}_{k,j'}^{\ell'}\ (1\leq k<j\leq n),\, e_n^\ell,\,
\mathcal{E}_{k,n}^{\ell}$; $\ \mathcal{F}_{k,j}^{\ell'}\ (1\leq k\le
j< n),\ \mathcal{F}_{k,j'}^{\ell'}\ (1\leq k<j\leq n),\, f_n^\ell,\
\mathcal{F}_{k,n}^{\ell}$; and $\omega_k^{\ell'}-1,\,
(\omega_k')^{\ell'}-1\ (1\leq k< n),\ \omega_n^{\ell}-1,\,
(\omega_n')^{\ell}-1$ are central in
$U_{r,s}(\mathfrak{so}_{2n+1})$. 

\section{Restricted two-parameter quantum groups}

\noindent{\it 4.1.} In what follows, {\it we assume that $\ell$ is
odd.}
\begin{defi}\label{3.17}\  The
\textit{restricted two-parameter quantum group} is the quotient
$$\mathfrak{u}_{r,s}(\mathfrak{so}_{2n+1}):=U_{r,s}(\mathfrak{so}_{2n+1})/I_n,$$
where  $I_n$ is the ideal of $U_{r,s}(\mathfrak{so}_{2n+1})$
generated by  $\mathcal{E}_{k,j}^{\ell}\ (1\leq k\leq j\leq n),\
\mathcal{E}_{k,j'}^{\ell}\ (1\leq k<j\leq n),\
\mathcal{F}_{k,j}^{\ell}\ (1\leq k\leq j\leq n),\
\mathcal{F}_{k,j'}^{\ell}\ (1\leq k<j\leq n)$ and
$\omega_k^{\ell}-1,\ (\omega_k')^{\ell}-1 \ (1\leq k\leq n)$.
\end{defi}

By Theorem 2.5 and Corollary 2.7,
$\mathfrak{u}_{r,s}(\mathfrak{so}_{2n+1})$ is an algebra of
dimension $\ell^{2n^2+2n}$ with linear basis
$$
\mathcal{E}_{1,1}^{a_1}\mathcal{E}_{1,2}^{a_2}\cdots
\mathcal{E}_{n-1,n'}^{a_{n^2-1}}\mathcal{E}_{n,n}^{a_{n^2}}
\omega_1^{b_1}\cdots\omega_{n}^{b_{n}}
(\omega'_1)^{b'_1}\cdots(\omega'_{n})^{b'_{n}}
\mathcal{F}_{1,1}^{a'_1}\mathcal{F}_{1,2}^{a'_2}\cdots\mathcal{F}_{n-1,n'}^{a'_{n^2-1}}
\mathcal{F}_{n,n}^{a'_{n^2}},
$$
where all powers range between $0$ and $\ell-1$.

\medskip \noindent{\it 4.2.} The main theorem of this section is
\begin{theorem}
The ideal $I_n$ is a Hopf ideal, so that
$\mathfrak{u}_{r,s}(\mathfrak{so}_{2n+1})$ is a finite dimensional
Hopf algebra.
\end{theorem}

\noindent{\it 4.3.} The proof of Theorem 4.2 will be done through a
sequence of Lemmas. To determine $\Delta
({\mathcal{E}}_{i,j}^{\ell})$, $\Delta(\mathcal{E}_{k,j'}^{\ell})$,
recall $j'=2n-j+1$ given in Section 2. We make conventions:
\begin{eqnarray}
\omega_{i,j}=\prod_{k=i}^j\omega_k\; (i\le j),
\quad\omega_{i,j'}=\omega_{i,n}\omega_n\cdots\omega_j\; (i<j),\quad \zeta=s^2{-}r^2;\\
\mathcal{E}_{n+1,n+1}=e_n,\quad \mathcal{E}_{n,n+1}=(s^2{-}rs)e_n^2;
\\
\mathcal{E}_{2n-k,2n-k}=r^{-2}s^{-2}e_{k+1}, \quad 1\leq k<n-1;\\
\mathcal{E}_{k,2n-k+1}=s^2\mathcal{E}_{k,2n-k}e_k-s^{-2}e_k\mathcal{E}_{k,2n-k},\
1\leq k<n; \\
\mathcal{E}_{2n-i+1,2n-k+1}=\mathcal{E}_{2n-i+1,2n-k}e_k-s^{-2}e_k\mathcal{E}_{2n-i+1,2n-k},
\  n \geq i>k; \\
\qquad
\mathcal{E}_{k+1,2n-k+1}=r^2\mathcal{E}_{k+1,2n-k}e_k
-s^{-2}e_k\mathcal{E}_{k+1,2n-k}+\mathcal{E}_{k,2n-k},
\ 1\leq k<n;\\
\mathcal{E}_{i,2n-k+1}=\mathcal{E}_{i,2n-k}e_k-s^{-2}e_k\mathcal{E}_{i,2n-k},
\qquad n\geq i> k+1.
\end{eqnarray}

It is easy to verify that
\begin{equation*}
\begin{split}
\mathcal{E}_{k,2n-k+1}&=\sum_{i=0}^{n-k-1}(-1)^ir^{-2(i+1)}s^{-2}\zeta
\mathcal{E}_{k,k+i}\mathcal{E}_{k,2n-k-i}\\
&\qquad+(-1)^{n-k}r^{-2(n-k)}(s^2-rs)\mathcal{E}_{k,n}^2 \in U^+,\\
\mathcal{E}_{2n-i+1,2n-k+1}&=\sum_{j=0}^{i-k}(-1)^jr^{-2j}\zeta^{i-k-j}
\Bigl(\sum_b\mathcal{E}_{I}\Bigr)\in U^+,
\end{split}
\end{equation*}
where $ b=i-k+1-j$ and
$\mathcal{E}_{I}=\mathcal{E}_{i_0,i_1}\mathcal{E}_{i_1+1,i_2}\mathcal{E}_{i_2+1,i_3}\cdots
\mathcal{E}_{i_{b-1}+1,i_b}$ for $k\leq i_0\leq i_1\leq
i_2-1\leq\cdots \leq i_b-1\leq i-1$.

A simple calculation shows that $\mathcal{E}_{k+1,2n-k+1}\in U^+$
for  $1\leq k<n$, and $\mathcal{E}_{i,2n-k+1}\in U^+$ for $k+1<i\le
n$. So the elements in (4.1)---(4.7) are in $U^+$.

\begin{lemm}
$\mathrm{(i)}\ \Delta (\mathcal{E}_{k,j})=\mathcal{E}_{k,j}\otimes 1
+ \omega_{k,j}\otimes \mathcal{E}_{k,j}+ \zeta
\sum\limits_{i=k}^{j-1}\mathcal{E}_{i+1,j}\omega_{k,i} \otimes
\mathcal{E}_{k,i}$ for $1\leq k\leq j\leq n;$

$\mathrm{(ii)}\ \Delta
(\mathcal{E}_{k,j'})=\mathcal{E}_{k,j'}\otimes 1 +
\omega_{k,j'}\otimes \mathcal{E}_{k,j'}+\zeta
\sum\limits_{i=k}^{2n-j}\mathcal{E}_{i+1,2n-j+1 }\omega_{k,i}
\otimes \mathcal{E}_{k,i}$ for $1\leq k<j\leq n$.
\end{lemm}
\begin{proof}
(i): The statement is true when $k=j$, owing to the coproduct of
$e_k$. Assume that it is true for $k+1$. Then we obtain
\begin{equation*}
\begin{split}
\Delta &(\mathcal{E}_{k,j})=\Delta (e_{k})\Delta
({\mathcal{E}}_{k+1,j})-r^2\Delta (\mathcal{E}_{k+1,j})\Delta
(e_{k})\\
&=(e_k{\otimes} 1{+}\omega_k{\otimes} e_k)
\Bigl({\mathcal{E}}_{k+1,j}{\otimes} 1 {+} \omega_{k+1,j}{\otimes}
{\mathcal{E}}_{k+1,j}{+}\zeta
\sum_{i=k+1}^{j-1}\mathcal{E}_{i+1,j}\omega_{k+1,i}
{\otimes} \mathcal{E}_{k+1,i}\Bigr) \\
&\quad{-}r^2\Bigl(\mathcal{E}_{k+1,j}{\otimes} 1 {+}
\omega_{k+1,j}{\otimes} \mathcal{E}_{k+1,j}{+}{\zeta}
\sum_{i=k+1}^{j-1}\mathcal{E}_{i+1,j}\omega_{k+1,i} {\otimes}
\mathcal{E}_{k+1,i}\Bigr)(e_k{\otimes} 1{+}\omega_k{\otimes} e_k) \\
&= e_k\mathcal{E}_{k+1,j}\otimes
1+\omega_k\mathcal{E}_{k+1,j}\otimes
e_k+e_k\omega_{k+1,j}\otimes\mathcal{E}_{k+1,j}+\omega_{k,j}\otimes
e_k\mathcal{E}_{k+1,j}\\
&
\quad+{\zeta}\sum_{i=k+1}^{j-1}e_k\mathcal{E}_{i+1,j}\omega_{k+1,i}
\otimes \mathcal{E}_{k+1,i} + \zeta
\sum_{i=k+1}^{j-1}\omega_k\mathcal{E}_{i+1,j}\omega_{k+1,i} \otimes
e_k\mathcal{E}_{k+1,i}\\
&\quad -r^2\mathcal{E}_{k+1,j}e_k\otimes
1-r^2\mathcal{E}_{k+1,j}\omega_k\otimes
e_k-r^2\omega_{k+1,j}e_k\otimes\mathcal{E}_{k+1,j}\\
&\quad-r^2\omega_{k+1,j}\omega_k\otimes \mathcal{E}_{k+1,j}e_k
-r^2\zeta\sum_{i=k+1}^{j-1}\mathcal{E}_{i+1,j}\omega_{k+1,i}e_k
\otimes \mathcal{E}_{k+1,i}\\
& \quad-r^2 \zeta
\sum_{i=k+1}^{j-1}\mathcal{E}_{i+1,j}\omega_{k+1,i}\omega_k \otimes
\mathcal{E}_{k+1,i}e_k \\
&=\mathcal{E}_{k,j}\otimes 1 + \omega_{k,j}\otimes
\mathcal{E}_{k,j}+\zeta \mathcal{E}_{k+1,j}\omega_k\otimes e_k+\zeta
\sum_{i=k+1}^{j-1}\mathcal{E}_{i+1,j}\omega_{k,i} \otimes
\mathcal{E}_{k,i}\\&=\mathcal{E}_{k,j}\otimes 1 +
\omega_{k,j}\otimes \mathcal{E}_{k,j}+\zeta
\sum_{i=k}^{j-1}\mathcal{E}_{i+1,j}\omega_{k,i} \otimes
\mathcal{E}_{k,i}
\end{split}
\end{equation*}
by using induction hypothesis.

(ii): Use induction on $j$. For $j=n$, the argument proceeds by
induction on $k$. We have
\begin{equation*}
\begin{split}
\Delta (\mathcal{E}_{k,n'})&=\mathcal{E}_{k,n'}\otimes 1 +
\omega_{k,n'}\otimes \mathcal{E}_{k,n'}+\zeta
\sum_{i=k}^{n-2}\mathcal{E}_{i+1,n' }\omega_{k,i}\otimes
\mathcal{E}_{k,i}\\ & \quad+(s^2-rs)\zeta e_n^2\omega_{k,n-1}
\otimes \mathcal{E}_{k,n-1}+\zeta e_n \omega_{k,n}\otimes
\mathcal{E}_{k,n}\\
&=\mathcal{E}_{k,n'}\otimes 1 + \omega_{k,n'}\otimes
\mathcal{E}_{k,n'}+\zeta \sum_{i=k}^{n}\mathcal{E}_{i+1,n+1
}\omega_{k,i} \otimes \mathcal{E}_{k,i}.
\end{split}
\end{equation*}
Assume that it is true for $j+1$, we obtain
\begin{equation*}
\begin{split}
 \Delta(\mathcal{E}_{k,j'})&=\Bigl(\mathcal{E}_{k,j+1'}\otimes 1 + \omega_{k,j+1'}\otimes
\mathcal{E}_{k,j+1'}+\zeta \sum_{i=k}^{2n-j-1}\mathcal{E}_{i+1,2n-j
}\omega_{k,i} \otimes \mathcal{E}_{k,i}\Bigr)\\
&\quad \times(e_j\otimes 1+\omega_j\otimes e_j) -s^{-2}(e_j\otimes
1+\omega_j\otimes e_j) \Bigl(\mathcal{E}_{k,j+1'}\otimes 1\\
&\quad + \omega_{k,j+1'}\otimes \mathcal{E}_{k,j+1'}+\zeta
\sum_{i=k}^{2n-j-1}\mathcal{E}_{i+1,2n-j }\omega_{k,i} \otimes
\mathcal{E}_{k,i}\Bigr)\\
\end{split}
\end{equation*}
\begin{equation*}
\begin{split}&=\mathcal{E}_{k,j+1'}e_j\otimes 1 + \omega_{k,j+1'}e_j\otimes
\mathcal{E}_{k,j+1'} +\zeta \sum_{i=k}^{2n-j-1}\mathcal{E}_{i+1,2n-j
}\omega_{k,i}e_j \otimes \mathcal{E}_{k,i}\\
&\quad+\mathcal{E}_{k,j+1'}\omega_j\otimes e_j +
\omega_{k,j+1'}\omega_j\otimes \mathcal{E}_{k,j+1'}e_j\\
&\quad+\zeta \sum_{i=k}^{2n-j-1}\mathcal{E}_{i+1,2n-j
}\omega_{k,i}\omega_j \otimes \mathcal{E}_{k,i}e_j\\
&\quad -s^{-2}e_j\mathcal{E}_{k,j+1'}\otimes 1 -s^{-2}e_j
\omega_{k,j+1'}\otimes \mathcal{E}_{k,j+1'}\\
&\quad-s^{-2}\zeta \sum_{i=k}^{2n-j-1}e_j\mathcal{E}_{i+1,2n-j
}\omega_{k,i} \otimes
\mathcal{E}_{k,i}\\
&\quad
-s^{-2}\omega_j\mathcal{E}_{k,j+1'}\otimes e_j -s^{-2}\omega_j
\omega_{k,j+1'}\otimes e_j\mathcal{E}_{k,j+1'}\\
&\quad-s^{-2}\zeta \sum_{i=k}^{2n-j-1}\omega_j\mathcal{E}_{i+1,2n-j
}\omega_{k,i}
\otimes e_j\mathcal{E}_{k,i}\\
&= \mathcal{E}_{k,j'}\otimes 1 + \omega_{k,j'}\otimes
\mathcal{E}_{k,j'}+\zeta \mathcal{E}_{2n-j+1,2n-j+1 }\omega_{k,2n-j}
\otimes
\mathcal{E}_{k,2n-j}\\
&\quad+\zeta \sum_{i=k}^{j-2}\mathcal{E}_{i+1,2n-j+1 }\omega_{k,i}
\otimes \mathcal{E}_{k,i}+\zeta \mathcal{E}_{j,2n-j+1
}\omega_{k,j-1} \otimes \mathcal{E}_{k,j-1}\\
&\quad+\zeta \mathcal{E}_{j+1,2n-j+1 }\omega_{k,j} \otimes
\mathcal{E}_{k,j}+\zeta \sum_{i=j+1}^{n}\mathcal{E}_{i+1,2n-j+1
}\omega_{k,i} \otimes \mathcal{E}_{k,i}\\
&\quad+\zeta \sum_{i=n+1}^{2n-j-1}\mathcal{E}_{i+1,2n-j+1
}\omega_{k,i} \otimes
\mathcal{E}_{k,i}\\
&= \mathcal{E}_{k,j'}\otimes 1 + \omega_{k,j'}\otimes
\mathcal{E}_{k,j'}+\zeta \sum_{i=k}^{2n-j}\mathcal{E}_{i+1,2n-j+1
}\omega_{k,i} \otimes \mathcal{E}_{k,i}
\end{split}
\end{equation*}
by induction hypothesis and notations in (4.1)---(4.7) and Lemma 3.2
(3).
\end{proof}

\noindent{\it 4.4.} We generalize Lemma 4.3 (i) to an expression for
$\Delta (\mathcal{E}_{k,j}^a)$ in Lemma \ref{3.21} below. The
formula involves certain exponents $p_m$ of $s$ that are defined
recursively on ordered $m$-tuples of nonnegative integers, as
follows (see \cite{BW3}):
\begin{equation*}
\begin{split}
p_m(0,\cdots,0)&:=0, \\
p_m(c_1+1,c_2,\cdots,c_m)&:=p_m(c_1,c_2,\cdots,c_m)-c_2-c_3-\cdots -c_m, \\
p_m(c_1,c_2,\cdots,c_{m-1},c_{m}+1)&:=p_m(c_1,c_2,\cdots,c_m)-c_1-c_2- \cdots -c_{m-1}, \\
p_m(c_1,c_2,\cdots,c_j+1,\cdots,c_m)&:=p_m(c_1,c_2,\cdots,c_m)-c_1-c_2-\cdots-c_{j-1}\\
&\quad +c_j-c_{j+1}-\cdots-c_m,\quad (1<j<m).
\end{split}
\end{equation*}
An  inductive argument on $c_1+\cdots +c_m$ shows that $p_m$ is
well-defined.

Set $C^a_m(r,s)=\frac{[c_1+\cdots+c_m]!}{[c_1]!\cdots[c_m]!}$, for
$a=c_1+\cdots+c_m$.
The lemma below is true, which is equivalent to checking
inductively: $
[\,c_1+\cdots+c_m\,]=\sum_{j=1}^mr^{2C_j}s^{2D_j}[\,c_j\,]$.
\begin{lemm}\label{3.20}
\ Let $c_1,\cdots,c_m$ be positive integers with $a=c_1+\cdots+c_m$. Then 
$$C^a_m(r,s)
=\sum\limits_{j=1}^{m}r^{2C_j}s^{2D_j}\left[c_1+\cdots+c_m-1 \atop
c_m\right]\cdots\left[c_1+\cdots+c_{j}-1 \atop
c_{j}-1\right]\cdots\left[c_1+c_{2} \atop c_{2}\right],
$$
where $C_j=c_{j+1}+c_{j+2}+\cdots+c_{m}, D_j=c_1+\cdots+c_{j-1}$,
and $C_m=0=D_1$.
\end{lemm}

In the above sum, the binomial coefficients have $1$ subtracted from
their top number to the one with $c_1+\cdots+c_j$ on top and not
thereafter, as in the $j$th summand we have replaced $c_j$ by
$c_j-1$.

\begin{lemm}\label{3.21}\  For $1\le k\leq j\le n$ and $m=j-k+2$. Then
\begin{equation*}
\begin{split}
& \Delta (\mathcal{E}_{k,j}^{a})=\sum
s^{2p_m(c_1,\cdots,c_m)}\zeta^{a-c_1-c_m}C^a_m(r,s)\mathcal{E}_{k,j}^{c_1}
\mathcal{E}_{k+1,j}^{c_2} \cdots
\mathcal{E}_{j,j}^{c_{m-1}}\omega_{k,k}^{c_2}\omega_{k,k+1}^{c_3}
\cdots \omega_{k,j}^{c_{m}} \\ & \qquad\qquad\qquad \otimes
\mathcal{E}_{k,k}^{c_2}\mathcal{E}_{k,k+1}^{c_3} \cdots
\mathcal{E}_{k,j}^{c_{m}},
\end{split}
\end{equation*}
the sum taken over all $m$-tuples $(c_1,\cdots,c_{m})$ with
$c_1+\cdots+c_{m}=a$.
\end{lemm}
\begin{proof}
Observe first that if $a=1$, the  above expression coincides with
that in Lemma 4.3 (i). To simplify notation, we will assume without
loss of generality that $k=1$ (and hence $m=j+1$). Now suppose that
the above formula holds for $a-1$, that is $\Delta
(\mathcal{E}_{1,j}^{a-1})=\sum
s^{2p_{j+1}(c_1,\cdots,c_{j+1})}\zeta^{a-1-c_1-c_{j+1}}C^{a-1}_{j+1}(r,s)
\mathcal{E}_{1,j}^{c_1}
\cdots
\mathcal{E}_{j,j}^{c_{j}}\omega_{1,1}^{c_2}
\cdots \omega_{1,j}^{c_{j+1}}$ $ \otimes
\mathcal{E}_{1,1}^{c_2}
\cdots \mathcal{E}_{1,j}^{c_{j+1}}$, where the sum is over all
$(j{+}1)$-tuples $(c_1,\cdots,c_{j+1})$ with
$c_1+\cdots+c_{j+1}=a-1$. Then
$\Delta(\mathcal{E}_{1,j}^{a})=\Delta(\mathcal{E}_{1,j}^{a-1})\Delta(\mathcal{E}_{1,j})$,
which by Lemma 4.3 (i) is equal to the above sum times
$$
\mathcal{E}_{1,j}\otimes 1
+ \omega_{1,j}\otimes \mathcal{E}_{1,j}+ \zeta
\sum_{i=1}^{j-1}\mathcal{E}_{i+1,j}\omega_{1,i} \otimes
\mathcal{E}_{1,i}.
$$
Expanding and re-summing over all $(d_1,\cdots,d_{j+1})$ with
$d_1+\cdots+d_{j+1}=a$, by Lemma 3.3 (4) \& (5), we obtain
\begin{equation*}
\begin{split}
&\text{\it the coefficient of }\;
\mathcal{E}_{1,j}^{d_1}\mathcal{E}_{2,j}^{d_2} \cdots
\mathcal{E}_{j,j}^{d_{j}}\omega_{1,1}^{d_2}\omega_{1,2}^{d_3} \cdots
\omega_{1,j}^{d_{j+1}}\otimes
\mathcal{E}_{1,1}^{d_2}\mathcal{E}_{1,2}^{d_3} \cdots
\mathcal{E}_{1,j}^{d_{j+1}}\\
&\quad=s^{2p_{j+1}(d_1,\cdots,d_{j+1})}\zeta^{a-d_1-d_{j+1}}\times\\
&\qquad\times\Bigl(\sum\limits_{i=1}^{j+1}r^{2C'_i}s^{2D'_i}\left[d_1+\cdots+d_{j+1}-1
\atop d_{j+1}\right]\cdots\left[d_1+\cdots+d_{i}-1 \atop
d_{i}-1\right]\cdots\left[d_1+d_{2} \atop d_{2}\right]\Bigr),
\end{split}
\end{equation*}
 where $C'_i=d_{i+1}+d_{i+2}+\cdots+d_{j+1},$
 $D'_i=d_1+\cdots+d_{i-1}$ and $C'_{j+1}=0=D'_1$. By Lemma \ref{3.20}, we
 get the desired result.
\end{proof}
As a result, we have the following formula for
$e_j=\mathcal{E}_{j,j}$ for $1\leq j\leq n$:
$$\Delta (e_{j}^a)= \sum_{i=0}^{a}s^{2i(i-a)}\left[a\atop i\right]
e_{j}^i\omega_{j}^{a-i} \otimes e_{j}^{a-i}.\leqno (4.8)$$

\smallskip
\noindent{\it 4.5.} We generalize Lemma 4.5 to a formula for $\Delta
(\mathcal{E}_{k,n'}^{\ell})$ in Lemma 4.6 below. The formula
involves certain exponents $p_m (m=n-k+3)$ of $s$ that are defined
recursively on ordered $m$-tuples of nonnegative integers, as
follows (which are different from those in Lemma \ref{3.21}):
\begin{equation*}
\begin{split}
p_m(0,\cdots,0)&:=0, \\
p_m(c_1+1,c_2,\cdots,c_m)&:=p_m(c_1,c_2,\cdots,c_m)-c_2-c_3-\cdots -c_m, \\
p_m(c_1,c_2,\cdots,c_{m-1},c_{m}+1)&:=p_m(c_1,c_2,\cdots,c_m)-c_1-c_2- \cdots -c_{m-1}, \\
\end{split}
\end{equation*}
\begin{equation*}
\begin{split}p_m(c_1,c_2,\cdots,c_j+1,\cdots,c_m)&:=p_m(c_1,c_2,\cdots,c_m)-c_1-c_2-\cdots-c_{j-1}\\
&\quad +c_j-c_{j+1}-\cdots-c_m,\quad (1<j<m-2),
\\p_m(c_1,c_2,\cdots,c_{m-2}+1,c_{m-1},c_m)&:=p_m(c_1,c_2,\cdots,c_m)-c_1-\cdots-c_{m-3}\\
&\quad +2c_{m-2}-c_m,\\
p_m(c_1,c_2,\cdots,c_{m-2},c_{m-1}+1,c_m)&:=p_m(c_1,c_2,\cdots,c_m)-c_1-\cdots-c_{m-3}
-c_m.\\
\end{split}
\end{equation*}
An inductive argument on $c_1+\cdots +c_m$ shows that $p_m$ is
well-defined.

\begin{lemm}\label{3.22}\  For $1\leq k< n$ and $m=n-k+3$,  we have
\begin{equation*}
\begin{split}
\Delta (\mathcal{E}_{k,n'}^{a})&=\sum\limits_{c_{m-1}=0,1}
s^{2p_m(c_1,\cdots,c_m)}\zeta^{a-c_1-c_m}C^a_m(r,s)\mathbf{E}^{(a)}_k\\
&\quad +\zeta(rs-r^2)[a]!\sum\limits_{2\leq c_{m-1}< a}
s^{2p_m(c_1,\cdots,c_m)}\mathbf{E}^{(a)}_k\\
&\quad +
r^{\frac{a(a-1)}{2}}(s-r)^{a-1}\zeta\Big(\prod_{i=2}^{a}(r^i+s^i)\Big)
e_n^a\omega_{k,n}^a\otimes\mathcal{E}_{k,n}^a,
\end{split}
\end{equation*}
where  $\mathbf{E}^{(a)}_k=
\mathcal{E}_{k,n+1}^{c_1}
\cdots
\mathcal{E}_{n+1,n+1}^{c_{n-k+2}}
\omega_{k,k}^{c_2}
\cdots \omega_{k,k+j-2}^{c_j}
\cdots
\omega_{k,n+1}^{c_{n-k+3}}\otimes
\mathcal{E}_{k,k}^{c_2}
\cdots\mathcal{E}_{k,k+j-2}^{c_j}\cdots\times$ 
$\times\, \mathcal{E}_{k,n+1}^{c_{n-k+3}}$ with
$c_1+\cdots+c_{m}=a$. 
\end{lemm}
\begin{proof}
Observe first that if $a=1$, the  above expression only has the
first term which coincides with that in Lemma 4.3 (ii). To simplify
notation, we will assume without loss of generality that $k=1$ (and
hence $m=n+2$). Now suppose that the above formula holds for $a-1$,
that is,
\begin{equation*}
\begin{split}
\Delta (\mathcal{E}_{1,n'}^{a-1})&=\sum\limits_{c_{n+1}=0,1\atop
c_1+\cdots+c_{n+2}=a-1}
s^{2p_{n+2}(c_1,\cdots,c_{n+2})}\zeta^{a-1-c_1-c_{n+2}}C^{a-1}_{n+2}(r,s)\mathbf{E}^{(a-1)}_1\\
&\quad+\zeta(rs-r^2)[a-1]!\sum\limits_{2\leq c_{n+1}< a-1,\atop
c_1+\cdots+c_{n+2}=a-1}
s^{2p_{n+2}(c_1,\cdots,c_{n+2})}\mathbf{E}^{(a-1)}_1\\
&\quad+
r^{\frac{(a-2)(a-1)}{2}}(s-r)^{a-2}\zeta\Big(\prod_{i=2}^{a-1}(r^i+s^i)\Big)e_n^{a-1}
\omega_{k,n}^{a-1}
\otimes\mathcal{E}_{k,n}^{a-1},
\end{split}
\end{equation*}
where the sum is over all $(n+2)$-tuples $(c_1,\cdots,c_{n+2})$ with
$c_1+\cdots+c_{n+2}=a-1$. Then
$\Delta(\mathcal{E}_{1,n'}^{a})=\Delta(\mathcal{E}_{1,n'}^{a-1})\Delta(\mathcal{E}_{1,n'})$,
which by Lemma 4.3 (ii) is equal to the above sum times
$$
\mathcal{E}_{1,n'}\otimes 1 + \omega_{1,n'}\otimes
\mathcal{E}_{1,n'}+\zeta \sum_{i=1}^{n}\mathcal{E}_{i+1,n'
}\omega_{1,i} \otimes \mathcal{E}_{1,i}.
$$
Expanding and re-summing over all $(d_1,\cdots,d_{n+2})$ with
$d_1+\cdots+d_{n+2}=a$, we get the desired result.
\end{proof}

\noindent {\it 4.6. Proof of Theorem 4.2.} \  We are now in the
position to prove that $I_n$ is a Hopf ideal. As $[\ell]=0$, it
follows from Lemma \ref{3.21} that the only nonzero terms of
$\Delta(\mathcal{E}_{k,j}^\ell)$ are those having $c_i=\ell$ for
some $i$. As a result,
$$\Delta(\mathcal{E}_{k,j}^\ell)=\mathcal{E}_{k,j}^\ell\otimes 1+
\omega_{k,j}^\ell\otimes\mathcal{E}_{k,j}^\ell+s^{\ell(\ell-1)}\zeta^\ell\sum_{i=k}^{j-1}
\mathcal{E}_{i+1,j}^\ell\omega_{k,i}^\ell\otimes\mathcal{E}_{k,i}^\ell,
\leqno (4.9)$$ which is clearly in $I_n\otimes U+U\otimes I_n$.
Applying the antipode property to (4.9), we obtain
$$
0=\varepsilon(\mathcal{E}_{k,j}^\ell)=\mathcal{E}_{k,j}^\ell+
\omega_{k,j}^\ell
S(\mathcal{E}_{k,j}^\ell)+s^{\ell(\ell-1)}\zeta^\ell\sum_{i=k}^{j-1}
\mathcal{E}_{i+1,j}^\ell\omega_{k,i}^\ell
S(\mathcal{E}_{k,i}^\ell).
$$ So we have
$$
S(\mathcal{E}_{k,j}^\ell)=-\omega_{k,j}^{-\ell}
\bigl(\mathcal{E}_{k,j}^\ell+s^{\ell(\ell-1)}\zeta^\ell\sum_{i=k}^{j-1}
\mathcal{E}_{i+1,j}^\ell\omega_{k,i}^\ell S(\mathcal{E}_{k,i}^\ell)
\bigr).\leqno (4.10)$$ If $j=k$, we have
$S(\mathcal{E}_{k,j}^\ell)=S(e_k^\ell)=(-\omega_k^{-1}e_k)^\ell$,
which is a scalar multiple of $\omega_k^{-\ell}e_k^\ell \in I_n$.
Using (4.10) and induction on $j-k$, we see that
$S(\mathcal{E}_{k,j}^\ell)\in I_n$ for all $k<j$ as well. From Lemma
\ref{3.22}, we have
\begin{equation*}
\begin{split}
\Delta(\mathcal{E}_{k,n'}^\ell)
&=\mathcal{E}_{k,n'}^\ell\otimes 1+
\omega_{k,n'}^\ell\otimes\mathcal{E}_{k,n'}^\ell+s^{2\ell(\ell-1)}\zeta^\ell
\mathcal{E}_{n,n'}^\ell\omega_{k,n-1}^\ell\otimes\mathcal{E}_{k,n-1}^\ell
\\
&\quad
+r^{\frac{\ell(\ell-1)}{2}}(s-r)^{\ell-1}\zeta\Big(\prod_{i=2}^{\ell}(r^i+s^i)\Big)e_n^\ell
\omega_{k,n}^\ell\otimes\mathcal{E}_{k,n}^\ell\\
&\quad
+ s^{\ell(\ell-1)}\zeta^\ell\sum_{i=k}^{n-2}
\mathcal{E}_{i+1,n'}^\ell\omega_{k,i}^\ell\otimes\mathcal{E}_{k,i}^\ell,
\end{split}\tag{4.11}
\end{equation*}
which is clearly in $I_n\otimes U+U\otimes I_n$. Applying the
antipode property to (4.11), we have $S(\mathcal{E}_{k,n'}^\ell) \in
I_n$. For $j<n$, using a similar method of Lemma \ref{3.22}, we have
\begin{gather*}
\begin{split}
\Delta(\mathcal{E}_{k,j'}^\ell)&=\mathcal{E}_{k,j'}^\ell\otimes 1+
\omega_{k,j'}^\ell\otimes\mathcal{E}_{k,j'}^\ell+s^{2\ell(\ell-1)}\zeta^\ell
\mathcal{E}_{j,j'}^\ell\omega_{k,j-1}^\ell\otimes\mathcal{E}_{k,j-1}^\ell
\\
&\quad + r^{\frac{\ell(\ell-1)}{2}}(s-r)^{\ell-1}\zeta\Big(
\prod_{i=2}^{\ell}(r^i+s^i)\Big)\mathcal{E}_{n+1,j'}^\ell
\omega_{k,n}^\ell\otimes\mathcal{E}_{k,n}^\ell\\
&\quad + s^{\ell(\ell-1)}\zeta^\ell\sum_{i=k}^{j-2}
\mathcal{E}_{i+1,j'}^\ell\omega_{k,i}^\ell\otimes\mathcal{E}_{k,i}^\ell\\
&\quad + r^{\ell(\ell-1)}s^{2\ell(\ell-1)}\zeta^\ell\sum_{i=j}^{n-1}
\mathcal{E}_{i+1,j'}^\ell\omega_{k,i}^\ell\otimes\mathcal{E}_{k,i}^\ell\\
&\quad
+ r^{-\ell(\ell-1)}\zeta^\ell\sum_{i=n+1}^{2n-j}
\mathcal{E}_{i+1,j'}^\ell\omega_{k,i}^\ell\otimes\mathcal{E}_{k,i}^\ell,
\end{split}\tag{4.12}
\end{gather*}
which is in $I_n\otimes U+U\otimes I_n$. Applying the antipode
property to this, we have $S(\mathcal{E}_{k,j'}^\ell) \in I_n$.

Using $\tau$, we find that $\Delta(\mathcal{F}_{k,j}),\
\Delta(\mathcal{F}_{k,j'}) \in I_n\otimes U+U\otimes I_n$. So $I_n$
is a Hopf ideal, and $\mathfrak{u}_{r,s}(\mathfrak{so}_{2n+1})$ is a
finite-dimensional Hopf algebra. \hfill\qed

\section{Isomorphisms of $\mathfrak
u_{r,s}(\mathfrak {so}_{2n+1})$, as well as $\mathfrak
u_{r,s}(\mathfrak {sl}_n)$}

\smallskip
\noindent{\it 5.1.} Write
$\mathfrak{u}_{r,s}=\mathfrak{u}_{r,s}(\mathfrak{so}_{2n+1})$. Let
$G$ denote the group generated by $\omega_i$, $\omega'_i\ (1\leq
i\leq n)$ in $\mathfrak{u}=\mathfrak{u}_{r,s}$. Define linear
subspace $\mathfrak{a}_k$ of $\mathfrak{u}$ by
\begin{gather*}
\begin{split}
\mathfrak{a}_0=\mathbb{K}G, \qquad \mathfrak{a}_1=\mathbb{K}G
+\sum\limits_{i=1}^{n}(\mathbb{K}e_iG+\mathbb{K}f_iG), \\
\mathfrak{a}_k=(\mathfrak{a}_1)^k \quad \textrm{for}\quad  k\geq 1.
\end{split}\tag{5.1}
\end{gather*}
Note that $1 \in \mathfrak{a}_0$, $\Delta(\mathfrak{a}_0)\subseteq
\mathfrak{a}_0\otimes \mathfrak{a}_0$, $\mathfrak{a}_1$ generates
$\mathfrak{u}$ as an algebra, and $\Delta(\mathfrak{a}_1)\subseteq
\mathfrak{a}_1\otimes \mathfrak{a}_0+\mathfrak{a}_0\otimes
\mathfrak{a}_1$. By \cite{M}, $\{\mathfrak{a}_k\}$ is a coalgebra
filtration of $\mathfrak{u}$ and
$\mathfrak{u}_0\subseteq\mathfrak{a}_0$, where the coradical
$\mathfrak{u}_0$ of $\mathfrak{u}$ is the sum of all the simple
subcoalgebras of $\mathfrak{u}$. Clearly, $\mathfrak{a}_0\subseteq
\mathfrak{u}_0$ as well, and so $\mathfrak{u}_0=\mathbb{K}G$. This
implies that $ \mathfrak{u}$ is \textit{pointed}, that is, every
simple subcoalgebra of $\mathfrak{u}$ is one-dimensional.

Let $\mathfrak{b}$ (resp., $\mathfrak{b}'$) be the Hopf subalgebra
of $\mathfrak{u}=\mathfrak{u}_{r,s}(\mathfrak{so}_{2n+1})$ generated
by $e_i, \omega_i^{\pm 1}$ (resp., $f_i, (\omega_i')^{\pm 1}$),
$1\leq i\leq n$. The same argument shows that $\mathfrak{b}$ and
$\mathfrak{b}'$ are pointed as well. Then we have
\begin{prop}\label{4.1}
$\mathfrak{u}_{r,s}(\mathfrak{so}_{2n+1})$ is a pointed Hopf
algebra, as are  $\mathfrak{b}$ and $\mathfrak{b}'$. \hfill\qed
\end{prop}

\noindent{\it 5.2.} It follows from [{\bf M}, Lemma 5.5.1] that
$\mathfrak{a}_k\subseteq \mathfrak{u}_k$ for all $k$, where
$\{\mathfrak{u}_k\}$ is the \textit{coradical filtration} of
$\mathfrak{u}$ defined inductively by
$\mathfrak{u}_k=\Delta^{-1}(\mathfrak{u}\otimes
\mathfrak{u}_{k-1}+\mathfrak{u}_0\otimes \mathfrak{u}).$ In
particular, $\mathfrak{a}_1\subseteq \mathfrak{u}_1.$ By [{\bf M},
Theorem 5.4.1], as $\mathfrak{u}$ is pointed, $\mathfrak{u}_1$ is
spanned by the set of group-like elements $G$, together with all the
skew-primitive elements of $\mathfrak{u}$.

Let $m$ be the smallest positive integer such that $r^{m}=s^{m}$.
Note that $[m]=0$ while $[1], [2],\cdots, [m-1]$ are all nonzero. It
follows from this observation and $(4.8)$ that $e_i^m, f_i^m \
(1\leq k\leq n)$ are skew-primitive if $m<\ell$. From now on, we
assume that $m=\ell$. Based on Lemmas 4.3, 4.5, 4.6 \& formula
(4.12), together with with the same argument of [BW3, Lemma 3.3], we
have

\begin{lemm}\label{4.2}
Assume that $rs^{-1}$ is a primitive $\ell$th root of unity. Then
$${\mathfrak u}_1=\mathbb{K}G
+\sum\limits_{i=1}^{n}(\mathbb{K}e_iG+\mathbb{K}f_iG).$$
\end{lemm}

\smallskip
\noindent{\it 5.3.} Given two group-like elements $g$ and $h$ in a
Hopf algebra $H$, let $P_{g,h}(H)$ denote the set of skew-primitive
elements of $H$ given by
$$
P_{g,h}(H)=\{x\in H\mid \Delta(x)=x\otimes g+h\otimes x\}.
$$

\begin{lemm} Assume that $rs^{-1}$ is a primitive $\ell$th root of
unity. Then
\smallskip

$({\textrm{\rm i}})$ \;\
$P_{1,\omega_i}(\mathfrak{u}_{r,s})=\mathbb{K}(1-\omega_i)+\mathbb{K}e_i;$
$\quad
P_{1,\omega_i'^{-1}}(\mathfrak{u}_{r,s})=\mathbb{K}(1-\omega_i'^{-1})+\mathbb{K}f_i\omega_i'^{-1};$

$\qquad P_{1,\sigma}(\mathfrak{u}_{r,s})=\mathbb{K}(1-\sigma), \quad
\textrm{for}\
 \sigma\not\in\{\omega_i, \, \omega_i'^{-1}\mid 1\le i\le n\}.$

\smallskip
$({\rm ii})$\;
$P_{\omega'_i,1}(\mathfrak{u}_{r,s})=\mathbb{K}(1-\omega'_i)+\mathbb{K}f_i;$
$\quad
P_{\omega_i^{-1},1}(\mathfrak{u}_{r,s})=\mathbb{K}(1-\omega_i^{-1})+\mathbb{K}e_i\omega_i^{-1};$

$\qquad P_{\sigma,1}(\mathfrak{u}_{r,s})=\mathbb{K}(1-\sigma), \quad
\textrm{for}\
 \sigma\not\in\{\omega_i^{-1}, \omega'_i\mid 1\le i\le n\}.$
\end{lemm}
\begin{proof} Similar to that of Lemma 4.3 in \cite{HW}.
\end{proof}

\noindent{\it 5.4.} Lemma 5.3 is crucial for determining the Hopf
isomorphisms of $\mathfrak{u}_{r,s}$.
\begin{theorem}\label{4.3}
Assume that $rs^{-1}$ and $r's'^{-1}$ are primitive $\ell$th roots
of unity with $\ell\ne 3, 4$, and $\zeta$ is a $2$nd root of unity.
Then $\varphi:\,\mathfrak{u}_{r,s}\cong\mathfrak{u}_{r',s'}$ as Hopf
algebras if and only if either
\begin{enumerate}
\item  $(r',s')=\zeta(r,s)$, $\varphi$ is a diagonal
isomorphism$:$ $\varphi(\omega_i)=\tilde{\omega}_i$,
$\varphi(\omega_i')=\tilde\omega_i'$, $\varphi(e_i)=a_i\tilde e_i$,
$\varphi(f_i)=\zeta^{\delta_{i,n}} a_i^{-1}\tilde f_i$; or

\item  $(r',s')=\zeta(s,r)$,
$\varphi(\omega_i)=\tilde\omega_i'^{-1}$,
$\varphi(\omega_i')=\tilde\omega_i^{-1}$, $\varphi(e_i)=a_i \tilde
f_i\tilde\omega_i'^{-1}$, $\varphi(f_i)=\zeta^{\delta_{i,n}}
a_i^{-1}\tilde\omega_i^{-1}\tilde e_i$, where $a_i\in\mathbb K^*$.
\end{enumerate}
\end{theorem}
\begin{proof}
Suppose $\varphi: \mathfrak{u}_{r,s}\longrightarrow
\mathfrak{u}_{r',s'}$ is a Hopf algebra isomorphism. Write
$\tilde{e}_i$, $\tilde{f}_i,\, \tilde\omega_i$, $\tilde\omega'_i$ to
distinguish the generators of $\mathfrak{u}_{r',s'}$. Because
$$
\Delta(\varphi(e_i))=(\varphi\otimes\varphi)
(\Delta(e_i))=\varphi(e_i)\otimes 1+\varphi(\omega_i)\otimes
\varphi(e_i),
$$
$\varphi(e_i) \in P_{1,\varphi(\omega_i)}(\mathfrak{u}_{r',s'})$. As
$\varphi$ is an isomorphism,
$\varphi(\omega_i)\in\mathbb{K}\tilde{G}$. Lemma 5.3 (i) implies
$\varphi(\omega_i)\in\{\,\tilde\omega_j, \tilde\omega_j'^{-1}\mid
1\le j\le n\,\}$. More precisely, we have either
$\varphi(\omega_i)=\tilde\omega_j$,
$\varphi(e_i)=a(1-\tilde\omega_j)+b\tilde e_j$, or
$\varphi(\omega_i)=\tilde\omega_j'^{-1}$,
$\varphi(e_j)=a(1-\tilde\omega_j'^{-1})+b\tilde
f_j\tilde\omega_j'^{-1}$ for some $j$ with $1\le j\le n$ and $a,
b\in\mathbb K$. But in both cases, we have $a=0$. The reasoning is
the following: Applying $\varphi$ to relation
$\omega_ie_i=r_is_i^{-1}e_i\omega_i$, we get
\begin{gather*}
\tilde\omega_j\bigl[\,a(1-\tilde\omega_j)+b\tilde
e_j\,\bigr]=r_is_i^{-1}\bigl[\,a(1-\tilde\omega_j)+b\tilde
e_j\,\bigr]\tilde\omega_j,
\quad\textit{or}\\
\tilde\omega_j'^{-1}\bigl[\,a(1-\tilde\omega_j'^{-1})+b\tilde
f_j\tilde\omega_j'^{-1}\,\bigr]=r_is_i^{-1}\bigl[\,a(1-\tilde\omega_j'^{-1})+b\tilde
f_j\tilde\omega_j'^{-1}\,\bigr]\tilde\omega_j'^{-1}.
\end{gather*}
Since $r_i\ne s_i$, $a=0$. Thus, $b\in\mathbb K^*$.

\smallskip
(I) We claim $\varphi(\omega_n)\in\{\tilde\omega_n,
\tilde\omega_n'^{-1}\}$. Indeed, let us consider the restriction
$\varphi|_\bullet$ of $\varphi$ to the Hopf subalgebra generated by
$e_i, f_i, \omega_i^{\pm1}, \omega_i'^{\pm1}$ for $i=n-1, n$, which
is isomorphic to $\mathfrak u_{r,s}(\mathfrak {so}_5)$. According to
the above discussion, we have $4$ possibilities to occur, namely,
case (i): For $i=n{-}1, n$, we have
$\varphi(\omega_i)=\tilde{\omega}_{j_i}$ and $\varphi(e_i)=a_i
\tilde{e}_{j_i}$ ($j_i\in\{n{-}1,n\}$, $a_i\in \mathbb{K}^*$); \
case (ii): For $i=n{-}1, n$, we have
$\varphi(\omega_i)=\tilde\omega_{j_i}'^{-1}$ and $\varphi(e_i)=a_i
\tilde{f}_{j_i}\tilde{\omega}_{j_i}'^{-1}$ ($j_i\in\{n{-}1, n\}$,
$a_i \in \mathbb{K}^*$); \ case (iii): $\varphi(\omega_{n-1})=\tilde
\omega_n$, but $\varphi(\omega_n)=\tilde\omega_{n-1}'^{-1}$; \ case
(iv): $\varphi(\omega_{n-1})=\tilde \omega_n'^{-1}$, but
$\varphi(\omega_n)=\tilde\omega_{n-1}$.

\smallskip
Case (i): In this case, we claim $\varphi(\omega_i)=\tilde\omega_i$
for $i=n{-}1, n$. Otherwise, if $\varphi(\omega_i)=\tilde\omega_j$
for $i\ne j\in\{n{-}1, n\}$, then applying $\varphi$ to relation
$\omega_ie_i=r_is_i^{-1}e_i\omega_i$ yields $\tilde{\omega}_j\tilde
e_j=r_is_i^{-1}\tilde e_j\tilde \omega_j$, it follows that
$r_j's_j'^{-1}=r_is_i^{-1}$ for $j\ne i\in\{n{-}1, n\}$. This means
$(rs^{-1})^3=1$, which contradicts the assumption. So we get
$\varphi(\omega_i)=\tilde{\omega}_i$ for $i=n{-}1, n$. Similarly, we
have $\varphi(\omega_i')=\tilde\omega_i'$ for $i=n{-}1, n$, as
$\varphi$ preserves $(B4)$. Thereby, we get $\varphi(e_i)=a_i\tilde
e_i$, $\varphi(f_i)=\zeta^{\delta_{i,n}} a_i^{-1}\tilde f_i$ for
$i=n{-}1, n$, where $a_i\in\mathbb K^*$.

Furthermore, applying $\varphi$ to relations $\omega_ie_j=\langle
\omega_j',\omega_i\rangle e_j\omega_i$ for $i, j\in\{n{-}1, n\}$, we
get $\langle \omega_j',\omega_i\rangle=\langle
\tilde\omega_j',\tilde\omega_i\rangle$ for $i, j\in\{n{-}1, n\}$,
i.e., $r'^2s'^{-2}=r^2s^{-2}$, $r's'^{-1}=rs^{-1}$, $r'^{-2}=r^{-2}$
and $s'^2=s^2$. Hence, $(r', s')=\zeta(r, s)$, $\varphi|_\bullet$ is
a diagonal isomorphism.

\smallskip
Case (ii): In this case, we claim
$\varphi(\omega_i)=\tilde\omega_i'^{-1}$ for $i=n{-}1, n$.
Otherwise, if $\varphi(\omega_i)=\tilde\omega_j'^{-1}$ for $i\ne
j\in\{n{-}1, n\}$, then $\varphi(e_i)=a_i \tilde
f_j\tilde\omega_j'^{-1}$ for $i\ne j\in\{n{-}1, n\}$. Using
$\varphi$ to relations $\omega_ie_i=\langle
\omega_i',\omega_i\rangle e_i\omega_i$ for $i\in\{n{-}1, n\}$, we
get $\tilde\omega_j'^{-1}\tilde f_j\tilde\omega_j'^{-1}=\langle
\omega_i',\omega_i\rangle \tilde
f_j\tilde\omega_j'^{-1}\tilde\omega_j'^{-1}$ so that $\langle
\tilde\omega_j',\tilde\omega_j\rangle^{-1}=\langle
\omega_i',\omega_i\rangle$ for $i\ne j\in\{n{-}1, n\}$. This means
$(rs^{-1})^3=1$. This contradicts the assumption. So we have
$\varphi(\omega_i)=\tilde\omega_i'^{-1}$ for $i=n{-}1, n$.
Similarly, $\varphi(\omega_i')=\tilde\omega_i^{-1}$ for $i=n{-}1,
n$, as $\varphi$ preserves $(B4)$. Hence,
$\varphi(\omega_i)=\tilde\omega_i'^{-1}$,
$\varphi(\omega_i')=\tilde\omega_i^{-1}$, $\varphi(e_i)=a_i \tilde
f_i\tilde\omega_i'^{-1}$ for $i=n{-}1, n$, where $a_i\in\mathbb
K^*$.

Using $\varphi$ to relations $\omega_ie_j=\langle
\omega_j',\omega_i\rangle e_j\omega_i$ for $i, j\in\{n{-}1, n\}$, we
get $\langle \tilde\omega_i',\tilde\omega_j\rangle^{-1}=\langle
\omega_j',\omega_i\rangle$ for $i, j\in\{n{-}1, n\}$, i.e.,
$r'^{-2}s'^2=r^2s^{-2}$, $(r's'^{-1})^{-1}=rs^{-1}$, $r'^2=s^2$ and
$s'^{-2}=r^{-2}$. Hence, we have $(r', s')=\zeta(s, r)$,
$\varphi(f_i)=\zeta^{\delta_{i,n}}a_i^{-1}\tilde\omega_i^{-1}\tilde
e_i$ for $i=n{-}1, n$. Besides, we can check that $\varphi$
preserves all of the $(r,s)$-Serre relations $(B5)$.

\smallskip
Case (iii): We claim that this case is impossible. First, we assume
that $\varphi(\omega_{n-1})=\tilde \omega_{n-1}$ and
$\varphi(\omega_n)=\tilde\omega_n'^{-1}$. Then we have
$\varphi(e_{n-1})=a_1\tilde e_{n-1}$, $\varphi(e_n)=a_2\tilde
f_n\tilde\omega_n'^{-1}$, where $a_1, a_2\in\mathbb K^*$. Using
$\varphi$ to relation $\omega_{n-1}e_n=s^2e_n\omega_{n-1}$, we get
$s^2=s'^{-2}$; to relation $\omega_ne_{n-1}=r^{-2}e_{n-1}\omega_n$
to get $r^{-2}=s'^2$. So we get $r^2=s^2$. This is a contradiction.

Next, we assume that $\varphi(\omega_{n-1})=\tilde \omega_n$, and
$\varphi(\omega_n)=\tilde\omega_{n-1}'^{-1}$. Then we have
$\varphi(e_{n-1})=a_1\tilde e_n$, $\varphi(e_n)=a_2\tilde
f_{n-1}\tilde\omega_{n-1}'^{-1}$, where $a_1, a_2\in\mathbb K^*$.
Using $\varphi$ to relation $\omega_{n-1}e_n=s^2e_n\omega_{n-1}$, we
 get $s^2=r'^2$; to relation
$\omega_ne_{n-1}=r^{-2}e_{n-1}\omega_n$ to get $r^{-2}=r'^{-2}$. So
we have $r^2=s^2$, which is contrary to the assumption.

\smallskip
Case (iv): Similarly to case (iii), it is also impossible.

\smallskip
(II) It suffices to consider the restriction $\varphi_i$ of
$\varphi$ to an $\mathfrak u_{r_i,s_i}(\mathfrak {sl}_3)$-copy in
$\mathfrak u_{r,s}(\mathfrak {so}_{2n+1})$ generated by $e_j, f_j,
\omega_j^{\pm 1}, \omega_j'^{\pm 1}$ for $j\in\{i, i{+}1\}$, where
$1\le i<n{-}1$. Note that $\varphi(\omega_j)\in\{\tilde\omega_i,
\tilde\omega_{i+1}, \tilde\omega_i'^{-1},
\tilde\omega_{i+1}'^{-1}\}$ for $j=i, i{+}1$.

Similarly to the proof of case (iii) or (iv) in (I), there don't
exist the possibilities: $(\varphi(\omega_i),
\varphi(\omega_{i+1}))=(\tilde \omega_{j_i},
\tilde\omega_{j_{i+1}}'^{-1})$, and $(\varphi(\omega_i),
\varphi(\omega_{i+1}))=(\tilde \omega_{j_i}'^{-1},
\tilde\omega_{j_{i+1}})$ for $j_i\ne j_{i+1}\in\{i, i{+}1\}$.
Thereby, we only need to consider the following possibilities: case
(1): $(\varphi(\omega_i), \varphi(\omega_{i+1}))=(\tilde \omega_i,
\tilde\omega_{i+1})$; case (2): $(\varphi(\omega_i),
\varphi(\omega_{i+1}))=(\tilde\omega_i'^{-1},
\tilde\omega_{i+1}'^{-1})$; case (3): $(\varphi(\omega_i),
\varphi(\omega_{i+1}))=(\tilde \omega_{i+1}, \tilde\omega_i)$; case
(4): $(\varphi(\omega_i),
\varphi(\omega_{i+1}))=(\tilde\omega_{i+1}'^{-1},
\tilde\omega_i'^{-1})$.

\smallskip
It is easy to treat the first two cases:

\smallskip
Case (1): As $(\varphi(\omega_i), \varphi(\omega_{i+1}))=(\tilde
\omega_i, \tilde\omega_{i+1})$, by the argument of the paragraph
above (I), we have $(\varphi(e_i), \varphi(e_{i+1}))=(a_i\tilde e_i,
a_{i+1}\tilde e_{i+1})$ and $(\varphi(f_i),
\varphi(f_{i+1}))=(a_i^{-1}\tilde f_i$, $a_{i+1}^{-1}\tilde
f_{i+1})$, where $a_i, a_{i+1}\in\mathbb K^*$. By a similar argument
of case (i) in (I), we have $(r'^2, s'^2)=(r^2, s^2)$.

\smallskip
Case (2): As $(\varphi(\omega_i),
\varphi(\omega_{i+1}))=(\tilde\omega_i'^{-1},
\tilde\omega_{i+1}'^{-1})$, by the argument of the paragraph above
(I), we have $(\varphi(e_i), \varphi(e_{i+1}))=(a_i\tilde
f_i\tilde\omega_i'^{-1}, a_{i+1}\tilde
f_{i+1}\tilde\omega_{i+1}'^{-1})$ and $(\varphi(f_i),
\varphi(f_{i+1}))$ $=(a_i^{-1}\tilde\omega_i^{-1}\tilde e_i,
a_{i+1}^{-1}\tilde\omega_{i+1}^{-1}\tilde e_{i+1})$, where $a_i,
a_{i+1}\in\mathbb K^*$. By a similar argument of case (ii) in (I),
we have $(r'^2, s'^2)=(s^2, r^2)$.

\smallskip
As for the latter two cases, we have

\smallskip
Case (3): As $(\varphi(\omega_i), \varphi(\omega_{i+1}))=(\tilde
\omega_{i+1}, \tilde\omega_i)$, using $\varphi$ to
$\omega_je_k=\langle \omega_k',\omega_j\rangle e_k\omega_j$ for $j,
k \in\{i, i{+}1\}$, we get $\langle \tilde\omega_j',
\tilde\omega_k\rangle=\langle \omega_k',\omega_j\rangle$. More
precisely, we have $r'^2s'^{-2}=r^2s^{-2}$, $r'^{-2}=s^2$,
$s'^2=r^{-2}$, i.e., $(r'^2, s'^2)=(s^{-2}, r^{-2})$. Hence,
$(\varphi(e_i), \varphi(e_{i+1}))=(a_{i+1}\tilde e_{i+1}, a_i\tilde
e_i)$ for $a_i, a_{i+1}\in\mathbb K^*$, $(\varphi(f_i),
\varphi(f_{i+1}))=(rs)^{-2}(a_{i+1}^{-1}\tilde f_{i+1},
a_i^{-1}\tilde f_i)$, since $\varphi$ preserves $(B4)$.

\smallskip
Case (4): As $(\varphi(\omega_i),
\varphi(\omega_{i+1}))=(\tilde\omega_{i+1}'^{-1},
\tilde\omega_i'^{-1})$, using $\varphi$ to $\omega_je_k=\langle
\omega_k',\omega_j\rangle e_k\omega_j$ for $j, k \in\{i, i{+}1\}$,
we get $\langle \tilde\omega_k', \tilde\omega_j\rangle^{-1}=\langle
\omega_k',\omega_j\rangle$. More precisely, we have
$r'^{-2}s'^2=r^2s^{-2}$, $r'^2=r^{-2}$, $s'^{-2}=s^2$, i.e., $(r'^2,
s'^2)=(r^{-2}, s^{-2})$. Hence, $(\varphi(e_i),
\varphi(e_{i+1}))=(a_{i+1}\tilde f_{i+1}\tilde\omega_{i+1}'^{-1},
a_i\tilde f_i\tilde\omega_i'^{-1})$ for $a_i, a_{i+1}\in\mathbb
K^*$, $(\varphi(f_i), \varphi(f_{i+1}))=(rs)^{-2}(a_{i+1}^{-1}\tilde
\omega_{i+1}^{-1}\tilde e_{i+1}$, $a_i^{-1}\tilde\omega_i^{-1}\tilde
e_i)$, since $\varphi$ preserves $(B4)$. Moreover, we can verify
that $\varphi$ preserves all of the $(r,s)$-Serre relations $(B5)$.

\smallskip
In view of the above discussions in cases (1)---(4), we derive the
complete description below on all the Hopf algebra isomorphisms of
$\mathfrak u_{r,s}(\mathfrak {sl}_n)$ (some of cases were missed in
\cite{BW3}).

\begin{theorem}
Assume that $rs^{-1}$ and $r's'^{-1}$ are primitive $\ell$th roots
of unity. Then
$\varphi:\,\mathfrak{u}_{r,s}(\mathfrak{sl}_n)\cong\mathfrak{u}_{r',s'}(\mathfrak{sl}_n)$
as Hopf algebras if and only if either
\begin{enumerate}
\item
 $(r',s')=(r,s)$, $\varphi$ is a diagonal isomorphism$:$
$\varphi(\omega_i)=\tilde{\omega}_i$,
$\varphi(\omega_i')=\tilde\omega_i'$, $\varphi(e_i)=a_i\tilde e_i$,
$\varphi(f_i)=a_i^{-1}\tilde f_i;$ or
\item
$(r',s')=(s,r)$, $\varphi(\omega_i)=\tilde\omega_i'^{-1}$,
$\varphi(\omega_i')=\tilde\omega_i^{-1}$, $\varphi(e_i)=a_i \tilde
f_i\tilde\omega_i'^{-1}$, $\varphi(f_i)=
a_i^{-1}\tilde\omega_i^{-1}\tilde e_i;$ or
\item
 $(r', s')=(s^{-1}, r^{-1})$,
$\varphi(\omega_i)=\tilde\omega_{n-i}$,
$\varphi(\omega_i')=\tilde\omega_{n-i}'$,
$\varphi(e_i)=a_{n-i}\tilde e_{n-i}$,
$\varphi(f_i)=(rs)^{-1}a_{n-i}^{-1}\tilde f_{n-i};$ or
\item
 $(r', s')=(r^{-1}, s^{-1})$,
$\varphi(\omega_i)=\tilde\omega_{n-i}'^{-1}$,
$\varphi(\omega_i')=\tilde\omega_{n-i}^{-1}$,
$\varphi(e_i)=a_{n-i}\tilde f_{n-i}\tilde\omega_{n-i}'^{-1}$,
$\varphi(f_i)=(rs)^{-1}a_{n-i}^{-1}\tilde\omega_{n-i}^{-1}\tilde
e_{n-i}$, $(a_i\in\mathbb K^*)$.
\end{enumerate}
\end{theorem}

Now continue to the proof of Theorem 5.4. According to the
discussions in (I), there are $2$ possibilities: (i)
$\varphi(\omega_{n-1})=\tilde \omega_{n-1}$; (ii)
$\varphi(\omega_{n-1})=\tilde \omega_{n-1}'^{-1}$. Combining with
the discussions in (II), we obtain two kinds of the Hopf algebra
isomorphisms of $\mathfrak u_{r,s}(\mathfrak {so}_{2n+1})$ as
required, which are only compatible with the cases (1) and (2) in
(II), respectively.
\end{proof}

\section{$\mathfrak{u}_{r,s}(\mathfrak{so}_{2n+1})$ is a Drinfel'd double}

\noindent{\it 6.1.} Let $\theta$ be a primitive $\ell$th root of
unity in $\mathbb{K}$, and write $r=\theta^y$, $s=\theta^z$.
\begin{lemm}\label{5.1}
Assume that $(2^{n-1}(y^n+(-1)^nz^n),\ell)=1$ and $rs^{-1}$ is a
primitive $\ell$th root of unity with $\ell\ne 3$. Then
$(\mathfrak{b}')^{coop}\cong \mathfrak{b}^*$ as Hopf algebras.
\end{lemm}
\begin{proof}
Define $\gamma_j\in \mathfrak{b}^*$ such that $\gamma_j$'s are
algebra homomorphisms with
$$\gamma_{j}(e_i)=0, \ \forall\; i, \textrm{ and } \ \gamma_{j}(\omega_i)=\langle \omega_j',\omega_i\rangle.$$
So they are group-like elements in $\mathfrak{b}^*$. Define
$\eta_j=\sum\limits_{g \in G(\mathfrak b)}(e_jg)^*$, where
$G(\mathfrak{b})$ is  the group generated by $\omega_i\ (1\leq i\leq
n)$ and the asterisk denotes the dual basis element relative to the
PBW-basis of $\mathfrak{b}$. The isomorphism $\phi:
(\mathfrak{b}')^{coop}\rightarrow \mathfrak{b}^*$ is defined by
$$
\phi(\omega'_j)=\gamma_j, \qquad \phi(f_j)=\eta_j.
$$
First, we will check that $\phi$ is a Hopf algebra homomorphism and
then we will show that it is a bijection.

Clearly, the $\gamma_j$'s are invertible elements in
$\mathfrak{b}^*$ that commute with one another and
$\gamma_j^\ell=1$. Note that $\eta_j^\ell=0$, as it is $0$ on any
basis element of $\mathfrak{b}$. We calculate
$\gamma_j\eta_i\gamma_j^{-1}$: It is nonzero only on basis elements
of the form $e_i\omega_1^{k_1}\cdots\omega_n^{k_n}$, and on such an
element it takes the value
\begin{gather*}
\begin{split}
& (\gamma_j\otimes\eta_i\otimes \gamma_j^{-1})((e_i\otimes 1\otimes
1+ \omega_i\otimes e_i\otimes 1+\omega_i\otimes \omega_i\otimes
e_i)(\omega_1^{k_1}\cdots\omega_n^{k_n})^{\otimes 3}) \\&=
\gamma_j(\omega_i\omega_1^{k_1}\cdots\omega_n^{k_n})\eta_i(e_i\omega_1^{k_1}\cdots\omega_n^{k_n})
\gamma_j^{-1}(\omega_1^{k_1}\cdots\omega_n^{k_n})\\
&=\gamma_{j}(\omega_i).
\end{split}
\end{gather*}
So, we have
$\gamma_j\eta_i\gamma_j^{-1}=\langle
\omega_j',\omega_i\rangle\,\eta_i$,
which corresponds to relation $(B3)$ for $\mathfrak{b}'$.

Next, with the same argument in the proof of Lemma 5.1 in \cite{HW}
(see, p. 263--264), we easily check relations in $(B5)$:
\begin{gather*}
r^2s^2\eta_i^2\eta_{i+1}-(r^2{+}s^2)\eta_i\eta_{i+1}\eta_i+\eta_{i+1}\eta_i^2=0.\\
\eta_i\eta_{i+1}^2-(r^{-2}{+}s^{-2})\eta_{i+1}\eta_{i}\eta_{i+1}+(rs)^{-2}\eta_{i+1}^2\eta_i=0,\qquad i<n{-}1.\\
 \eta_{n-1}\eta_n^{3}-(r_n^{-2}{+}r_n^{-1}s_n^{-1}
{+} s_n^{-2})\,\eta_n\eta_{n-1}\eta_n^{2} \\  +\,
(r_ns_n)^{-1}(r_n^{-2} {+} r_n^{-1}s_n^{-1}{+}
s_n^{-2})\,\eta_{n}^2\eta_{n-1}\eta_n -
(r_ns_n)^{-3}\,\eta_{n}^3\eta_{n-1} = 0.
\end{gather*}
Hence, $\phi $ is an algebra homomorphism.

We have already seen that $\gamma_i$ is a group-like element in
$\mathfrak{b}^*$. Now using the same argument in the proof of Lemma
5.1 in \cite{HW} (see, p. 264--265), we have 
$$\Delta(\eta_i)=\eta_i\otimes 1+\gamma_i\otimes \eta_i,$$
which means that $\phi$ is a Hopf algebra homomorphism.

Finally, we prove $\phi$ is bijective. As
$\dim\,\mathfrak{b}^*=\dim\,\mathfrak{b}'^{coop}$, it suffices to
show $\phi$ is injective. By [{\bf M}], we need only to show
$\phi|_{(\mathfrak{b}')_1^{coop}}$ is injective. Lemma \ref{4.2}
yields $(\mathfrak{b}')_1^{coop}=\mathbb{K}G(\mathfrak{b}')
+\sum_{i=1}^{n}\mathbb{K}f_iG(\mathfrak{b}')$, where
$G(\mathfrak{b}')$ is the group generated by $\omega'_i\ (1\leq
i\leq n)$. First, we claim that
$$
\textrm{span}_\mathbb{K}\Bigl\{\left.
\gamma_1^{k_1}\cdots\gamma_n^{k_n}\,\right|\, 0\le k_i<\ell\,\Bigr\}
=\textrm{span}_\mathbb{K}\Bigl\{\left.(\omega_1^{k_1}\cdots\omega_n^{k_n})^*\,\right|\,
0\le k_i<\ell\,\Bigr\}. \quad \leqno (6.1)
$$
This is equivalent to saying that the
$\gamma_1^{k_1}\cdots\gamma_n^{k_n}$'s span the space of characters
over $\mathbb{K}$ of the finite group $\mathbb{Z}/\ell
\mathbb{Z}\times \cdots\times \mathbb{Z}/\ell \mathbb{Z}$ generated
by $\omega_1,\cdots,\omega_n$. Since we assumed that $\mathbb{K}$
contains a primitive $\ell$th root of unity, the irreducible
characters of this group are the functions $\chi_{i_1,\cdots,i_n}$
given by
$$\chi_{i_1,\cdots,i_n}(\omega_1^{k_1}\cdots\omega_n^{k_n})=\theta
^{i_1k_1+\cdots+i_nk_n},$$ where $\theta$ is a primitive $\ell$th
root of unity in $\mathbb{K}$. Note that
\begin{gather*}
\begin{split}
\gamma_1&=\chi_{2(y-z),-2y,0,\cdots,0},\\
\gamma_2&=\chi_{2z,2(y-z),-2y,0,\cdots,0},\\
\vdots\\
\gamma_{n-1}&=\chi_{0,\cdots,0,2z,2(y-z),-2y},\\
\gamma_{n}&=\chi_{0,\cdots,0,2z,y-z}.
\end{split}
\end{gather*}
We must show that, given $i_1,\cdots,i_n$, there are $k_1, \cdots,
k_n$ such that
$$
\chi_{i_1,\cdots,i_n}=\gamma_1^{k_1}\cdots\gamma_n^{k_n},
$$
which is equivalent to the existence of a solution to the matrix
equation
$$AK=I$$
in $\mathbb{Z}/\ell \mathbb{Z}$ (as these are powers of
$\theta$), where
\[A= \begin{pmatrix} 2(y-z) & 2z &0 & 0& \cdots &0& 0\\
-2y & 2(y-z)& 2z &0&  \cdots&0&0 \\
0&-2y & 2(y-z)& 2z & \cdots&0&0\\
\vdots&\vdots&\vdots&\ddots&\ddots&\vdots&\vdots\\
0&0&0&0&\cdots&2z& 0\\
0&0&0&0&\cdots&\quad 2(y-z)& 2z\\
0&0&0&0&\cdots&-2y& y-z
\end{pmatrix},\]
$K$ is the transpose of $(k_1,\cdots,k_n)$ and $I$ is the transpose
of $(i_1,\cdots,i_n)$. The determinant of  the coefficient matrix
$A$
 is $2^{n-1}(y^n+(-1)^nz^n)$, which is invertible in $\mathbb{Z}/\ell
 \mathbb{Z}$ by the hypothesis in the Lemma. So, (6.1) holds.
 In particular, this implies that the matrix
 $$
 \big((\gamma_1^{k_1}\cdots\gamma_n^{k_n})(\omega_1^{j_1}\cdots\omega_n^{j_n})\big)_{\bar
 k\times \bar j} \leqno (6.2)
 $$
 is invertible, and that $\phi$ is
 bijection on group-like elements.

 Next, we will show for each $i\ (1\leq i\leq n)$ that the following
 matrix is invertible:
$$
\big(
 (\eta_i\gamma_1^{k_1}\cdots\gamma_n^{k_n})(e_i\omega_1^{j_1}\cdots\omega_n^{j_n})\big)_{\bar
 k\times \bar j}\,. \leqno (6.3)
 $$
This will complete the proof that
 $\phi$ is injective on $(\mathfrak{b}')_1^{coop}$, as desired. We
 will show that the matrix is block upper-triangular. Each matrix
 entry is
\begin{gather*}
\begin{split}
&(\eta_i\otimes\gamma_1^{k_1}\cdots\gamma_n^{k_n})(\Delta(e_i)
\Delta(\omega_1^{j_1}\cdots\omega_n^{j_n}))\\
&\quad=(\eta_i\otimes\gamma_1^{k_1}\cdots\gamma_n^{k_n})(e_i\omega_1^{j_1}\cdots\omega_n^{j_n}
\otimes\omega_1^{j_1}\cdots\omega_n^{j_n}).
\end{split}
\end{gather*}
Thus, (6.3) is precisely the invertible matrix (6.2).
\end{proof}

\noindent{\it 6.2.} With the same argument of Theorem 4.8 in
\cite{BW3}, we have
\begin{theorem}\label{5.2}
Assume $(2^{n-1}(y^n+(-1)^nz^n),\ell)=1$, then $D(\mathfrak{b})\cong
\mathfrak{u}_{r,s}(\mathfrak{so}_{2n+1})$ as Hopf algebras.
\hfill\qed
\end{theorem}

\section{Integrals}

\noindent{\it 7.1.} {\it Recall some definitions}. Let $H$ be a
finite-dimensional Hopf algebra. An element $y\in H$ is a
\textit{left} (resp., \textit{right}) \textit{integral} if
$ay=\varepsilon(a)y$ (resp., $ya=\varepsilon(a)y$) for all $a \in
H$. The left (resp., right) integrals form a one-dimensional ideal
$\int^l_H$ (resp., $\int^r_H$) of $H$, and $S_H(\int^r_H)=\int^l_H$
under the antipode $S_H$ of $H$.

\smallskip
When $y\neq 0$ is a left integral of $H$, there exists a unique
group-like element $\gamma$ in the dual algebra $H^*$ (the so-called
\textit{distinguished group-like element} of $H^*$) such that
$ya=\gamma(a)y$.  If we had begun instead with a right integral
$y'\in H$, then we would have $ay'=\gamma^{-1}(a)y'$. This follows
from the fact that group-like elements are invertible, and the
following calculation, which can be found in [\textbf{M}, p.22]:
When $y' \in \int^r_H$, $S_H(y')$ is a nonzero multiple of $y$, so
that $S_H(y')S_H(a)=\gamma(S_H(a))S_H(y')$ for all $a \in H$.
Applying $S_H^{-1}$, we find $ay'=\gamma(S_H(a))y'$. As $\gamma$ is
group-like, we have
$\gamma(S_H(a))=S_{H^*}(\gamma)(a)=\gamma^{-1}(a).$

\smallskip
Now if $\lambda \neq 0$ is right integral of $H^*$, then there
exists a unique group-like element $g$ of $H$ (the
\textit{distinguished group-like element} of $H$) such that
$\xi\lambda=\xi(g)\lambda$ for all $\xi \in H^*$. The algebra $H$ is
\textit{unimodular} (i.e., $\int^l_H=\int^r_H$) if and only if
$\gamma=\varepsilon$; and the dual algebra $H^*$ is unimodular if
and only if $g=1$.

\smallskip
The left and right $H^*$-module actions on $H$ are given by
$$
\xi\rightharpoonup a=\sum a_{(1)}\xi(a_{(2)}),\qquad
a\leftharpoonup\xi=\sum \xi(a_{(1)})a_{(2)}, \leqno(7.1)
$$
for all
$\xi \in H^*$ and $a\in H$. In particular,
$\varepsilon\rightharpoonup a=a=a\leftharpoonup\varepsilon$ for all
$a \in H$. 

\smallskip
\noindent{\it 7.2.} The first aim of this section is to construct
the left and right integrals in the Borel subalgebra $\mathfrak{b}$
of $\mathfrak{u}_{r,s}(\mathfrak{so}_{2n+1})$. To this end, we
introduce the elements:
$$
t=\prod\limits_{i=1}^{n}(1+\omega_i+\cdots+\omega_i^{\ell-1}),
\quad  x=\prod\limits_{i=1}^{n}\mathcal{E}_{\widehat{i,i+1}},
$$
where
$\mathcal{E}_{\widehat{i,i+1}}=\mathcal{E}_{i,i}^{\ell-1}\mathcal{E}_{i,i+1}^{\ell-1}
\cdots\mathcal{E}_{i,n}^{\ell-1} \mathcal{E}_{i,n'}^{\ell-1}\cdots
\mathcal{E}_{i,i+1'}^{\ell-1}$ ($1\le i<n$),
$\mathcal{E}_{\widehat{n,n+1}}=\mathcal{E}_{n,n}^{\ell-1}$.

\begin{theorem}\label{6.1}
The element $y=tx$ is a left integral in $\mathfrak{b}$.
\end{theorem}
\begin{proof}
We need to show that $by=\varepsilon(b)y$ for all $b \in
\mathfrak{b}$. It suffices to show this for the generators
$\omega_k$ and $e_k$, as the counit $\varepsilon$ is an algebra
homomorphism.

Observe that $\omega_kt=t=\vn(\om_k)t$ for all $k=1,\cdots,n$, as
the $\om_i$'s commute and
$\om_k(1+\om_k+\cdots+\om_k^{\ell-1})=1+\om_k+\cdots+\om_k^{\ell-1}$.
From that, relation $\om_k y=\vn(\om_k)y$ is clear,  for all $k$.

Next we compute $e_ky$. By a direct calculation, we get
\begin{gather*}
\begin{split}
e_kt&=\prod\limits_{i=1}^{n-1}\bigl(1+r^{-2(\ep_i,\a_k)}s^{-2(\ep_{i+1},\a_k)}\om_i+
\cdots+r^{-2(\ell-1)(\ep_i,\a_k)}s^{-2(\ell-1)(\ep_{i+1},\a_k)}\om_i^{\ell-1}\bigr)\cdot\\
& \quad\cdot\Bigl[\de_{kn}(1+r^{-1}s\om_n+\cdots+r^{-(\ell-1)}s^{\ell-1}\om_n^{\ell-1})\\
&\quad
+(1-\de_{kn})(1+r^{-2(\ep_n,\a_k)}\om_n+\cdots
+r^{-2(\ell-1)(\ep_n,\a_k)}\om_n^{\ell-1})\Bigr]e_k.
\end{split}
\end{gather*}
So it suffices to show that $e_kx=0=\vn(e_k)x$.

Now
\begin{gather*}
x=\mathcal{E}_{\widehat{{1,2}}}\mathcal{E}_{\widehat{{2,3}}}\cdots
\mathcal{E}_{\widehat{{k-1,k}}}(\mathcal{E}_{k,k}^{\ell-1}\cdots
\mathcal{E}_{k,k+1'}^{\ell-1})\cdots
\mathcal{E}_{\widehat{{n,n+1}}}, \\
e_1x=0.
\end{gather*}
We want to show $e_kx=0$, that is, $e_k\ (2\leq k \leq n)$ can be
moved across the terms $\mathcal{E}_{\widehat{{i,i+1}}}\ (1\leq i<
k)$ next to $ \mathcal{E}_{k,k}^{\ell-1}$. By Lemma 3.1 (1), $e_k$
commutes with
$\mathcal{E}_{i,i}^{\ell-1}\cdots\mathcal{E}_{i,k-2}^{\ell-1}$ when
$i\le k{-}2\le n{-}2$. By Lemmas 3.9 (2) \& 3.3 (5), we obtain
$$
e_k\mathcal{E}_{i,k-1}^{\ell-1}=r^2\mathcal{E}_{i,k-1}^{\ell-1}e_k+
s^2\mathcal{E}_{i,k-1}^{\ell-2}\mathcal{E}_{i,k}, \qquad(k\le n).
$$
The second term reaches $\mathcal{E}_{i,k}^{\ell-1}$ and we get
$\mathcal{E}_{i,k}^{\ell}=0$. Hence, it suffices to treat the first
term $r^2\mathcal{E}_{i,k-1}^{\ell-1}e_k$ $(k\le n)$.

\smallskip
(i) \ For $k<n$:  by Lemmas 3.2 (3), 3.3 (1) \& (6), $e_k$
quasi-commutes with
$\mathcal{E}_{i,k}^{\ell-1}\mathcal{E}_{i,k+1}^{\ell-1}\cdots\mathcal{E}_{i,k+2'}^{\ell-1}$
(up to a factor $s^2$). By Lemmas 3.10 (3)  \& 3.6 (5), we obtain
$$
e_k\mathcal{E}_{i,k+1'}^{\ell-1}=s^{-2}\mathcal{E}_{i,k+1'}^{\ell-1}e_k+
r^{-2}\mathcal{E}_{i,k+1'}^{\ell-2}\mathcal{E}_{i,k'}, \qquad(k\le
n{-}1).
$$
The second term meets $\mathcal{E}_{i,k'}^{\ell-1}$ and yields
$\mathcal{E}_{i,k'}^{\ell}=0$. By Lemma 3.2 (5), we have
$e_k\mathcal{E}_{i,k'}^{\ell-1}=r^{-2}\mathcal{E}_{i,k'}^{\ell-1}e_k$
$(i<k<n)$. This implies that $e_k$ quasi-commutes with $\mathcal
E_{\widehat{i,i+1}}$ for all $i<k$. Since
$e_k\mathcal{E}_{\widehat{k,k+1}}=0$, $e_kx=0$ for $k<n$.

\smallskip
(ii) \ For $k=n$: by Lemmas 3.9 (4) \& 3.5 (3), we have $e_n\mathcal
E_{i,n}^{\ell-1}=(rs)\mathcal E_{i,n}^{\ell-1}e_n+s^2\mathcal
E_{i,n}^{\ell-2}\mathcal E_{i,n'}$. The second term meets $\mathcal
E_{i,n'}^{\ell-1}$ and yields $\mathcal E_{i,n'}^\ell=0$, so it
remains to observe that $e_n$ quasi-commutes with $\mathcal
E_{i,n'}^{\ell-1}\mathcal E_{i,n-1'}^{\ell-1}\cdots\mathcal
E_{i,i+1'}^{\ell-1}$ (by Lemmas 3.2 (4) \& 3.4 (4)), that is, $e_n$
quasi-commutes with those $\mathcal E_{\widehat{i,i+1}}$ for $i<n$
and finally meets $\mathcal E_{\widehat{n,n+1}}$ and yields $0$.

\medskip
This completes the proof.
\end{proof}

\begin{theorem}\label{6.2}
The element $y'=xt$ is a right integral in $\mathfrak{b}$.
\end{theorem}
\begin{proof}
Arguing as in Theorem 7.1, we see that $t\omega_j=t=\vn(\om_j)t$,
and hence that
 $y'\om_j=\vn(\om_j)y'$ for all
$j$. As
\begin{gather*}
\begin{split}
te_j&=e_j\prod\limits_{i=1}^{n-1}(1+r^{2(\ep_i,\a_j)}s^{2(\ep_{i+1},\a_j)}\om_i+
\cdots+r^{2(\ell-1)(\ep_i,\a_j)}s^{2(\ell-1)(\ep_{i+1},\a_j)}\om_i^{\ell-1})\cdot\\
& \quad
\cdot\Bigl[\de_{jn}e_j(1+rs^{-1}\om_n+\cdots+r^{\ell-1}s^{-(\ell-1)}\om_n^{\ell-1})\\
&\quad +(1-\de_{jn})\,e_j(1+r^{2(\ep_n,\a_j)}\om_n+\cdots
+r^{2(\ell-1)(\ep_n,\a_j)}\om_n^{\ell-1})\Bigr],
\end{split}
\end{gather*}
so it suffices to show that $xe_j=0$.

Now
\begin{gather*}
x=\mathcal{E}_{\widehat{{1,2}}}\cdots
\mathcal{E}_{\widehat{{j,j+1}}}\Bigl(\mathcal{E}_{j+1,j+1}^{\ell-1}\cdots\mathcal
E_{j+1,k'}^{\ell-1}\cdots\mathcal E_{j+1,j+2'}^{\ell-1}\Bigr)
\cdots\mathcal{E}_{\widehat{n,n+1}}, \\
xe_n=0.
\end{gather*}
We need to prove $xe_j=0$ for $1\le j<n$. By Lemmas 3.9 (3), (4),
(6) \& 3.4 (3), it is easy to see that
$\mathcal{E}_{\widehat{{n-1,n}}}\mathcal{E}_{\widehat{{n,n+1}}}e_{n-1}=0$.
It suffices to show
$\mathcal{E}_{\widehat{{j,j+1}}}\mathcal{E}_{\widehat{{j+1,j+2}}}e_j=0$
since $e_j$ commutes with $\mathcal{E}_{\widehat{j+2,j+3}}
\cdots\mathcal E_{\widehat{n,n+1}}$ (by Lemma 3.1 (1)) for $j< n-1$.

\medskip
(I) \ Claim: $\mathcal{E}_{\widehat{{j,j+1}}}e_j=0$ for $j<n{-}1$.

Write
\begin{equation*}
\begin{split}
\Pi_j^{(k)}:&=\mathcal E_{j,j}^{\ell-1}\mathcal
E_{j,j+1}^{\ell-1}\cdots\mathcal E_{j,k}^{\ell-1},\quad (j\le k\le n),\\
\Pi_j^{(k')}:&=\Pi_j^{(n)}\,\mathcal E_{j,n'}^{\ell-1}\cdots\mathcal
E_{j,k'}^{\ell-1},\quad (n'\le k'\le j{+}1').
\end{split}
\end{equation*}
Thus, $\Pi_j^{(j{+}1')}=\mathcal E_{\widehat{j,j+1}}$, and
$\Pi_j^{(j{+}2')}\,e_j=0$, as $e_j\mathcal E_{j,k}=s^2\mathcal
E_{j,k}e_j$ for $j<k\le n$, $e_j\mathcal E_{j,k'}=s^2\mathcal
E_{j,k'}e_j$ for $j{+}1<k\le n$, and $e_j^\ell=0$.

\smallskip
For $j<n{-}1$, by Lemma 3.10 (4) \& (5) and $\mathcal
E_{j,j+2'}^{\ell}=\cdots=\mathcal E_{j,n'}^{\ell}=\mathcal
E_{j,n}^\ell=0$, we get
\begin{equation*}
\begin{split}
\Pi_j^{(j{+}1')}e_j&=\Pi_j^{(j{+}2')}\Bigl(\mathcal{E}_{j,j+1'}^{\ell-1}e_j\Bigr)\hskip3.5cm\,\;
\text{(by Lemma 3.10 (4))}\\
&=*_1\Pi_j^{(j{+}2')}\mathcal E_{j+1,j+1}^{\ell-1}\cdots\mathcal
E_{j+1,n-2}^{\ell-1}\cdots\,
\mathcal{E}_{j,j+1}\mathcal{E}_{j,j+2'}\mathcal E_{j,j+1'}^{\ell-2}\\
&=*_2\Pi_j^{(j{+}3')}\,\Bigl(\mathcal E_{j,j+2'}^{\ell-1}\mathcal
E_{j,j+1}\Bigr)\mathcal{E}_{j,j+2'}\mathcal{E}_{j,j+1'}^{\ell-2}\qquad
\text{(by Lemma 3.10 (4))}\\
&=*_2\Pi_j^{(j{+}3')}\,\Bigl(\mathcal E_{j,j+1}\mathcal
E_{j,j+2'}^{\ell-1}-(rs)^{-2}\mathcal E_{j,j+2}\mathcal
E_{j,j+3'}\mathcal
E_{j,j+2'}^{\ell-2}\Bigr)\mathcal{E}_{j,j+2'}\mathcal{E}_{j,j+1'}^{\ell-2}\\
&=*_3\Pi_j^{(j{+}4')}\Bigl(\mathcal E_{j,j+3'}^{\ell-1}\mathcal
E_{j,j+2}\Bigr)\mathcal{E}_{j,j+3'}\mathcal
E_{j,j+2'}^{\ell-1}\mathcal{E}_{j,j+1'}^{\ell-2}\\
&=\cdots \qquad\,\hskip5cm\text{(by Lemma 3.10 (4))}\\
&=*\,\Pi_j^{(n)}\Bigl(\mathcal E_{j,n'}^{\ell-1}\mathcal
E_{j,n-1}\Bigr)\mathcal E_{j,n'}\mathcal
E_{j,n-1'}^{\ell-1}\cdots\mathcal E_{j,j+2'}^{\ell-1}\mathcal
E_{j,j+1'}^{\ell-2}\ \text{(by Lemma 3.10 (5))}\\
&=*'\Pi_j^{(n)}\Bigl(\mathcal E_{j,n}^2\mathcal
E_{j,n'}^{\ell-2}\Bigr)\mathcal E_{j,n'}\mathcal
E_{j,n-1'}^{\ell-1}\cdots\mathcal E_{j,j+2'}^{\ell-1}\mathcal
E_{j,j+1'}^{\ell-2}\\
&=0, \qquad\text{(since $\Pi_j^{(n)}\mathcal E_{j,n}=0$)}.
\end{split}
\end{equation*}

(II) \ Write $\Delta_j^{(k)}:=\Pi_j^{(j{+}1')}\Pi_{j+1}^{(k)}$ for
$j{+}1\le k\le n$. Then
$\Delta_j^{(n)}=\Pi_j^{(j{+}1')}\Pi_{j+1}^{(n)}$. We make a
convention: $\Delta_j^{(j)}:=\Pi_j^{(j{+}1')}$.

\smallskip
(i) \ Assume $\Delta_j^{(k)}e_j=0$ for $j\le k<n$, we want to prove
$\Delta_j^{(k+1)}e_j=0$.

\smallskip
(1) \ For $k+1<n$: indeed, by Lemmas 3.9 (1) \& 3.3 (4),  we have
\begin{equation*}
\begin{split}
\Delta_j^{(k+1)}e_j&=\Delta_j^{(k)}\Bigl(\mathcal
E_{j+1,k+1}^{\ell-1}e_j\Bigr)\\
&=\Delta_j^{(k)}\Bigl(r^2e_j\mathcal E_{j+1,k+1}+s^2\mathcal
E_{j,k+1}\Bigr)\mathcal
E_{j+1,k+1}^{\ell-2}\quad\text{(by Lemma 3.3 (3))}\\
&=s^2\Pi_j^{(j{+}1')}\mathcal E_{j,k+1}\Pi_{j+1}^{(k)}\mathcal
E_{j+1,k+1}^{\ell-2}.
\end{split}
\end{equation*}
Lemmas 3.7 (8), 3.5 (3), 3.3 (9) show that $\mathcal E_{j,k+1}$
quasi-commutes with any element of $\{\,\mathcal E_{j,l'}\mid j+1\le
l\le n, \ l\ne k+2\,\}$. Hence, by a similar argument in (I) (using
Lemma 3.10 (4) \& (5)), we can get
\begin{equation*}
\begin{split}
\Pi_j^{(j{+}1')}\mathcal
E_{j,k+1}&=*_1\,\Pi_j^{(k{+}3')}\Bigl(\mathcal
E_{j,k+2'}^{\ell-1}\mathcal E_{j,k+1}\Bigr)\mathcal
E_{j,k+1'}^{\ell-1}\cdots\mathcal E_{j,j+1'}^{\ell-1}\\
&=*_2\,\Pi_j^{(k{+}3')}\Bigl(\mathcal E_{j,k+1}\mathcal
E_{j,k+2'}^{\ell-1}{-}(rs)^{-2}\mathcal E_{j,k+2}\mathcal
E_{j,k+3'}\mathcal E_{j,k+2'}^{\ell-2}\Bigr)\cdots\mathcal E_{j,j+1'}^{\ell-1}\\
&=0+*_3\,\Pi_j^{(k{+}4')}\Bigl(\mathcal E_{j,k+3'}^{\ell-1}\mathcal
E_{j,k+2}\Bigr)\mathcal E_{j,k+3'}\mathcal
E_{j,k+2'}^{\ell-2}\mathcal
E_{j,k+1'}^{\ell-1}\cdots\mathcal E_{j,j+1'}^{\ell-1}\\
&=0,
\end{split}
\end{equation*}
here in the 1st summand $0$ we used the fact that $\mathcal
E_{j,k+1}$ quasi-commutes with $\mathcal
E_{j,n'}^{\ell-1}\cdots\mathcal E_{j,k+3'}^{\ell-1}$, as well as
$\mathcal E_{j,l}$ (by Lemma 3.3 (5)) for any $k{+}2\le l\le n$, and
$\mathcal E_{j,k+1}^\ell=0$; while the 2nd summand $0$ used
repeatedly the proof-procedure of (I).

\smallskip
(2) \ For $k+1=n$: at first we note that
$$
\Delta_j^{(n-1)}\mathcal E_{j,n}=0, \qquad \Delta_j^{(n-1)}\mathcal
E_{j+1,n-1}=0,
$$
since Lemmas 3.3 (3) \& 3.5 (3) say: $\mathcal E_{j,n}$ commutes
with $\mathcal E_{j+1,l}$ for $j{+}1\le l\le n{-}1$, and
quasi-commutes with $\mathcal E_{j,l'}$ for $j{+}1\le l\le n$. Thus,
by Lemma 3.4 (1), we have
$$
\mathcal E_{j+1,n}\mathcal E_{j,n}\equiv rs^{-1}\mathcal
E_{j,n'}\mathcal E_{j+1,n-1}\quad \mod \,(\mathcal E_{j,n}\mathcal
E_{j+1,n}+\mathcal E_{j+1,n-1}\mathcal E_{j,n'}).\leqno(7.2)
$$

Furthermore, we note that
$$
\Delta_j^{(k-2)}\mathcal E_{j,k'}=0, \qquad \Delta_j^{(k-2)}\mathcal
E_{j{+}1,k-2}=0, \qquad (j{+}2< k\le n).\leqno(7.3)
$$
Hence, by Lemma 3.4 (2), for $j{+}2<k\le n$, we get
$$
\mathcal E_{j+1,k-1}\mathcal E_{j,k'}\equiv (rs)^2\mathcal
E_{j,k-1'}\mathcal E_{j+1,k-2}\ \mod\,(\mathcal E_{j,k'}\mathcal
E_{j+1,k-1}{-}\mathcal E_{j+1,k-2}\mathcal E_{j,k-1'}).\leqno(7.4)
$$

Using (7.2), (7.3) \& (7.4), by Lemmas 3.9 (5), 3.3 (5) \& 3.10 (6),
we get
\begin{equation*}
\begin{split}
\Delta_j^{(n)}e_j&=\Delta_j^{(n-1)}\left(\mathcal
E_{j+1,n}^{\ell-1}e_j\right) =\Delta_j^{(n-1)}\Bigl[r^2e_j\mathcal
E_{j+1,n}^{\ell-1}+(rs)^{-1}\mathcal E_{j+1,n}^{\ell-2}\mathcal
E_{j,n}\\
&\quad+(rs^3)^{-1}\mathcal E_{j+1,n}^{\ell-3}\Bigl(r^2\mathcal
E_{j,n'}\mathcal
E_{j+1,n-1}{-}\mathcal E_{j+1,n-1}\mathcal E_{j,n'}\Bigr)\Bigr] \ \text{(by Lemma 3.3 (5))}\\
&=(rs)^{-1}\Delta_j^{(n-1)}\Bigl[\mathcal E_{j+1,n}^{\ell-2}\mathcal
E_{j,n}{+}s^4\mathcal E_{j,n'}\mathcal E_{j+1,n-1}\mathcal
E_{j+1,n}^{\ell-3}\Bigr] \ \;\text{(by Lemma 3.5 (1))}\\
\end{split}
\end{equation*}
\begin{equation*}
\begin{split}
&=\Delta_j^{(n-1)}\mathcal E_{j,n'}\mathcal
E_{j+1,n-1}\Bigl(*_1\mathcal E_{j+1,n-1}^{\ell-3}+*_2\mathcal
E_{j+1,n}^{\ell-3}\Bigr)\\
&=\Delta_j^{(n-2)}\Bigl(\mathcal E_{j+1,n-1}^{\ell-1}\mathcal
E_{j,n'}\Bigr)\mathcal E_{j+1,n-1}\bigl(\cdots\bigr)\\
&=*_3\Delta_j^{(n-3)}\Bigl(\mathcal E_{j+1,n-2}^{\ell-1}\mathcal
E_{j,n-1'}\Bigr)\mathcal E_{j+1,n-2}^{\ell-1}\mathcal
E_{j+1,n-1}\bigl(\cdots\bigr)\\
&=\cdots \,\hskip5cm\quad \text{(using (7.3) \& (7.4) repeatedly)}\\
&=*\Delta_j^{(j)}\Bigl(\mathcal E_{j+1,j+1}^{\ell-1}\mathcal
E_{j,j+2'}\Bigr)\Pi_{j+1}^{(n-2)}\mathcal
E_{j+1,n-1}\bigl(\cdots\bigr)
\qquad\ \, \text{(by Lemma 3.10 (6))}\\
&=*\Pi_j^{(j{+}1')}\Bigl(s^{-2}\mathcal
E_{j,j+2'}e_{j+1}+r^{-2}\mathcal
E_{j,j+1'}\Bigr)e_{j+1}^{\ell-2}\,\Pi_{j+1}^{(n-2)}\mathcal
E_{j+1,n-1}\bigl(\cdots\bigr)\\
&=0,\qquad\text{(since $e_{j+1}\Pi_{j+1}^{(n-2)}=0$)}.
\end{split}
\end{equation*}

(ii) \ Write $\Delta_j^{(k')}:=\Pi_j^{(j{+}1')}\Pi_{j+1}^{(k')}$ for
$n'\le k'\le j{+}2'$, and make a convention:
$\Delta_j^{(n{+}1')}:=\Pi_j^{(j{+}1')}$. Assume
$\Delta_j^{(k{+}1')}e_j=0$ for $n{+}1'\le k{+}1'< j{+}2'$, we want
to prove $\Delta_j^{(k')}e_j=0$.

By Lemmas 3.10 (1) \& 3.5 (6) or 3.6 (7), we obtain
\begin{equation*}
\begin{split}
\mathcal{E}_{j+1,n'}^{\ell-1}e_j&=r^2e_j\mathcal{E}_{j+1,n'}^{\ell-1}+
(r^2s)^2\mathcal{E}_{j,n'}\mathcal{E}_{j+1,n'}^{\ell-2},\\
\mathcal{E}_{j+1,k'}^{\ell-1}e_j&=r^2e_j\mathcal{E}_{j+1,k'}^{\ell-1}+
s^2\mathcal{E}_{j,k'}\mathcal{E}_{j+1,k'}^{\ell-2},\qquad j{+}2\leq
k< n.
\end{split}\tag{7.5}
\end{equation*}
By Lemmas 3.6 (6) \& 3.5 (5), $\mathcal{E}_{j,k'}$ quasi-commutes
with $\mathcal E_{j+1,n}^{\ell-1}\mathcal
E_{j+1,n'}^{\ell-1}\cdots\mathcal E_{j+1,k+1'}^{\ell-1}$ (up to a
factor of some power of $(rs)^2$). Again, Lemma 3.7 (7) indicates
that $\mathcal E_{j,k'}$ quasi-commutes with $\mathcal
E_{j+1,k}^{\ell-1}\cdots\mathcal E_{j+1,n-1}^{\ell-1}$ (up to a
factor of some power of $(rs)^{-2}$). Based on these data, together
with (7.3) \& (7.4), using the same proof-procedure as the above
(i): (2), we get
\begin{equation*}
\begin{split}
\Delta_j^{(k')}e_j&=\Delta_j^{(k+1')}\Bigl(\mathcal E_{j+1,k'}^{\ell-1}e_j\Bigr) \qquad\hskip4cm\ \;\text{(by (7.5))}\\
&=*\,\Delta_j^{(k-2)}\Bigl(\mathcal E_{j+1,k-1}^{\ell-1}\mathcal
E_{j,k'}\Bigr)\mathcal
E_{j+1,k}^{\ell-1}\cdots\mathcal E_{j+1,k'}^{\ell-1}\qquad\hskip1cm\text{(by (7.4))}\\
&=\cdots \hskip3cm\, \text{(using the same procedure of (i):
(2))}\\
&=0.
\end{split}
\end{equation*}

The proof is complete.
\end{proof}

\noindent{\it 7.3.} A finite-dimensional Hopf algebra $H$ is
semisimple if and only if $\vn(\int_{H}^l)\neq 0$ or equivalently,
$\vn(\int_{H}^r)\neq 0$. For the algebra $\mathfrak{b}$ above, $y$
gives a basis for $\int_{\mathfrak{b}}^l$ and $y'$ a basis for
$\int_{\mathfrak{b}}^r$. As $\vn(y)=0=\vn(y')$, we have
\begin{prop}\label{7.3}
The Hopf subalgebra $\mathfrak{b}$ is not semisimple.\hfill\qed
\end{prop}

\noindent{\it 7.4.} The second aim of this section is to single out
the distinguished group-like elements of $\mathfrak{b}$ and
$\mathfrak{b}^*$. Since the group-like elements of $\mathfrak{b}^*$
are exactly the algebra homomorphisms in
$\textrm{Alg}_{\mathbb{K}}(\mathfrak{b},\mathbb{K})$, it suffices to
compute its values on the generators.
\begin{prop}\label{6.4}
Write $2\rho=\sum_{j=1}^{n}j(2n-j)\a_j$, where $\rho$ is the half
sum of positive roots of $\mathfrak{so}_{2n+1}$. Let $\gamma \in
\textrm{Alg}_\mathbb{K}(\mathfrak{b},\mathbb{K})$ be defined by
$$
\gamma(e_k)=0, \qquad
\gamma(\omega_k)=\langle\omega_{2\rho}',\omega_k\rangle. \leqno
(7.6)
$$
Then $\gamma$ is the distinguished group-like element of
$\mathfrak{b}^*$.
\end{prop}
\begin{proof}
It suffices to argue that $\gamma$ as in (7.6) satisfies
$ya=\gamma(a)y$ for $a=e_k$ and $a=\omega_k$, $1\leq k\leq n$, and
for $y=tx$ given in Theorem 7.1. Recall from the proof of Theorem
7.2 that $xe_k=0$. Thus, $ye_k=txe_k=0=\gamma(e_k)y$. We have
$$
y\omega_k
=t\Bigl(\prod\limits_{i=1}^{n}\mathcal{E}_{\widehat{i,i+1}}\Bigr)\omega_k=
\Pi_{i=1}^n\langle \omega_i',\omega_k\rangle^{i(2n-i)}
\,t\omega_kx=\gamma(\omega_k)y.$$

This completes the proof.
\end{proof}

Under the assumptions of Lemma 6.1,
$(\mathfrak{b}')^{coop}\cong\mathfrak{b}^*$ as Hopf algebras, via
the map $\phi:(\mathfrak{b}')^{coop}\rightarrow \mathfrak{b}^*,
\phi(\omega'_j)=\gamma_j, \phi(f_j)=\eta_j$ (for definition, see the
proof of Lemma 6.1.). This allows us to define a Hopf pairing
$\mathfrak{b}'\times\mathfrak{b}\rightarrow \mathbb{K}$ whose values
on generators are given by
$$
(f_j\,|\, e_i)=\delta_{ij},\qquad (\omega'_j\,|\, \omega_i)=\langle
\om'_j,\om_i\rangle, \leqno(7.7)
$$
and are zero on all other pairs of generators.

\smallskip
By (7.6) \& (7.7), we get $(\omega'_{2\rho}\,|\,b)=\gamma(b)$, for
all $b \in \mathfrak{b}$.

\smallskip
Note that $\mathfrak b_{s^{-1},r^{-1}}\cong (\mathfrak
b')^{\textit{coop}}\cong \mathfrak b^*$ as Hopf algebras. Under the
isomorphism $\phi\psi^{-1}$ (where $\psi(f_i)=e_i$,
$\psi(\omega_i')=\omega_i$), a nonzero left (resp., right) integral
of $\mathfrak{b}$ maps to a nonzero left (resp., right) integral of
$\mathfrak{b}^*$. Thus, we have

\begin{prop}\label{6.5}
Let $\lambda =\nu\eta$ and $\lambda'=\eta\nu \in \mathfrak{b}^*$,
where $$
\nu=\prod_{i=1}^{n}(1+\gamma_i+\cdots+\gamma_i^{\ell-1}),\quad
\eta=\prod\limits_{i=1}^{n}\eta_{\widehat{i,i+1}},$$  where
$\eta_{\widehat{i,i+1}}=\eta_{i,i}^{\ell-1}\eta_{i,i+1}^{\ell-1}
\cdots\eta_{i,n}^{\ell-1} \eta_{i,n'}^{\ell-1}\cdots
\eta_{i,i+1'}^{\ell-1}$ $(i<n)$,
$\eta_{\widehat{n,n+1}}=\eta_{n,n}^{\ell-1}$, $\eta_{i,i}=\eta_i$,
$\eta_{i,j}=[\eta_{i+1,j},\eta_i]_{s^{2}}$,
$\eta_{i,n'}=[\eta_{n},\eta_{i,n}]_{rs}$,
$\eta_{i,j'}=[\eta_{j},\eta_{i,j+1'}]_{r^{-2}}$.

\medskip Then $\lambda$ is a left integral and $\lambda'$ is a right
integral of $\mathfrak{b}^*$. \hfill\qed
\end{prop}

\begin{coro}
The element $g=\omega_{2\rho}^{-1}$ is the distinguished group-like
element of $\mathfrak{b}$, and under the Hopf pairing in $(7.7)$,
$$
(\omega'_i\,|\, g)= \langle
\omega_i',\omega_{2\rho}^{-1}\rangle=\gamma_i(g).
$$
\end{coro}
\begin{proof} Let
$\mathcal{F}=\prod\limits_{i=1}^{n}\mathcal{F}_{\widehat{i,i+1}}$,
where
$\mathcal{F}_{\widehat{i,i+1}}=\mathcal{F}_{i,i}^{\ell-1}\mathcal{F}_{i,i+1}^{\ell-1}
\cdots\mathcal{F}_{i,n}^{\ell-1} \mathcal{F}_{i,n'}^{\ell-1}\cdots
\mathcal{F}_{i,i+1'}^{\ell-1}$ for $i \leq n{-}1$, and
$\mathcal{F}_{\widehat{n,n+1}}=\mathcal{F}_{n,n}^{\ell-1}$. Then we
have
$$
\omega'_k \mathcal{F}=\langle
\omega_k',\omega_{2\rho}\rangle^{-1}\mathcal{F}\omega'_k.
$$
Because $ \phi^{-1}(\lambda')
=\mathcal{F}\Bigl(\prod_{i=1}^{n}(1+\omega'_i+\cdots+(\omega'_i)^{\ell-1})\Bigr)$,
it follows that
$$
\gamma_k \lambda'=\langle
\omega_k',\omega_{2\rho}\rangle^{-1}\lambda', \ \textrm{ and }\
\eta_k \lambda'=0.
$$
Taking $g=\omega_{2\rho}^{-1}$, we have $\xi\lambda'=\xi(g)\lambda'$
for all $\xi \in \mathfrak{b}^*$.
\end{proof}

\section{To be a ribbon Hopf algebra}
\noindent {\it 8.1.} A finite-dimensional Hopf algebra $H$ is
\textit{quasitriangular} if there is an invertible element $R=\sum
x_i\ot y_i$ in $H\ot H$ such that $\Delta^{op}(a)=R\Delta(a)R^{-1}$
for all $a \in H$, and $R$ satisfies the relations $(\Delta\ot
id)R=R_{1,3}R_{2,3}$, $(id\ot\Delta )R=R_{1,3}R_{1,2}$, where
$R_{1,2}=\sum x_i\ot y_i\ot 1$, $R_{1,3}=\sum x_i\ot 1\ot y_i$, and
$R_{2,3}=\sum 1\ot x_i\ot y_i$.

Write $u=\sum S(y_i)x_i$. Then $c=uS(u)$ is central in $H$ (cf.
\cite{Ka}) and is referred to as the \textit{Casimir element}.

An element $v\in H$ is a \textit{quasi-ribbon element} of
quasitriangular Hopf algebra $(H, R)$ if $\textrm{(i)}$\ $v^2=c$, \
$\textrm{(ii)}$\ $S(v)=v$, \ $\textrm{(iii)}$\ $\vn(v)=1$,\
$\textrm{(iv)}$\ $\Delta(v)=(R_{2,1}R)^{-1}(v\ot v)$, where
$R_{2,1}=\sum y_i\ot x_i$. If moreover $v$ is central in $H$, then
$v$ is a \textit{ribbon element}, and $(H, R, v)$ is called a
\textit{ribbon Hopf algebra}. Ribbon elements provide a very
effective means of constructing invariants of knots or links. 

\medskip \noindent{\it 8.2.} The Drinfel'd double $D(A)$ of a
finite-dimensional Hopf algebra $A$ is quasitriangular, and
Kauffman-Radford (\cite{KR}) proved a criterion for $D(A)$ to be
ribbon.

\begin{theorem} {\rm([KR, Thm. 3])}
Assume $A$ is a finite-dimensional Hopf algebra, let $g$ and
$\gamma$ be the distinguished group-like elements of $A$ and $A^*$
respectively. Then

$(\mathrm{i})$ \ $(D(A),R)$ has a quasi-ribbon element if and only
if there exist group-like elements $h\in A$, $\delta \in A^*$ such
that $h^2=g$, $\delta^2=\gamma$.

$(\mathrm{ii})$ \ $(D(A),R)$ has a ribbon element if and only if
there exist  $h$  and $\delta$  as in $(\mathrm{i})$ such that
$S^2(a)=h(\delta\rightharpoonup a\leftharpoonup \delta^{-1})h^{-1},
\ \forall\; a\in A$.
\end{theorem}

\medskip
\noindent{\it 8.3.} As a consequence of Theorem 8.1, we have
\begin{theorem}\label{7.2}
Assume $r,\,s$ are $\ell$th roots of unity. 
Then for the Hopf subalgebra $\mathfrak{b}$ of
$\mathfrak{u}_{r,s}(\mathfrak{so}_{2n+1})$, the following are
equivalent

$\mathrm{(i)}$ \; $D(\mathfrak{b})$ has a quasi-ribbon element;

$\mathrm{(ii)}$ \ $D(\mathfrak{b})$ has a ribbon element;

$\mathrm{(iii)}$ $\, \ell$ is odd.
\end{theorem}
\begin{proof}
By Corollary 7.6, $g=\omega_{2\rho}^{-1}$ is the distinguished
group-like element of $\mathfrak{b}$. There exists a group-like
element $h=\prod_{j=1}^n\omega_j^{a_j}\in \mathfrak{b}$ with
$a_j\in\mathbb Z$ such that $h^2=g$ if and only if the equations
$2a_j\equiv-j(2n-j)\;\mod\, \ell$ can be solved for $j=1,\cdots, n$.
Solutions exist if and only if $\ell$ is odd. Suppose now $\ell$ is
odd, then we can take $h=\omega_\rho^{-1}\in \mathfrak{b}$. Because
$\gamma=(\omega_{2\rho}'\,|\,\cdot)$ corresponds to
$\omega'_{2\rho}$ under the isomorphism $\phi^{-1}:
\mathfrak{b}^*\rightarrow (\mathfrak{b}')^{coop}$, there exists a
$\delta=(\omega_\rho'\,|\,\cdot) \in \mathfrak{b}^*$ such that
$\delta^2=\gamma$, where $\delta \in
\textrm{Alg}_{\mathbb{K}}(\mathfrak{b},\mathbb{K})$ is defined by
$$
\delta(e_k)=0,\qquad \delta(\omega_k)=\langle
\omega_\rho',\omega_k\rangle.
$$
Then using (7.1), we get
\begin{gather*}
\begin{split}h(\delta\rightharpoonup\omega_k\leftharpoonup\delta^{-1})h^{-1}&=\delta(\omega_k)
\delta^{-1}(\omega_k)h\omega_kh^{-1}=\omega_k=S^2(\omega_k).\\
h(\delta\rightharpoonup
e_k\leftharpoonup\delta^{-1})h^{-1}&=\delta(1)
\delta^{-1}(\omega_k)he_kh^{-1}\\
&=\langle \omega_\rho',\omega_k\rangle^{-1} 
he_kh^{-1}\\
&=\langle \omega_\rho',\omega_k\rangle^{-1}\langle
\omega_k',\omega_\rho\rangle^{-1} e_k\\&=
r_k^{-1}s_ke_k=\omega_k^{-1}e_k\omega_k=S^2(e_k). 
\end{split}
\end{gather*}

By Theorem 8.1, we complete the proof.
\end{proof}

By Theorem \ref{5.2}, we have
\begin{coro}
Assume that $r=\theta^y$, $s=\theta^z$, where $\theta$ is a
primitive $\ell$th root of unity and
$(2^{n-1}(y^n+(-1)^nz^n),\ell)=1$. Then the following are equivalent

$\mathrm{(i)}$ \;  $\mathfrak{u}_{r,s}(\mathfrak{so}_{2n+1})$ has a
quasi-ribbon element;

$\mathrm{(ii)}$ \  $\mathfrak{u}_{r,s}(\mathfrak{so}_{2n+1})$ has a
ribbon element;

$\mathrm{(iii)}$ $\,\ell$ is odd.\hfill\qed
\end{coro}

\bigskip
\bibliographystyle{amsalpha}

\begin{thebibliography}{A}
\medskip

\bibitem[B]{B} J. Beck, \textit{Convex bases of PBW type for quantum affine algebras}, Comm. Math. Phys. \textbf{165} (1994), 193--199.

\bibitem [BGH1]{BGH1} N. Bergeron, Y. Gao and N. Hu, \textit{Drinfel'd doubles and
Lusztig's symmetries of two-parameter quantum groups},
(arXiv:math/0505614), J. Algebra \textbf{301} (2006), 378--405.

\bibitem [BGH2]{BGH2} N. Bergeron, Y. Gao and N. Hu, \textit{Representations of
two-parameter quantum orthogonal groups and symplectic groups},
(arXiv math. QA/0510124), AMS/IP, Studies in Advanced Mathematics,
vol. \textbf{39}, pp. 1--21, 2007.

\bibitem [BH]{BH}X. Bai and N. Hu, \textit{Two-parameter quantum group of exceptional type
$E$-series and convex \textrm{PBW} type basis},
(arXiv.Math.QA/0605179), Algebra Colloq.,
 \textbf{15} (4) (2008), 619--636.



\bibitem [BW1]{BW1} G. Benkart and S. Witherspoon, \textit{Two-parameter quantum groups (of type
$A$) and Drinfel'd doubles}, Algebra Represent. Theory, \textbf{7}
(2004), 261--286.

\bibitem [BW2]{BW2} G. Benkart and S. Witherspoon, \textit{Representations of two-parameter quantum
groups (of type $A$) and Schur-Weyl duality}, in ``Hopf Algebras",
Lecture Notes in Pure and Appl. Math., \textbf{237}, pp. 65--92,
Dekker, New York, 2004.

\bibitem [BW3]{BW3} G. Benkart and S. Witherspoon, \textit{Restricted two-parameter quantum groups (of type
$A$)}, Fields Institute Communications, ``Representations of Finite
Dimensional Algebras and Related Topics in Lie Theory and
Geometry",vol. \textbf{40}, pp. 293--318, Amer. Math. Soc.,
Providence, RI, 2004.

\bibitem[CX]{CX} V. Chari and N. Xi, \textit{Monomial bases of
quantized enveloping algebras}, in ``Recent developments in quantum
affine algebras and related topics",  Contemp. Math. vol.
\textbf{216}, pp. 23--57, Amer. Math. Soc., 2001.



\bibitem [H]{H} N. Hu, \textit{Lyndon words, convex PBW bases and their $R$-matrices for the
two-parameter quantum groups of $B_2$, $C_2$, $D_4$ types},
manuscript 2005.

\bibitem [HP]{HP} N. Hu and Y. Pei, \textit{Notes on two-parameter quantum groups, (I)}, (arXiv.math.QA/0702298),
Sci. in China, Ser. A. \textbf{51} (6) (2008), 1101--1110.

\bibitem [HRZ]{HRZ} N. Hu, M. Rosso and H. Zhang, \textit{Two-parameter affine quantum group
$U_{r,s}(\widehat{\frak{sl}_n})$, Drinfel'd realization and quantum
affine Lyndon basis}, (arXiv:0812.3107), Comm. Math. Phys.
\textbf{278} (2) (2008), 453--486.

\bibitem [HS]{HS} N. Hu and Q. Shi, \textit{The two-parameter quantum group of exceptional type $G_2$
and Lusztig symmetries}, (arXiv:math/0601444), Pacific J. Math.,
\textbf{230} (2) (2007), 327--346.

\bibitem[HW]{HW} N. Hu and X. Wang, \textit{Convex PBW-type Lyndon basis
and restricted two-parameter quantum groups of type $G_2$},
(arXiv:0811.0209), Pacific J. Math. \textbf{241} (2) (2009),
243--273.

\bibitem[HZ]{HZ} N. Hu and H. Zhang, \textit{Vertex representations of two-parameter quantum affine algebras
$U_{r,s}(\widehat{\frak{g}})$: the simply-laced cases}, Preprint
2006-2007.

\bibitem [Ka]{Ka} C. Kassel, \textit{Quantum Groups}, GTM
\textbf{155}, Springer-Verlag Berlin/Heidelberg/New York, 1995.

\bibitem [KR] {KR} L.H. Kauffman and D.E. Radford, \textit{A necessary and
sufficient condition for a finite-dimensional Drinfel'd double to be
a ribbon Hopf algebra}, J. Algebra \textbf{159} (1993), 98--114.


\bibitem [K]{K} V.K. Kharchenko, \textit{A combinatorial approach to
the quantification of Lie algebras}, Pacific J. Math.  \textbf{23}
(2002), 191--233.




\bibitem[LR]{LR} M. Lalonde and A. Ram, \textit{Standard Lyndon bases of Lie algebras and enveloping algebras}, Trans. Amer. Math. Soc.,
\textbf{347} (5) (1995), 1821--1830.

\bibitem [Lu]{Lu} G. Lusztig, \textit{Introduction to Quantum Groups}, Birkh\"auser Boston, 1993.

\bibitem [M]{M} S. Montgomery, \textit{Hopf Algebras and Their Actions on Rings}, CBMS Conf.
Math. Publ., \textbf{82}, Amer. Math. Soc., Providence, 1993.



\bibitem[Ro]{Ro}M. Rosso,  \textit{Lyndon bases and the multiplicative formula for $R$-matrices},
   (2002),  preprint.



\end{thebibliography}

\newpage
$$\text{\Large\bf Appendix: the proof of Lemma 3.6}
$$

\begin{proof}
 (1): The proof is divided into 3 steps:

$(\textrm{i})$ If $j>i+3$, then taking decomposition $\mathcal
E_{i,j'}=[\,\mathcal E_{i,i+2},\mathcal E_{i+3,j'}]_{\bullet}$,
combining with $\mathcal E_{i,i+2}e_{i+1}=e_{i+1}\mathcal E_{i,i+2}$
(by Lemma 3.3 (1)), $\mathcal E_{i+3,j'}e_{i+1}=e_{i+1}\mathcal
E_{i+3,j'}$ (by Lemma 3.1 (2)), we get $\mathcal
E_{i,j'}e_{i+1}=e_{i+1}\mathcal E_{i,j'}$.

$(\textrm{ii})$ If $j=i+3=n$, then $\mathcal
E_{n-3,n'}e_{n-2}=e_{n-2}\mathcal E_{n-3,n'}$ (by Lemma 3.3 (6)). If
$j=i+3<n$, then the decomposition $\mathcal E_{i,i+3'}=[\,\mathcal
E_{i,i+4'},e_{i+3}]_{\bullet}$ shows $\mathcal E_{i,i+3'}e_{i+1}$
$e_{i+1}\mathcal E_{i,i+3'}$, as $e_{i+1}$ commutes with $e_{i+3}$
and $\mathcal E_{i,i+4'}=[\,\mathcal E_{i,i+2},\mathcal
E_{i+3,i+4'}]_{\bullet}$ (by Lemma 3.1 (4), Lemma 3.3 (1) \& Lemma
3.1 (2)).

$(\textrm{iii})$ If $j=i+1=n$ (as $j\ne i+2$ by assumption), then
the case is just (2.7). If $j=i+1=n-1$, then by (2.4), Lemma 3.3 (1)
\& $(B5)$, we have
\begin{equation*}
\begin{split}
&[\,\mathcal E_{n-2,n-1'},e_{n-1}]=\mathcal
E_{n-2,n-1'}e_{n-1}{-}r^{-2}e_{n-1}\mathcal E_{n-2,n-1'}\\
&\quad=(\mathcal E_{n-2,n'}e_{n-1}{-}s^{-2}e_{n-1}\mathcal
E_{n-2,n'})e_{n-1}{-}r^{-2}e_{n-1}(\mathcal
E_{n-2,n'}e_{n-1}{-}s^{-2}e_{n-1}\mathcal E_{n-2,n'})\\
&\quad=(\mathcal E_{n-2,n}e_n{-}rse_n\mathcal E_{n-2,n})e_{n-1}^2
{-}s^{-2}e_{n-1}(\mathcal E_{n-2,n}e_n{-}rse_n\mathcal E_{n-2,n})e_{n-1}\\
&\qquad-r^{-2}e_{n-1}(\mathcal E_{n-2,n}e_n{-}rse_n\mathcal
E_{n-2,n})e_{n-1} {+}(rs)^{-2}e_{n-1}^2(\mathcal
E_{n-2,n}e_n{-}rse_n\mathcal E_{n-2,n})\\
&\quad=\mathcal
E_{n-2,n}\Bigl(e_ne_{n-1}^2-(r^{-2}{+}s^{-2})e_{n-1}e_ne_{n-1}{+}(rs)^{-2}e_{n-1}^2e_n\Bigr)\\
&\qquad{-}(rs)\Bigl(e_ne_{n-1}^2-(r^{-2}{+}s^{-2})e_{n-1}e_ne_{n-1}{+}(rs)^{-2}e_{n-1}^2e_n\Bigr)\mathcal
E_{n-2,n}\\
&\quad=0.
\end{split}
\end{equation*}
If $j=i+1<n-1$, then by (ii), we have
\begin{equation*}
\begin{split}
&[\,\mathcal E_{i,i+1'},e_{i+1}]=\mathcal
E_{i,i+1'}e_{i+1}{-}r^{-2}e_{i+1}\mathcal E_{i,i+1'}\\
&\quad=(\mathcal E_{i,i+2'}e_{i+1}{-}s^{-2}e_{i+1}\mathcal
E_{i,i+2'})e_{i+1}{-}r^{-2}e_{i+1}(\mathcal
E_{i,i+2'}e_{i+1}{-}s^{-2}e_{i+1}\mathcal E_{i,i+2'})\\
&\quad=(\mathcal E_{i,i+3'}e_{i+2}{-}s^{-2}e_{i+2}\mathcal
E_{i,i+3'})e_{i+1}^2{-}s^{-2}e_{i+1}(\mathcal
E_{i,i+3'}e_{i+2}{-}s^{-2}e_{i+2}\mathcal E_{i,i+3'})e_{i+1}\\
&\ {-}r^{-2}e_{i+1}(\mathcal
E_{i,i+3'}e_{i+2}{-}s^{-2}e_{i+2}\mathcal
E_{i,i+3'})e_{i+1}{+}(rs)^{-2}e_{i+1}^2(\mathcal
E_{i,i+3'}e_{i+2}{-}s^{-2}e_{i+2}\mathcal E_{i,i+3'})\\
&\quad=\mathcal
E_{i,i+3'}\Bigl(e_{i+2}e_{i+1}^2{-}(r^{-2}{+}s^{-2})e_{i+1}e_{i+2}e_{i+1}{+}(rs)^{-2}e_{i+1}^2e_{i+2}\Bigr)\\
&\qquad
{-}s^{-2}\Bigl(e_{i+2}e_{i+1}^2{-}(r^{-2}{+}s^{-2})e_{i+1}e_{i+2}e_{i+1}{+}(rs)^{-2}e_{i+1}^2e_{i+2}\Bigr)\mathcal
E_{i,i+3'}\\
&\quad=0,
\end{split}
\end{equation*}
where in the $1$st $``="$ we used the definition, in the $2$nd and
$3$rd $``="$ we used (2.4), in the $4$th $``="$ we used
$e_{i+1}\mathcal E_{i,i+3'}=\mathcal E_{i,i+3'}e_{i+1}$ in (ii),
while in the last $``="$, we used the $(r,s)$-Serre relations
$(B5)$.

\smallskip
(2): From the proof of step (iii) in (1), we have $\mathcal
E_{j-1,j'}e_j=r^{-2}e_j\mathcal E_{j-1,j'}$ for $j<n$. Hence, for
$i<j-1$, we have $\mathcal E_{i,j'}e_j=r^{-2}e_j\mathcal E_{i,j'}$
since $\mathcal E_{i,j'}=[\,\mathcal E_{i,j-2},\mathcal
E_{j-1,j'}]_{r^2}$.

\smallskip
(3): If $i+2=n$, then this is the special case of Lemma 3.4 (4).

If $i+2<n$, then using (2), we obtain
\begin{equation*}
\begin{split}
&[\,\mathcal E_{i,i+1'},e_{i+2}]=\mathcal
E_{i,i+1'}e_{i+2}{-}e_{i+2}\mathcal E_{i,i+1'}\\
&\quad=\Bigl(\mathcal E_{i,i+2'}e_{i+1}{-}s^{-2}e_{i+1}\mathcal
E_{i,i+2'}\Bigr)e_{i+2}{-}e_{i+2}\Bigl(\mathcal
E_{i,i+2'}e_{i+1}{-}s^{-2}e_{i+1}\mathcal E_{i,i+2'}\Bigr)\\
&\quad=\mathcal
E_{i,i+2'}\Bigl(e_{i+1}e_{i+2}{-}r^2e_{i+2}e_{i+1}\Bigr)
{-}(rs)^{-2}\Bigl(e_{i+1}e_{i+2}{-}r^2e_{i+2}e_{i+1}\Bigr)\mathcal
E_{i,i+2'}
\\
&\quad=\mathcal E_{i,i+2'}\mathcal E_{i+1,i+2}{-}(rs)^{-2}\mathcal
E_{i+1,i+2}\mathcal E_{i,i+2'},
\end{split}
\end{equation*}
so (3) is equivalent to prove $$\mathcal E_{i+1,i+2}\mathcal
E_{i,i+2'}=(rs)^{2}\mathcal E_{i,i+2'}\mathcal E_{i+1,i+2}, \qquad
i+2<n.$$

If $i+2=n-1$, then using Lemmas 3.4 (1) and 3.2 (3), we obtain
\begin{equation*}
\begin{split}
\mathcal E_{n-2,n-1}\mathcal E_{n-3,n-1'}&=\mathcal
E_{n-2,n-1}\mathcal E_{n-3,n'}e_{n-1}-s^{-2}\mathcal
E_{n-2,n-1}e_{n-1}\mathcal E_{n-3,n'}\\
&=\Bigl(r^2\mathcal E_{n-3,n'}\mathcal E_{n-2,n-1}+\mathcal
E_{n-3,n}\mathcal E_{n-2,n}-rs\mathcal E_{n-2,n}\mathcal
E_{n-3,n}\Bigr)e_{n-1}\\& \quad -
e_{n-1}\mathcal E_{n-2,n-1}\mathcal E_{n-3,n'}\\
&=(rs)^2\mathcal
E_{n-3,n'}e_{n-1}\mathcal E_{n-2,n-1}\\
&\quad +\underline{\mathcal E_{n-3,n}\mathcal
E_{n-2,n}e_{n-1}-rs\mathcal E_{n-2,n}\mathcal
E_{n-3,n}e_{n-1}}\\
&\quad -r^2e_{n-1}\mathcal E_{n-3,n'}\mathcal
E_{n-2,n-1}\\
&\quad -\underline{e_{n-1}\mathcal E_{n-3,n}\mathcal
E_{n-2,n}+rse_{n-1}\mathcal E_{n-2,n}\mathcal E_{n-3,n}}\\
&= (rs)^2\mathcal E_{n-3,n-1'}\mathcal E_{n-2,n-1},
\end{split}
\end{equation*}
where the sum of terms underlined is $0$ since $e_{n-1}$ commutes
with $\mathcal E_{n-2,n}, \mathcal E_{n-3,n}$. The argument for
$i+2<n-1$ is similar by using Lemmas 3.3 (8) and 3.4 (2).

\smallskip
(4) For $i<j<l< n$: first noting the case when $j=l{-}1$ and
$i=j{-}1=l{-}2$, by (3), we have
$\mathcal{E}_{l-2,l-1'}e_l=e_l\mathcal{E}_{l-2,l-1'}$ for $l< n$.
For the case when $j=l{-}1$ and $i\le l-3$, using $\mathcal
E_{i,l-1'}=[\,\mathcal E_{i,l-3},\mathcal E_{l-2,l-1'}]_{r^2}$ and
$e_l\mathcal E_{i,l-3}=\mathcal E_{i,l-3}e_l$, we obtain
$e_l\mathcal E_{i,l-1'}=\mathcal E_{i,l-1'}e_l$.

For the case when $j<l-1$, noting that $\mathcal
E_{i,j'}=[\cdots[\,\mathcal
E_{i,l-1'},e_{l-2}]_{s^{-2}},\cdots,e_j]_{s^{-2}}$ and $e_l
e_k=e_ke_l$ for $j\le k\le l-2$, we get
$e_l\mathcal{E}_{i,j'}=\mathcal{E}_{i,j'}e_l$ for $i<j<l< n$.

\smallskip
(5): First, we claim: $\mathcal E_{i,n'}\mathcal
E_{i,n-1'}=r^{-2}\mathcal E_{i,n-1'}\mathcal E_{i,n'}$ for $i<n-1$.
Indeed, by Lemma 3.4 (1), $e_{n-1}\mathcal E_{i,n'}=r^2\mathcal
E_{i,n'}e_{n-1}+\mathcal E_{i,n}\mathcal E_{n-1,n}-rs\mathcal
E_{n-1,n}\mathcal E_{i,n}$. Again by Lemma 3.5 (3) \& Lemma 3.4 (5),
$\mathcal E_{i,n}\mathcal E_{i,n'}=s^2\mathcal E_{i,n'}\mathcal
E_{i,n}$, $\mathcal E_{i,n'}\mathcal E_{n-1,n}=\mathcal
E_{n-1,n}\mathcal E_{i,n'}$. So using (2.4), we have
\begin{equation*}
\begin{split}
&[\,\mathcal E_{i,n'},\mathcal E_{i,n-1'}]=
 \mathcal E_{i,n'}\mathcal E_{i,n-1'}{-}r^{-2}\mathcal
E_{i,n-1'}\mathcal E_{i,n'}\\
&\quad=\mathcal E_{i,n'}^2e_{n-1}{-}(s^{-2}{+}r^{-2})\mathcal
E_{i,n'}e_{n-1}\mathcal E_{i,n'}{+}(rs)^{-2}(e_{n-1}\mathcal
E_{i,n'})\mathcal E_{i,n'}\\
&\quad=\mathcal E_{i,n'}^2e_{n-1}{-}r^{-2}\mathcal
E_{i,n'}(e_{n-1}\mathcal E_{i,n'}){+}(rs)^{-2}(\mathcal
E_{i,n}\mathcal E_{n-1,n}{-}rs\mathcal E_{n-1,n}\mathcal
E_{i,n})\mathcal E_{i,n'}\\
&\quad=-r^{-2}\mathcal E_{i,n'}\Bigl(\mathcal E_{i,n}\mathcal
E_{n-1,n}{-}rs\mathcal E_{n-1,n}\mathcal
E_{i,n}\Bigr){+}(rs)^{-2}\Bigl(\mathcal E_{i,n}\mathcal
E_{n-1,n}{-}rs\mathcal E_{n-1,n}\mathcal
E_{i,n}\Bigr)\mathcal E_{i,n'}\\
&\quad=0,
\end{split}
\end{equation*}
this means $\mathcal E_{i,n'}\mathcal E_{i,n-1'}=r^{-2}\mathcal
E_{i,n-1'}\mathcal E_{i,n'}$ for $i<n-1$.

Next, we claim: $\mathcal E_{i,n'}\mathcal E_{i,j'}=r^{-2}\mathcal
E_{i,j'}\mathcal E_{i,n'}$ for $i<j<n$. In fact, noting $\mathcal
E_{i,j'}=[\,\cdots[\,\mathcal
E_{i,n-1'},e_{n-2}]_{s^{-2}},\cdots,e_j]_{s^{-2}}$ and $e_k\mathcal
E_{i,n'}=\mathcal E_{i,n'}e_k$ for any $k$ with $i<k<n-1$,  we
arrive at the required result.

Finally, we claim: $\mathcal E_{i,l'}\mathcal
E_{i,j'}=r^{-2}\mathcal E_{i,j'}\mathcal E_{i,l'}$ for $j<l<n$. This
follows from (4) and $\mathcal E_{i,l'}=[\,\cdots[\,\mathcal
E_{i,n'},e_{n-1}]_{s^{-2}},\cdots,e_{l}]_{s^{-2}}$.

\smallskip
(6) follows from (4) and Lemma 3.5 (6).

\smallskip
(7) follows from (2) \& (6).

\smallskip
(8): By Lemma 3.3 (5), we have $\mathcal E_{i,n-1}\mathcal
E_{i,n}=s^2\mathcal E_{i,n}\mathcal E_{i,n-1}$ and
\begin{equation*}
\begin{split}
&\mathcal E_{i,n-1}\mathcal E_{i,n'}{-}(rs)^2\mathcal
E_{i,n'}\mathcal
E_{i,n-1}\\
&\quad=\mathcal E_{i,n-1}(\mathcal E_{i,n}e_n{-}rse_n\mathcal
E_{i,n}){-}(rs)^2(\mathcal E_{i,n}e_n{-}rse_n\mathcal
E_{i,n})\mathcal
E_{i,n-1}\\
&\quad=s^2\mathcal E_{i,n}(\mathcal E_{i,n-1}e_n){-}rs(\mathcal
E_{i,n-1}e_n)\mathcal E_{i,n}{-}(rs)^2\mathcal E_{i,n}(e_n\mathcal
E_{i,n-1}){+}r^3s(e_n\mathcal
E_{i,n-1})\mathcal E_{i,n}\\
&\quad=s(s-r)\mathcal E_{i,n}^2.
\end{split}
\end{equation*}

This completes the proof.
\end{proof}

\end{document}